\documentclass[11pt]{article}
\usepackage{amsfonts,latexsym}
\title{Logarithmic version of the Parshin symbol for surfaces}
\author{Ivan Emilov Horozov}
\date{April 30, 2010}
\newcommand{\beq}{\begin{equation}}
\newcommand{\eeq}{\end{equation}}
\newcommand{\beqa}{\begin{eqnarray}}
\newcommand{\eeqa}{\end{eqnarray}}
\newcommand{\beaa}{\begin{eqnarray*}}
\newcommand{\ben}{\begin{eqnarray*}}
\newcommand{\eaa}{\end{eqnarray*}}
\newcommand{\een}{\end{eqnarray*}}

\newcommand{\text}{\textrm}
\newcommand \nc {\newcommand}
\nc \proof {\noindent {\em{Proof.\/ }}}
\nc \qed {$\Box$\hfill}
\newtheorem{theorem}{Theorem}[section]
\newtheorem{lemma}[theorem]{Lemma}
\newtheorem{proposition}[theorem]{Proposition}
\newtheorem{corollary}[theorem]{Corollary}
\newtheorem{definition}[theorem]{Definition}
\newtheorem{example}[theorem]{Example}
\newtheorem{remark}[theorem]{Remark}
\newtheorem{conjecture}[theorem]{Conjecture}
\newtheorem{question}[theorem]{Question}
\nc \bth[1] {\begin{theorem}\label{t#1} }
\nc \ble[1] {\begin{lemma}\label{l#1} }
\nc \bpr[1] {\begin{proposition}\label{p#1} }
\nc \bco[1] {\begin{corollary}\label{c#1} }
\nc \bde[1] {\begin{definition}\label{d#1}\rm }
\nc \bex[1] {\begin{example}\label{e#1}\rm }
\nc \bre[1] {\begin{remark}\label{r#1}\rm }
\nc \bcon[1] {\begin{conjecture}\label{con#1}\rm }
\nc \bque[1] {\begin{question}\label{que#1}\rm }
\nc {\eth} { \end{theorem} }
\nc {\ele} { \end{lemma} }
\nc {\epr}{ \end{proposition} }
\nc {\eco} { \end{corollary} }
\nc {\ede} {\end{definition} }
\nc {\eex} { \end{example} }
\nc {\ere} {\end{remark} }
\nc {\econ} { \end{conjecture} }
\nc {\eque} {\end{question} }
\nc \eqref[1] {{\rm{(\ref{#1})}}}
\nc \thref[1]{Theorem \ref{t#1}}
\nc \leref[1]{Lemma \ref{l#1}}
\nc \prref[1]{Proposition
\ref{p#1}} \nc \coref[1]{Corollary \ref{c#1}}
\nc \deref[1]{Definition \ref{d#1}}
\nc \exref[1]{Example \ref{e#1}}
\nc \reref[1]{Remark \ref{r#1}}
\nc \conref[1]{Conjecture\ref{con#1}}
\def\a{\alpha}
\def\b{\beta}

\def\e{\epsilon}

\def \d {{\mathrm d}}

\def \Q {{\mathbb Q}}

\def \C {{\mathbb C}}
\def \mod { {\mathrm{mod}} }

\def \d { {\partial}}

\nc \Wr {Wr} \nc \GRN { \Gr^{(N)} }
\nc \GRA[1] { \Gr_A^{(#1)} }   %% Gr_A
\nc \GRAN { \GRA{N} } \nc \GrA[1] { \Gr_A(#1) }\nc \GrAa {
\GrA{\alpha} }
\nc \GRB[1] { \Gr_B^{(#1)} }   %% Gr_B
\nc \GRBN { \GRB{N} } \nc \GrB[1] { \Gr_B(#1) } \nc \GrBb {
\GrB{\beta} }
\nc \GRMB[1] { \Gr_{MB}^{(#1)} }   %% Gr_{MB}
\nc \GRMBN { \GRMB{N} } \nc \GrMB[1] { \Gr_{MB}(#1) } \nc \GrMBb {
\GrMB{\beta} }
\begin{document}

%%%%%%%%%%%%%%%%%%%%%%%%%%%%%%%%%%%%%%%%%%%%%%%%%%%%%%%%%%%%%%%%%%%%%%%
%%%%%%%%%%%%%%%%%%%%%%%    Title    %%%%%%%%%%%%%%%%%%%%%%%%%%%%%%%%%%%%%%
\title{{\LARGE\bf{Refinement of the Parshin symbol for surfaces}}}

\author{
I. ~Horozov
\thanks{E-mail: ivan.horozov@uni-tuebingen.de}
\\ \hfill\\ \normalsize \textit{Mathematisches Institut,}\\
\normalsize \textit{Universit\"at T\"ubingen, Auf der Morgenstelle 10,}\\
\normalsize \textit {72076 T\"ubingen, Germany }  \\
}
\date{}
\maketitle
\begin{abstract}
On an algebraic curve there are Tate symbols, which satisfy Weil
reciprocity law. The analogues in higher dimensions are the
Parshin symbols, which satisfy Kato-Parshin reciprocity laws. We
give a refinement of the Parshin symbol for surfaces, using
iterated integrals in the sense of Chen. The product of the
refined symbol over the cyclic permutations of the functions
recovers the Parshin symbol. Also, we construct a logarithmic
version of the Parshin symbol. We prove reciprocity laws for both
the refined symbol and a logarithm of the Parshin symbol.
\end{abstract}
%%%%%%%%%%%%%%%%%%%
%%%%%%%%%%%%%%%%%%%%%%%%%%%%%%%%%%%%%%%%%%%%%%%%%%%%%%%%%%%%%%%%%%%%%%%%%%%%%%
%%%%%%%%%%%%%%%%%%%%%   Introduction    %%%%%%%%%%%%%%%%%%%%%%%%%%%%%%%%%%%%%%
\tableofcontents
\setcounter{section}{-1}%%%%%%%%%%%%%%%%%%%%%%%%%%%%%%%%%%%%%%%%%%%%%%%%%%%%%
\section{Introduction}
This paper is the second one from a series of papers on
reciprocity laws on complex varieties. Instead of using a homology
or a cohomology theory, we use properties of certain fundamental
groups. We capture structures of these fundamental groups by
examining iterated integrals.  In this way, we have proven various
reciprocity laws on a Riemann surface in \cite{H}. We use some of
these reciprocity laws here. But more importantly, we develop
further the use of iterated integrals in establishing new
reciprocity laws.

Parshin has considered iterated integrals \cite{P3} at the same
time as Chen \cite{Ch}. However, the ones by Chen are more general
and we use some of his constructions.

Let us recall the Weil reciprocity for the Tate symbol. Let $f_1$
and $f_2$ be two non-zero rational functions on a Riemann surface
$C$. At a point $P$ on $C$,
let $x$ be a rational function on $C$,
which has a zero of order $1$ at $P$.
Let $n_k$ be the (vanishing) order of $f_k$ at $P$.
Define $g_k$ as
$$g_k=x^{-n_k}f_k,$$
for $k=1,2$. Note that $g_k$ is a rational function, depending on the
choices of $P$ and $x$, which has no zero and no pole at the point
$P$. Then the Tate symbol is defined as
$$\{f_1,f_2\}_P=(-1)^{n_1n_2}\left(\frac{f_1^{n_2}}{f_2^{n_1}}\right)(P)
=(-1)^{n_1n_2}\frac{g_1(P)^{n_2}}{g_2(P)^{n_1}}.
$$

The Weil reciprocity for the Tate symbol is
$$\prod_{P\in C}\{f_1,f_2\}_P=1.$$

We are going to use a local coordinate, which is a rational function $x$,
as opposed to a uniformizer in the corresponding local field.

Let us explain what is the relation between a local coordinate and
a uniformizer. By a local coordinate at a point $P$ on a curve $C$
we mean a rational function $x$ on $C$, which has a zero of order
$1$ at $P$. The rational function $x$ might have other zeros or
poles; but this will not play an essential role. Let ${\cal{O}}_C$
be the structure sheaf on $C$. Let $U$ be a Zariski open set in
$C$, where $x$ is defined. Denote also by ${\cal{O}}_{C,P}$ the
local ring of stalks at $P$. Define $m_P$ to be the maximal ideal
in ${\cal{O}}_{C,P}$. And let
$$\hat{\cal{O}}_{C,P}=\lim_{\leftarrow}{\cal{O}}_{C,P}/m_p^N$$
be the completion of the local ring ${\cal{O}}_{C,P}$.
Then
$$x\in {\cal{O}}(U)\subset {\cal{O}}_{C,P}\subset \hat{\cal{O}}_{C,P}.$$
Thus, we can think of the rational function $x$ as a local
uniformizer in the complete local ring $\hat{\cal{O}}_{C,P}$.

There is analogous statement for surfaces. We recall the Parshin
symbol for a surface. Then we give a definition in terms of
rational functions as opposed to uniformizers, which is not as
general as the algebraic definition but it allows us to use
integrals. We need the new definition in order to construct a
refinement of the Parshin symbol together with several logarithmic
symbols.

Let $f_1,f_2,f_3$ be three non-zero rational functions on a smooth
complex surface $X$. Let ${\cal{O}}_X$ be the structure sheaf of
the surface $X$. Let $C$ be a non-singular curve in $X$ and $P$
be a point on $C$. Let ${\cal{I}}_C$ be the sheaf of ideals,
defining the curve $C$. Let $\widehat{({\cal{O}}_X)}_C$ be the
completion of the structure sheaf ${\cal{O}}_X$ with respect to
the sheaf of ideals ${\cal{I}}_C$. Let $x$ be an element of
$\widehat{({\cal{O}}_X)}_C{\cal{I}}_C-\widehat{({\cal{O}}_X)}_C{\cal{I}}_C^2$,
defining the curve $C$ in a Zariski neighborhood of $P$.
Let $m_k=ord_C (f_k)$. Let $\bar{h}_k$ be in
$\widehat{({\cal{O}}_X)}_C/{\cal{I}}_C={\cal{O}}_C$ such that
$$\bar{h}_k \equiv x^{-m_k}f_k\mbox{ } \mod\mbox{ }  {\cal{I}}_C.$$
Let $y$ be a uniformizer at $P$ in the completion
$\widehat{({\cal{O}}_C)}_P$. Let $m_P$ be the maximal ideal in
$\widehat{({\cal{O}}_C)}_P$, defining $P$, and let
$$n_k=ord_P(\bar{h}_k)=ord_P(x^{-m_k} f_k\mbox{ }\mod \mbox{ }{\cal{O}}_C).$$
Define $\bar{g}_k\in \C$ so that
$$\bar{g}_k\equiv y^{-n_k}\bar{h}_k\mbox{ }\mod
\mbox{ } m_P\equiv y^{-n_k}(x^{-m_k}f_k \mbox{ }\mod\mbox{
}{\cal{O}}_C)\mbox{ } \mod\mbox{ } m_p.$$

Then the Parshin symbol (computed  by Fesenko and Vostokov
\cite{FV}) is defined as
$$\{f_1,f_2,f_3\}_{C,P}=(-1)^K g_1^{D_{1}}
g_2^{D_{2}} g_3^{D_{3}},
$$
where
$$D_{1}=\left|\begin{tabular}{ll}
$m_2$ & $n_2$\\
$m_3$ & $n_3$
\end{tabular}\right|,\mbox{ }
D_{2}=\left|\begin{tabular}{ll}
$m_3$ & $n_3$\\
$m_1$ & $n_1$
\end{tabular}\right|,\mbox{ }
D_{3}=\left|\begin{tabular}{ll}
$m_1$ & $n_1$\\
$m_2$ & $n_2$
\end{tabular}\right|$$
and
$$K=n_1n_2m_3+n_2n_3m_1+n_3n_1m_2 -m_1m_2n_3-m_2m_3n_1-m_3m_1n_2.$$

One of the Kato-Parshin reciprocity laws is
$$\prod_{P}\{f_1,f_2,f_3\}_{C,P}=1,$$
where the product is over all points $P$ on $C$.

We need an alternative definition of the Parshin symbol, for which
we can apply analytic techniques. Given three non-zero functions
on $X$, $f_1$, $f_2$ and $f_3$, we want the decomposition
$f_k=x^{i_k}y^{j_k}g_k$, for $k=1,2,3$, so that $g_k$ is well
defined and non-zero at $P$ and $x$ and $y$ are rational functions
on $X$, not just uniformizers. This is not always possible. So we
make assumptions on the divisors of $f_1$, $f_2$ and $f_3$. We
assume that the divisors of the functions  $f_1$, $f_2$ and $f_3$
have normal crossings (that is, at intersection point only two
components meet transversally). Also, we assume that each
component of these divisors is a non-singular curve. These
properties can be achieved for any non-zero functions $f_1$, $f_2$
and $f_3$, after successive blow-ups of the surface $X$.

Now we consider carefully the algebraic definition of the Parshin
symbol. Denote by $C_i$ for $i=0,1,\dots,N$ the components of the
divisors of $f_1$, $f_2$ and $f_3$. Instead of choosing $x\in
\widehat{({\cal{O}}_X)}_{C_0}{\cal{I}}_{C_0}
-\widehat{({\cal{O}}_X)}_{C_0}{\cal{I}}_{C_0}^2$, defining $C$ in
a neighborhood of $P$, we choose a rational function $x$, such
that $ord_{C_0} x=1$ and $x$ defines $C$ in a neighborhood of $P$.
In particular, for such a rational function, we have $x\in
\widehat{({\cal{O}}_X)}_{C_0}{\cal{I}}_{C_0}
-\widehat{({\cal{O}}_X)}_{C_0}{\cal{I}}_{C_0}^2$. Now we use the
divisors of $f_1$, $f_2$ and $f_3$ have smooth normal crossing.
Let $P\in C_0\cap C_j$. Since $x$ defines $C$ in a neighborhood of
$P$, we have that $ord_{C_j}x=0$. Define $m_k=ord_{C_0}f_k$.
Similarly ,we define a rational function $y$ such that
$ord_{C_j}y=1$ and $ord_{C_0}y=0$. Let $n_k=ord_{C_j}f_k$. Let
$h_k=x^{-m_k}f_k$. Then
$$\bar{h}_k\sim h_k\mbox{ }mod\mbox{ }{\cal{I}}_{C_0}.$$
Define $$\bar{y}\sim y \mbox{ }mod\mbox{ }{\cal{I}}_{C_0}$$ on
$C_0$. Let $\tilde{g}_k=\bar{y}^{-n_k}\bar{h}_k$ on $C_0$. Then
similar to the one dimensional case, we have
$$\bar{g}_k\sim \tilde{g}_k\mbox{ }mod\mbox{ }m_P,$$
where $m_P$ is the maximal ideal in the local ring ${\cal{O}}_{C_0,P}$.
Moreover, $\bar{g}_k= \tilde{g}_k(Q)=g_k(Q)$.

Let
$$div(f_k)=\sum_{i=0}^N n_{ki} C_i.$$
Let $x_i$ be a rational function on $X$,
which has zero along $C_i$ of order $1$
(this is the order of vanishing of $x_i$ at the generic point),
and has no zeroes or poles along $C_j$ for $j\neq i$. That is,
$$ord_{C_j}x_i=\delta_{ij}.$$
The rational functions $x_i$ might have other zeros or poles.
However, this poses only a technical consideration and it does not
contribute to the symbols that we consider.

Note that in the algebraic definition, the Parshin symbol
$$\{f_1,f_2,f_3\}_{C,P}=1,$$
if $C$ is not a divisor of any of the functions $f_1,f_2,f_3$, or
if $P$ is not an intersection of two curves from the divisors of
the functions. Thus, it is enough to consider only the case, when
$C$ is a curve from the divisors of the functions $f_1$, $f_2$ and
$f_3$ and $P$ is a point of intersection of two divisors.

Let $C=C_0$ and let $P\in C_0\cap C_j$ for $j>0$. We are going to
define $\{f_1,f_2,f_3\}_{C_0,P},$ using the rational functions
$x_0$ and $x_j$. Let $x=x_0$ and $y=x_j$. For their corresponding
exponents, define $m_k=n_{k0}=ord_{C_0}f_k$ and
$n_k=n_{kj}=ord_{C_j}f_k$. Let $$f_k=x^{m_k}y^{n_k}g_k.$$ Note
that $g_k$ is a rational function, which is well defined and
non-zero at $P$.

The definition in terms of rational function $g_k$, which we are going to use, is
$$\{f_1,f_2,f_3\}_{C_0,P}=(-1)^K g_1(P)^{D_{1}}
g_2(P)^{D_{2}} g_3(P)^{D_{3}},$$
where $D_{1}$, $D_{2}$, $D_{3}$ and $K$
are defined as above in terms of $m_k$ and $n_k$ for $k=1,2,3$.

Now we can define a refinement of the Parshin symbol
$\{f_1,f_2,f_3\}_{C_0,P}$. First, we need one coherence condition:
For each point $P\in C_0\cap C_j$  and for each curve $C_j$ we are
going to use the same rational function $x=x_0$, which has order
$1$ at $C_0$. The refinement of the Parshin symbol will be
invariant with respect to choices of the rest of the rational
functions $x_i$ for $i>1$.
%%%%%%%%%%%%%%%%%%%%%%%%%%%%%%%%%%%%%%%%%%%
\begin{definition} With the above notation, we define
a refinement of the Parshin symbol
$$(f_1,f_2,f_3)^{x_0}_{C_0,P}=(-1)^{n_1n_3m_2-m_1m_3n_2}
\left(\frac{g_1(P)^{n_3}}{g_3(P)^{n_1}}\right)^{m_2}.$$
\end{definition}
%%%%%%%%%%%%%%%%%%%%%%%%%%%%%%%%%%%%%%%%%%%%%%%%%%%%%
The new symbol resembles a power of the Tate symbol for curves.
However, there is a simple relation among the Parshin symbol and
the refinement of the Parshin symbol.
%%%%%%%%%%%%%%%%%%%%%%%%%%%%%%%%%%%%%%%%%%%%%%%%%%
\begin{theorem} (The Parshin symbol in terms of the refinement of the Parshin symbol)
$$\{f_1,f_2,f_3\}_{C_0,P}=\prod_{cycl}(f_1,f_2,f_3)^{x_0}_{C_0,P},$$
where the product is taken over
cyclic permutations of the indexes of $f_1$, $f_2$ and $f_3$.
\end{theorem}
%%%%%%%%%%%%%%%%%%%%%%%%%%%%%%%%%%%%%%%%%%%%%%%%%%%%%%
Moreover, we have a reciprocity law for the refinement of the Parshin symbol.
%%%%%%%%%%%%%%%%%%%%%%%%%%%%%%%%%%%%%%%%%%%%
\begin{theorem} A reciprocity law for the refinement of the Parshin symbol is
$$\prod_{P}(f_1,f_2,f_3)^{x_0}_{C_0,P}=1,$$
where the product is taken over all points $P$ on $C_0$.
\end{theorem}
%%%%%%%%%%%%%%%%%%%%%%%%%%%%%%%%%%%%%%

As we mentioned in the beginning of the introduction, we use
iterated integrals of differential $1$-forms with logarithmic
poles. An example of such differential form is $df_1/f_1$. Its
integral is $log(f_1)$. With this approach, first we obtain
logarithmic symbols with additive reciprocity laws. After
exponentiating, we obtain the above symbols and the corresponding
multiplicative reciprocity laws.

In order to define the logarithmic symbols we need to integrate
over certain paths $\gamma_i$.
%%%%%%%%%%%%%%%%%%%%%%%%%%%%%%%%%%%%%%%%%%%%%%%%%%%%%%%%%%%%%%%%%%
\begin{definition} (Points, loops and paths on $C_0$)
Let $\{P_1,\dots,P_M\}=C_0\cap(\bigcup_j C_j)$ be all intersection
points of $C_0$ with the rest of the curves from the divisors of
$f_1,f_2,f_3$. Let $Q\in C_0$ be a base point, different from
$P_1,\dots,P_M$. Let $R_{i}$ be a point on $C_0$, which is within
an $\e$-neighborhood of the intersection point $P_i$. Let
$\sigma_{i}$ be a simple loop on $C_0$ around $P_i$, based at $Q$.
Let $\gamma_{i}$ be a path from $Q$ to $R_{i}$. Define also
$\sigma^0_i$ to be a (small) simple loop around $P_i$, based at
$R_i$, so that
$$\sigma_{i}=\gamma_i\sigma^0_i\gamma_i^{-1}.$$
Choose the paths $\gamma_{i}$ so that the loop $\delta$,
defined by
$$\delta=\left(\prod_{i=1}^g [\a_i,\b_i]\right)
\left(\prod_{i=1}^N \sigma_{i}\right),$$
is homotopic to the trivial loop at $Q$,
where $g$ is the genus of the curve $C_0$ and
$[\a_i,\b_i]$ is the commutator of loops $\a_i$ and $\b_i$
around the handles of $C_0$.
\end{definition}
%%%%%%%%%%%%%%%%%%%%%%%%%%%%%%%%%%%%%%%%%%%%%%%%%%%%%%%%
At the point $P_i$ define $g_k$ so that
$$f_k=x_0^{m_k}x_j^{n_k}g_k$$
on the surface $X$, using the rational functions $x_0$ and $x_j$
for some $j$ as local coordinates as opposed to uniformizers. Then
$g_k$ is well-defined and non-zero at $P_i$.
%%%%%%%%%%%%%%%%%%%%%%%%%%%%%%%%%%%%%%%%%%%%%%%%%%%
\begin{definition} Logarithm of the refinement of the Parshin symbol is defined by
$$Log(f_1,f_2,f_3)^{x_0,\gamma_{i}}_{C_0,P_{i}}
=(2\pi i)^2\left(\pi i(m_2n_1n_3-n_2m_1m_3)
+m_2n_3\int_{\gamma_{i}}\frac{dg_{1}}{g_{1}}
-m_2n_1\int_{\gamma_{i}}\frac{dg_{3}}{g_{3}},
\right)$$
\end{definition}
%%%%%%%%%%%%%%%%%%%%%%%%%%%%%%%%%%%%%%%%%%%%%
Let $n_{kj}=ord_{C_j}f_k$. Let also $L_j$ be the number of
intersection points of $C_0$ with $C_j$. Define $h_k$ by
$h_k=x_0^{-m_k}f_k$. Define also
$$D_{1}(j)=\left|\begin{tabular}{ll}
$m_2$ & $n_{2j}$\\
$m_3$ & $n_{3j}$
\end{tabular}\right|,\mbox{ }
D_{2}=\left|\begin{tabular}{ll}
$m_3$ & $n_{3j}$\\
$m_1$ & $n_{1j}$
\end{tabular}\right|,\mbox{ }
D_{3}=\left|\begin{tabular}{ll}
$m_1$ & $n_{1j}$\\
$m_2$ & $n_{2j}$
\end{tabular}\right|$$
%%%%%%%%%%%%%%%%%%%%%%%%%%%%%%%%%%%%%%%%%%%
\begin{theorem} A reciprocity law for the logarithm of the refinement of the
Parshin symbol is
$$\sum_{i}Log(f_1,f_2,f_3)^{x_0,\gamma_{i}}_{C_0,P_{i}}=(2\pi i)^3(M+N),$$
where
$$
\begin{tabular}{ll}
$M=$
&$\sum_{j_1<j_2}
(n_{1j_1}D_1(j_2)-n_{3j_1}D_3(j_2))L_{j_1}L_{j_2}$\\
\\
&$+\sum_{j_2=1}^N
(n_{1j_2}D_1(j_2)-n_{3j_2}D_3(j_2))\frac{1}{2}L_{j_2}(L_{j_2}-1)$\\
\\
$N=$&$(2\pi i)^{-2} m_1\sum_{j=1}^g\left(\int_{\a_j}\frac{dh_3}{h_3}\int_{\b_j}\frac{dh_2}{h_2}
-\int_{\b_j}\frac{dh_3}{h_3}\int_{\a_j}\frac{dh_2}{h_2}\right)+$\\
\\
&$ +(2\pi i)^{-2} m_3\sum_{j=1}^g\left(\int_{\a_j}\frac{dh_1}{h_1}\int_{\b_j}\frac{dh_2}{h_2}
-\int_{\b_j}\frac{dh_1}{h_1}\int_{\a_j}\frac{dh_2}{h_2}\right),$
\end{tabular}
$$
\end{theorem}
%%%%%%%%%%%%%%%%%%%%%%%%%%%%%%%%%%%%%%%%%%%%%%%%%%%%%%%%%%
The relation between the refinement of the Parshin symbol and the
logarithm of the refinement of the Parshin symbol is given by the
following theorem.
\begin{theorem}
$$(f_1,f_2,f_3)^{x_0}_{C_0,P_i}=
exp((2\pi i)^{-2}Log(f_1,f_2,f_3)^{x_0,\gamma_i}_{C_0,P_i}
g_1(Q)^{m_2n_3}g_3(Q)^{-m_2n_1},$$
where $Q$ is the base point.
\end{theorem}
%%%%%%%%%%%%%%%%%%%%%%%%%%%%%%%%%%%%%%%%%%%%%%%%%%
\begin{definition} We define logarithm of the Parshin symbol as
$$Log\{f_1,f_2,f_3\}^{\gamma_i}_{C_0,P_i}
=\prod_{cycl}Log(f_1,f_2,f_3)^{x_0,\gamma_i}_{C_0,P_i},$$
where the product is taken over
cyclic permutations of the indexes of $f_1$, $f_2$ and $f_3$.
\end{definition}
%%%%%%%%%%%%%%%%%%%%%%%%%%%%%%%%%%%%%%%%%%%%%%%%%%%%%%
\begin{theorem} The relation to the Parshin symbol is
$$\{f_1,f_2,f_3\}_{C_0,P_i}
=exp((2\pi i)^{-2}Log\{f_1,f_2,f_3\}^{\gamma_i}_{C_0,P_i})
g_1(Q)^{D_1}g_2(Q)^{D_2}g_3(Q)^{D_3},$$
where $Q$ is the base point.
\end{theorem}
%%%%%%%%%%%%%%%%%%%%%%%%%%%%%%%%%%%%%%%%%%%%%%%%%%%%%%%%%%%
A reciprocity law for the logarithm of the Parshin symbol is given
by the following theorem.
\begin{theorem}
$$\sum_i Log\{f_1,f_2,f_3\}^{\gamma_i}_{C_0,P_i}=0.$$
\end{theorem}
%%%%%%%%%%%%%%%%%%%%%%%%%%%%%%%%%%%%%%%%%%%%%%%%%%%%%%%
It is interesting to compare the logarithmic symbols on a complex
surface $Log\{f_1,f_2,f_3\}^{\gamma_i}_{C_0,P_i}$ and
$Log(f_1,f_2,f_3)^{x_0,\gamma_i}_{C_0,P_i}$ with logarithmic
symbol on a complex curve. We call the logarithmic symbol on a
complex curve a logarithm of the Tate symbol.
%%%%%%%%%%%%%%%%%%%%%%%%%%%%%%%%%%%%%%%%%%%%%%%%%%%%%%
\begin{definition} (Logarithm of the Tate symbol)
Let $f_1$ and $f_2$ be two non-zero rational functions on $C$. Let
$P_1,\dots,P_M$ be the points in the divisors of $f_1$ and $f_2$.
Denote by $m_i$ the order of $f_1$ at $P_i$ and by $n_i$ the order
of $f_2$ at $P_i$.

Near $P_i$, let $x$ be a rational function on $C$,which has order $1$ at $P_i$.
Define $g_1$ and $g_2$ by
$$f_1=x^{m_i}g_1\mbox{ and }f_2=x^{n_i}g_2.$$
Let $Q\in C$ be a base point, different from $P_1,\dots,P_M$. Let
$R_{i}$ be a point on $C$, which is within an $\e$-neighborhood of
the intersection point $P_i$. Denote by $\sigma_{i}$ a simple loop
on $C$ around $P_i$, based at $Q$. Let $\gamma_{i}$ be a path from
$Q$ to $R_{i}$. Define also $\sigma^0_i$ to be a (small) simple
loop around $P_i$, based at $R_i$, so that
$$\sigma_{i}=\gamma_i\sigma^0_i\gamma_i^{-1}.$$
Choose the paths $\gamma_{i}$ so that the loop $\delta$,
defined by
$$\delta=\left(\prod_{i=1}^g [\a_i,\b_i]\right)
\left(\prod_{i=1}^N \sigma_{i}\right),$$
is homotopic to the trivial loop at $Q$,
where $g$ is the genus of the curve $C$ and
$[\a_i,\b_i]$ is the commutator of loops $\a_i$ and $\b_i$
around the handles of $C$.
We define the logarithmic Tate symbol by
$$Log\{f_1,f_2\}^{\gamma_i}_{P_i}
=(2\pi i)\left(\pi i m_in_i+
n_i\int_{\gamma_i}\frac{dg_1}{g_1}-m_i\int_{\gamma_i}\frac{dg_2}{g_2}\right).$$
\end{definition}
%%%%%%%%%%%%%%%%%%%%%%%%%%%%%%%%%%%%%%%%%%%%%%%%%%%%%%
\begin{theorem} The relation to the Parshin symbol is
$$\{f_1,f_2\}_{P_i}
=exp((2\pi i)^{-1}Log\{f_1,f_2\}^{\gamma_i}_{P_i})
g_1(Q)^{n_i}g_2(Q)^{m_i},$$
where $Q$ is the base point.
\end{theorem}
%%%%%%%%%%%%%%%%%%%%%%%%%%%%%%%%%%%%%%%%%%%%%%%%%%%%%%%%%%%
A reciprocity law for the logarithm of the Parshin symbol is given by the following theorem.
\begin{theorem}
$$\sum_i Log\{f_1,f_2\}^{\gamma_i}_{P_i}=(2\pi i)^2(M+N),$$
where
$$\begin{tabular}{ll}
$M=$&$\sum_{i} m_{i}n_{i}$\\
\\
$N=$
&$\sum_{i=1}^g\left(\int_{\a_i}\frac{dg_1}{g_1}\int_{\b_i}\frac{dg_2}{g_2}
-\int_{\b_i}\frac{dg_1}{g_1}\int_{\a_i}\frac{dg_2}{g_2}\right).$
\end{tabular}
$$
\end{theorem}
%%%%%%%%%%%%%%%%%%%%%%%%%%%%%%%%%%%%%%%%%%%%%%%%%%%%%
Note that the reciprocity law for the logarithm of the Tate symbol
(Theorem 0.13) resembles the reciprocity law for the logarithm of
the refinement of the Parshin symbol (Theorem 0.6), when we look
at the right hand side of the corresponding equalities. However,
in the reciprocity law for the logarithm of the Parshin symbol
(Theorem 0.10) the right hand side is simply zero, which is
different from the reciprocity law for the logarithm of the Tate
symbol.

The refinement of the Parshin symbol can be defined over any field
$K$. However, our proof of the reciprocity law works only over a
field $K$ of characteristic zero. Let us explain how we relate the
field $K$ to the field of complex numbers $\C$. First, the
rational numbers $\Q$ can be embedded in $K$. The surface $X$ over
$K$ is defined by finitely many polynomials over $K$. We can
adjoin all the coefficient of these polynomials to $\Q$. Then we
obtain a finitely generated algebra over $\Q$. Denote its field of
fractions by $K_0$. Note that $X$ can be defined over $K_0$. Since
the field $K_0$ is of finite transcendence degree over $Q$, we can
embed $K_0$ in $\C$. Let $\bar{K}_0$ be the algebraic closure of
$K_0$. Since $\C$ is algebraically closed, we have that
$\bar{K}_0\subset \C$. Now we can use integration over the complex
numbers to define the refinement. After we have defined the
refinement of the Parshin symbol, we notice that its values are in
$\bar{K}_0$. The reciprocity law is a product over all points of
intersection of $C_0$ with the remaining components of the
divisors. If we consider first the product over the conjugate
points with respect to the Galois group $Gal(\bar{K}_0/K_0)$, we
obtain a symbol which is defined over $K_0$.

We have a higher dimensional analogue of the new symbol and a
reciprocity law for it. It will appear in another paper. Also, my
student Zhenbin Luo \cite{L} has generalized the new symbol on a
surface over a nilpotent extension of the complex numbers. Luo
gives a generalization of Contou-Carrere symbol to surfaces. An
alternative generalization of Contou-Carrere symbol to surfaces is
given by Romo \cite{R}.

An alternative approach to the Tate symbol and the Parshin symbol
come from logarithmic functionals (see \cite{Kh}). The direction
we take is closer to the one in the papers by Deligne \cite{D} and
Brylinski and McLaughlin \cite{BrMc1}. The common point is that
certain connection, in geometric sense, is the key structure. The
difference in our situation is that the connection is not flat.

It will be interesting to find characteristic classes that give
the refinement of the Parshin symbol. For analogues of the Parshin
symbol in terms of characteristic classes (see \cite{BrMc2}). The
use of such more conceptual approach might give a proof for the
reciprocity law of the new symbol for a surface over a finite
field. Note that the refinement of the Parshin symbol can be
defined in the same way for a surface over a finite field. There
must be a version of the new symbol in the arithmetic case for
curves over a number ring.

%%%%%%%%%%%%%%%%%%%%%%%%%%%%%%%%%%%%%%%%%%%%%%%%%%%%%%%%%%%%%%%%%%

Now let us say a few words about the structure of the paper. In
Section 1 we recall basic properties of iterated integrals in the
sense of Chen. The section ends with analogues of the Stokes
formula for iterated integrals in dimension $2$ and $3$ (Theorems
1.5 and 1.9).

In Section 2, we make the key geometric construction, which is a
foundation for the rest of the paper. Besides the geometric
construction, Subsection 2.1 contains most of the definitions
needed for the rest of the paper. In Subsection 2.2, we
construct an abstract reciprocity law,
which we use in the later sections
in order to prove the explicit reciprocity laws,
which we stated in the introduction.

In Section 3, we combine the abstract reciprocity from Section 2
and the properties of iterated integrals, written in Section 1, in
order to obtain a new logarithmic  symbol
$Log[f_1,f_2,f_3]^{x_0,\gamma}_{C_0,P}$. This is the first
logarithmic symbol that we construct. The other symbols are
constructed later in the paper. Also, we prove a reciprocity law
for the logarithm of the new symbol. This logarithmic symbol is
obtained as a limit of an iterated integral over a particular
loop. A key part of this construction is a differential equation
written in Subsection 3.3. This differential equation has no local
solutions. The way we should solve the differential equation is by
restricting it to a path. One may think of it as a connection,
which is not flat. A solution of this differential equation is a
generating series of iterated integrals. However, we do not need
the whole generating series. We need only specific term of it,
which correspond to one iterated integral. This is where the
analogy with a non-flat connection breaks. For that reason we call
it a differential equation.

In Section 4, we construct a logarithm of the Parshin symbol and
prove a reciprocity law for it. Besides many computations of
iterated integrals based on the definitions in Subsection 2.4, we
use one of the differential equations from subsection 3.3. In that
differential equation we put an equivalence among the formal
non-commuting variables. This technical condition corresponds to
considering a linear combination of iterated integrals as opposed
to one iterated integral. A particular linear combination of
iterated integrals gives us the logarithm of the Parshin symbol
$Log\{f_1,f_2,f_3\}^{\gamma}_{C_0,P}$, which we discussed earlier
in the introduction. From the logarithm of the Parshin symbol, we
can recover the Parshin symbol by exponentiation.

In Section 5, we define the logarithm of the refinement of the
Parshin symbol as a difference
$$Log(f_1,f_2,f_3)^{x_0,\gamma}_{C_0,P}
=Log[f_1,f_2,f_3]^{x_0,\gamma}_{C_0,P}-Log\{f_1,f_2,f_3\}^{\gamma}_{C_0,P}.$$
Using the reciprocity laws for the other two symbols, we obtain a
reciprocity law for the logarithm of the refinement of the Parshin
symbol. After exponentiation of the logarithm of the refinement of
the Parshin symbol, we obtain a refinement of the Parshin symbol,
which is roughly speaking $1/3$ of the Parshin symbol. We end the
paper with an example of the reciprocity law for the
(multiplicative) refinement of the Parshin symbol.

{\bf Acknowledgments.} I would like to thank Alexey Parshin and
Pierre Deligne for the useful commentaries on the earlier version
of the paper. I would like to thank Ivan Fesenko for fruitful
conversation on tame symbols. Also, I would like to thank Zhenbin
Luo for lengthy discussions on my approach as well as to thank
Anton Deitmar for the interest in this paper.

I would like to thank the University of Durham for the kind
hospitality and also to thank the Arithmetic Algebraic Geometry
Marie Curie Network for the financial support. I would like to
thank both Brandeis University, where this work was developed and
Universit\"at T\"ubingen, where this work was finished.

%%%%%%%%%%%%%%%%%%%%%%%%%%%%%%%%%%%%%%%%%%%%%%%%%%%%%%%%
\section{Background on iterated integrals}
In this section we recall known properties of iterated integrals,
which we are going to use heavily in the rest of the paper. One
can look at \cite{Ch} and \cite{G} for more properties of iterated
integrals. We include this section, because it is essential for
the rest of the paper. This section establishes both the notation
and the main properties of iterated integrals, which we are going
to use throughout the paper.
\subsection{Definition of iterated integrals over a path}
\begin{definition}
\label{def1} Let $\omega_1,\dots,\omega_n$ be holomorphic 1-forms
on a simply connected open subset $U$ of the complex plane
${\mathbb C}$. Let
$$\gamma:[0,1]\rightarrow U$$ be a path.
We define an iterated integral of the forms
$\omega_1,\dots,\omega_n$ over the path $\gamma$ to be
$$\int_\gamma \omega_1\circ\dots\circ\omega_n=
\int\dots\int_{0\leq t_1\leq\dots\leq t_n\leq 1}
\gamma^*\omega_1(t_1)\wedge\dots\wedge\gamma^*\omega_n(t_n).$$
\end{definition}

It is called iterated because it can be defined inductively by
$$\int_\gamma \omega_1\circ\dots\circ\omega_n=
\int_0^1
(\int_{\gamma|[0,t]}\omega_1\circ\dots\circ\omega_{n-1})\gamma^*\omega_n(t).$$
%%%%%%%%%%%%%%%%%%%%%%%%%%%%%%%%%%%%%%%%%%%%%%%%%%%%%%%%%%%%%%%%%%%%%%%%%%%%%%%%%%

%%%%%%%%%%%%%%%%%%%%%%%%%%%%%%%%%%%%%%%%%%%%%%%%%%%%%%%%%%%%%%%%%%%%%%%%%%%%%%%%
\subsection{Differential equation and generating series of
iterated integrals. One-dimensional case}

When we consider an iterated integral, we can let the end point
vary in a small neighborhood. Then the iterated integral becomes
an analytic function.

Let $\omega_1,\dots,\omega_n$ be differentials of 3rd kind on a
Riemann surface $X$. Following the idea of Manin \cite{M}, we
consider the differential equation
$$dF=F\sum_{i=1}^n A_i\omega_i,$$
where $A_1,\dots,A_n$ are non-commuting formal variables, which
commute with the differentiation $d$. Let $P$ be a point of $X$
such that none of the differential forms has a pole at $P$. It is
easy to check that the function
$$F(z)=1+\sum_{i}(A_i\int_P^z \omega_i) +
\sum_{i,j}A_iA_j\int_P^z\omega_i\circ\omega_j +
\sum_{i,j,k}A_iA_jA_k\int_P^z\omega_i\circ\omega_j\circ\omega_k+\dots.$$
is a solution to the differential equation with initial condition
$F(P)=1$ at the point $P$.  The summation continues so that every
iterated integral of the given $n$ 1-forms is present in the
summation. Note that
$$d\int_P^z \omega_{i_1}\circ\dots\circ\omega_{i_{k-1}}\circ\omega_{i_k} =
(\int_P^z\omega_{i_1}\circ\dots\circ\omega_{i_{k-1}})\omega_{i_k}.$$
Note that the coefficient of the above integrals in the solution
$F$ is
$$A_{i_1}\dots A_{i_{k-1}}A_{i_k},$$
whose indices enumerate the order of iteration of the differential
forms. Each of the indices $i_1,\dots,i_{k-1},i_k$ is an integer
among $\{1,2,\dots,n\}$ and repetitions of indices is allowed.
%%%%%%%%%%%%%%%%%%%%%%%%%%%%%%%%%%%%%%%%%%%%%%%%%%%%%%%%%%%%%%
\subsection{Multiplication formulas}
We can take a path $\gamma$ from $P$ to $z$. We denote the
solution of the differential equation by $F_{\gamma}$. If
$\gamma_1$ is a path that ends at $Q$ and $\gamma_2$ is a path
that starts at $Q$ we can compose them. Denote the composition by
$\gamma_1\gamma_2$. %%%%%%%%%%%%%%%%%%%%%%%%%%%%%%%%%%%%%%%%%%%%%%%%%
%%%%%%%%%%%%%%%%%%%%%%%%%%%%%%%%%%%%%%%%%%%%%%%%%%%%%%%%%%%%%%%%%%%%%
\begin{theorem} \label{CompPath} (Composition of paths) With the above notation, we have
$$F_{\gamma_1}F_{\gamma_2}=F_{\gamma_1\gamma_2}.$$
\end{theorem}
%%%%%%%%%%%%%%%%%%%%%%%%%%%%%%%%%%%%%%%%%%%%%%%%%%%%%%%%%%%%%%%%%%
\begin{corollary} \label{co1}(Composition of paths) On the level
of iterated integrals we have
$$\int_{\gamma_1\gamma_2} \omega_1\circ\dots\circ\omega_n=
\sum_{i=0}^n \int_{\gamma_1}\omega_1\circ\dots\circ\omega_i
\int_{\gamma_2} \omega_{i+1}\circ\dots\circ\omega_n,$$ where for
$i=0$ we define
$\int_{\gamma_1}\omega_1\circ\dots\circ\omega_i=1$, and similarly,
for $i=n$ we define $\int_{\gamma_2}
\omega_{i+1}\circ\dots\circ\omega_n=1$.
 \end{corollary}
%%%%%%%%%%%%%%%%%%%%%%%%%%%%%%%%%%%%%%%%%%%%%%%%%%%%%%%%%%%%%%%%%%%%%%%%
\subsection{Two-dimensional iterated integral over a surface}

In this subsection we recall certain type of Chen's iterated
integrals \cite{Ch}, which we call $2$-dimensional iterated
integrals. These integrals are defined by integration over a
(real) $2$-dimensional region of a variety. After we recall the
definition and give notation, we state an analogue of Stokes
theorem for $2$-dimensional iterated integrals, which expresses
$1$-dimensional iterated integrals over the boundary of a two
dimensional region in terms of $2$-dimensional iterated integrals
over the same region.

Let $X$ be a smooth complex manifold of dimension at least $2$.
Let $\omega_1,\dots,\omega_n$ be closed holomorphic differential
$1$-forms on $X$. Let $\omega^{(2)}$ be a closed holomorphic
differential $2$-form on $X$. Let
$$\gamma_\bullet:[0,1]\times[0,1]\rightarrow X$$
be a homotopy of the path
$$\gamma_0:[0,1]\rightarrow X$$
and
$$\gamma_1:[0,1]\rightarrow X,$$
which fixes the end points, such that
$\gamma_\bullet(t,0)=\gamma_0(t)$ and
$\gamma_\bullet(t,1)=\gamma_1(t)$. Define also a domain
$$\Delta^n_\bullet=\{(t_1,\cdots,t_n,s)\in [0,1]^{n+1}
|0\leq t_1\leq\cdots\leq t_{n},\mbox{ }0\leq s\leq 1 \}.$$

\begin{definition} A $2$-dimensional iterated integral over
$\gamma_\bullet$ is defined by
$$\begin{tabular}{ll}
$I^i_\bullet=\int\int_{\gamma_\bullet}\omega_1\circ\cdots\circ
\omega_{i-1}\circ\omega^{(2)}\circ\omega_{i+1}\circ\cdots\circ\omega_n=$\\
\\
$=\int_{\Delta_\bullet^n}
\gamma_\bullet^*\omega_1(t_1,s)\wedge\cdots\wedge
\gamma_\bullet^*\omega_{i-1}(t_{i-1},s)\wedge
\gamma_\bullet^*\omega^{(2)}(t_i,s)\wedge
\gamma_\bullet^*\omega_{i+1}(t_{i+1},s) \wedge\cdots\wedge
\gamma_\bullet^*\omega_n(t_{n},s)$.
\end{tabular}
$$
\end{definition}
For fixed $s$ let
$$I_s=\int_{\gamma_s}\omega_1\circ\cdots\omega_n.$$
%%%%%%%%%%%%%%%%%%%%%%%%%%%%%%%%%%%%%%%%%%%%%%%%%%%%%%%%%%%%%%%%%%%%%%%%%%%%%%
Now let us recall Stokes theorem for two dimensional iterated
integrals.
\begin{theorem}
Let
$$\widetilde{\omega}_i^{(2)}=\omega_i\wedge\omega_{i+1}.$$
The $2$-form $\widetilde{\omega}_i^{(2)}$ will be used in the
definition of $I^i_\bullet$. Then
$$I_0-I_1=\sum_{i=1}^{n-1} I^i_\bullet.$$
\end{theorem}

%%%%%%%%%%%%%%%%%%%%%%%%%%%%%%%%%%%%%%%%%%%%%%%%%%%%%%%%%%%%%%%%%%%%%%

\subsection{Slicing a membrane}
In this subsection, we examine what happens to the $2$-dimensional
iterated integrals, when the $2$-dimensional region is cut into
two pieces. We write an analogue of Theorem \ref{CompPath} and
Corollary \ref{co1} for $1$-dimensional iterated integrals, when
the path is cut in two pieces. Both of these formulas are given on
the level of iterated integrals and on the level of generating
series of iterated integrals.

Consider the square $[0,1]\times [0,1]$. We will separate it into
two domains in the following way.

Let $\gamma_0$ be the path that follows the lower edge and the
right edge of the square. Namely,
$$\gamma_0:[0,2]\rightarrow [0,1]^2$$
$$\gamma_0(t)=\left\{
\begin{tabular}{ll}
$(t,0)$ & for $0\leq t \leq 1$\\
$(1,t-1)$ & for $1\leq t \leq 2$
\end{tabular}
\right.$$ Let $\gamma_1$ be the path that follows the left edge
and the upper edge of the square. Namely,
$$\gamma_1:[0,2]\rightarrow [0,1]^2$$
$$\gamma_1(t)=\left\{
\begin{tabular}{ll}
$(0,t)$ & for $0\leq t \leq 1$\\
$(t-1,1)$ & for $1\leq t \leq 2$
\end{tabular}
\right.$$ We want to consider a homotopy $\gamma_s$ between the
two paths $\gamma_0$ and $\gamma_1$. Let
$$\gamma_s(t)=\left\{
\begin{tabular}{ll}
$(0,t)$ & for $0\leq t \leq s$\\
$(t-s,s)$ & for $s\leq t \leq 1+s$\\
$(1,t-1)$ & for $1+s\leq t \leq 2$
\end{tabular}
\right.$$

Consider the two rectangles $S_0=[0,1/2]\times[0,2]$ and
$S_1=[1/2,1]\times[0,2]$, where $s$ varies in the first interval
and $t$ varies in the second interval. Let $\sigma_0(t)=(t,0)$ for
$0\leq t \leq 1/2$ and $\sigma_1(t)=(t,1)$ for $1/2\leq t \leq 1$.
Then $\gamma_s$ for $0\leq s \leq 1/2$ has domain
$S_0\cup\sigma_1$ and $\gamma_s$ for $1/2\leq s \leq 1$ has domain
$\sigma_0\cup S_1$. We are going to write $S_0\sigma_1$ instead of
$S_0\cup\sigma_1$, when the end of the paths in $S_0$ is the
beginning of the path $\sigma_1$.

Let $$\gamma_\bullet(s,t)=\gamma_s(t).$$ Then we can slice the
square $[0,1]\times [0,1]$ into the following two domains:
$$\gamma_\bullet|S_0\sigma_1$$ and
$$\gamma_\bullet|\sigma_0S_1.$$ Let also $$S=S_0\sigma_1\cup
\sigma_0S_1.$$

As always we denote by $s$ the variable, which parameterizes the
variation of paths. Note that in the $2$-dimensional iterated
integrals the variables $t_1,t_2,\cdots$ cannot be permuted.
However, in direction of increasing $s$, we have "commutativity"
in the sense that we can integrate first $s$ in the interval
$[1/2,1]$ and then add this to the integral of $s$ in the interval
$[0,1/2]$, where in both integrals we have the same domain for
$t_1,t_2,\cdots$. Let
$$\tau_\bullet:S\rightarrow X$$
be a variation of paths on $X$ and let $\omega_1,\cdots ,\omega_n$
be closed 1-forms on $X$. Let
$$f_idt=\tau_\bullet^*\omega_i$$
on $S$. Then we have the following lemma.
%%%%%%%%%%%%%%%%%%%%%%%%%%%%%%%%%%%%%%%%%%%%%%%%%%%%%%%%%%%%%%%%%%%%%
\begin{lemma}
$$\left(\int_S - \int_{S_0\sigma_1}-\int_{\sigma_0S_1}\right)f_1dt_1\wedge\cdots
\wedge f_{i-1}dt_{i-1}\wedge(f_if_{i+1}ds\wedge t_i)\wedge
f_idt_{i+2}\wedge\cdots \wedge f_ndt_n=0.$$
\end{lemma}
%%%%%%%%%%%%%%%%%%%%%%%%%%%%%%%%%%%%%%%%%%%%%%%%%%%%%%%%%%%%%%%%%%%%%%%%%%%%%%%%%%%
\subsection{Composition of a path and a $2$-dimensional region} Let $\gamma_\bullet$ be a
homotopy of paths, fixing the end points. The $2$-dimensional
region is the one covered by $\gamma_\bullet$. Let $\sigma$ be a
path, whose starting point is the ending point of the paths in
$\gamma_\bullet$. Let $\gamma_\bullet \sigma$ be the composition
of a path from $\gamma_\bullet$ with the path $\sigma$. We want to
express a $2$-dimensional iterated integral over $\gamma_\bullet
\sigma$ in terms of iterated integrals over $\sigma$ and
$2$-dimensional iterated integrals over $\gamma_\bullet$.

These integrals are related by the following theorem.
%%%%%%%%%%%%%%%%%%%%%%%%%%%%%%%%%%%%%%%%%%%%%%%%%%%%%%%%
\begin{theorem}
Let $\omega_1,\cdots,\omega_n$ be closed $1$-forms. Then
$$\int_{\gamma_\bullet\sigma}
\omega_1\circ\cdots\circ(\omega_i\wedge\omega_{i+1})\circ\cdots\circ\omega_n=
\sum_{j=i+1}^n \int_{\gamma_\bullet}
\omega_1\circ\cdots\circ(\omega_i\wedge\omega_{i+1})\circ\cdots\circ\omega_j
\int_\sigma \omega_{j+1}\circ\cdots\circ\omega_n.$$
\end{theorem}
\proof The proof is the same as in the case of $1$-dimensional
iterated integrals. The only difference is that if $j<i$ then the
domain for the first integral will have dimension $j+1$, which is
larger than the number of differential $1$-forms, which is $j$.
Thus, the first integral will be zero.

%On the level of generating series, the above theorem becomes more
%compact.
%%%%%%%%%%%%%%%%%%%%%%%%%%%%%%%%%%%%%%%%%%%%%%%%%%%%%%%%%%%%%%%%%%%%%%%%%%%%%%%%%%%%%%%%%%
%\begin{theorem}
%$$F_{\gamma_\bullet\sigma}=F_{\gamma_\bullet}F_{\sigma}.$$
%\end{theorem}
%%%%%%%%%%%%%%%%%%%%%%%%%%%%%%%%%%%%%%%%%%%%%%%%%%
%%%%%%%%%%%%%%%%%%%%%%%%%%%%%%%%%%%%%%%%%%%%%%%%%%%%%%%%%%%%%%%%%%%%%%
\subsection{Three-dimensional iterated integrals}

In this section we recall certain type of Chen's iterated
integrals, which we call $3$-dimensional iterated integrals. These
integrals are defined by integration over a (real) $3$-dimensional
region of a variety. After we recall the definition and give
notation, we state an analogue of Stokes theorem for
$3$-dimensional iterated integrals, which expresses
$2$-dimensional iterated integrals over the boundary of a
$3$-dimensional region in terms of $2$-dimensional iterated
integrals over the same region.

Let $X$ be a smooth complex manifold of dimension at least $3$.
Let $\omega_1,\dots,\omega_n$ be closed holomorphic differential
$1$-forms on $X$. Let $\omega_1^{(2)},\cdots,\omega_n^{(2)}$ be
closed holomorphic differential $2$-forms on $X$. Let
$\omega^{(3)}$ be a closed holomorphic differential $2$-form on
$X$. Let
$$\gamma_{\bullet 0}:[0,1]\times[0,1]\rightarrow X$$ and
$$\gamma_{\bullet 1}:[0,1]\times[0,1]\rightarrow X$$
be two homotopies of paths, which fix the end points. Let
$$\gamma_{\bullet\bullet}:[0,1]^3\rightarrow X$$
be a homotopy of the two homotopies of the path
$$\gamma_{\bullet 0}:[0,1]^2\rightarrow X$$
and
$$\gamma_{\bullet 1}:[0,1]^2\rightarrow X,$$
which fixes the end points, such that
$\gamma_{\bullet\bullet}(t,s_1,0)=\gamma_{\bullet 0}(t,s_1)$ and
$\gamma_{\bullet\bullet}(t,s_1,1)=\gamma_{\bullet 1}(t,s_1)$.
Define also the following homotopies of paths
$$\gamma_{0\bullet}(t,s_2)=\gamma_{\bullet\bullet}(t,0,s_2)$$ and
$$\gamma_{1\bullet}(t,s_2)=\gamma_{\bullet\bullet}(t,1,s_2).$$
Define also a domain
$$\Delta^n_{\bullet\bullet}=\{(t_1,\cdots,t_{n},s_1,s_2)\in [0,1]^{n+2}
|0\leq t_1\leq\cdots\leq t_{n},\mbox{ }0\leq s_1\leq 1,\mbox{ }
0\leq s_2\leq 1\}.$$ We will omit the variables $s_1$ and $s_2$
from the notation of
$\gamma_{\bullet\bullet}^*\omega_i(s_1,s_2,t_i),$
$\gamma_{\bullet\bullet}^*\omega_i^{(2)}(s_1,s_2,t_i)$ and
$\gamma_{\bullet\bullet}^*\omega_i^{(3)}(s_1,s_2,t_i).$ We will
write $\gamma_{\bullet\bullet}^*\omega_i(t_i),$
$\gamma_{\bullet\bullet}^*\omega_i^{(2)}(t_i)$ and
$\gamma_{\bullet\bullet}^*\omega_i^{(3)}(t_i),$ respectively, in
their place.
%%%%%%%%%%%%%%%%%%%%%%%%%%%%%%%%%%%%%%%%%%%%%%%%%%%%%%%%%%%%%%%%%%%%%%
\begin{definition} A $3$-dimensional iterated integral over
$\gamma_{\bullet\bullet}$ can be of two types. One of the types is
defined by
$$\begin{tabular}{ll}
$I^{i,i+1}_{\bullet\bullet}=\int\int\int_{\gamma_{\bullet\bullet}}\omega_1\circ\cdots\circ
\omega_{i-1}\circ\omega^{(3)}\circ\omega_{i+1}\circ\cdots\circ\omega_n=$\\
\\
$=\int_{\Delta_{\bullet\bullet}^n}
\gamma_{\bullet\bullet}^*\omega_1(t_1)\wedge\cdots\wedge
\gamma_{\bullet\bullet}^*\omega_{i-1}(t_{i-1})\wedge
\gamma_{\bullet\bullet}^*\omega^{(3)}(t_i)\wedge
\gamma_{\bullet\bullet}^*\omega_{i+1}(t_{i+1}) \wedge\cdots\wedge
\gamma_{\bullet\bullet}^*\omega_n(t_{n})$.
\end{tabular}
$$
And the other type is defined for $i<j-1$
$$\begin{tabular}{ll}
$I^{i,j}_{\bullet\bullet}=\int\int\int_{\gamma_{\bullet\bullet}}\omega_1\circ\cdots\circ
\omega_{i-1}\circ\omega_i^{(2)}\circ\omega_{i+1}\circ\cdots\circ\omega_{j-1}
\circ\omega_j^{(2)}\circ\omega_{j+1}\circ\cdots\circ\omega_{n}=$\\
\\
$=\int_{\Delta_{\bullet\bullet}^n}
\gamma_{\bullet\bullet}^*\omega_1(t_1)\wedge\cdots\wedge
\gamma_{\bullet\bullet}^*\omega_{i-1}(t_{i-1})\wedge
\gamma_{\bullet\bullet}^*\omega_i^{(2)}(t_i)\wedge
\gamma_{\bullet\bullet}^*\omega_{i+1}(t_{i+1}) \wedge\cdots$\\
\\
$\cdots\wedge\gamma_{\bullet\bullet}^*\omega_{j-1}(t_{j-1})\wedge
\gamma_{\bullet\bullet}^*\omega_j^{(2)}(t_{j})\wedge
\gamma_{\bullet\bullet}^*\omega_{j+1}(t_{j+1}) \wedge \cdots\wedge
\gamma_{\bullet\bullet}^*\omega_n(t_{n})$.
\end{tabular}
$$
And set
$$I^{j,i}_{\bullet\bullet}=I^{i,j}_{\bullet\bullet}.$$
\end{definition}
%%%%%%%%%%%%%%%%%%%%%%%%%%%%%%%%%%%%%%%%%%%%%%%%%%%%%%%%%%%%%%%%%%%%%%%%%%%%%%
Now let us recall Stokes theorem for $3$-dimensional iterated
integrals. Let
$$\widetilde{\omega}_j^{(2)}=\omega_j\wedge\omega_{j+1}.$$ Let
$$\widetilde{\omega}_i^{(3)}=\omega^{(2)}_i\wedge\omega_{i+1}$$
and let
$$\widetilde{\omega}_{i-1}^{(3)}=\omega_{i-1}\wedge\omega^{(2)}_i$$
%%%%%%%%%%%%%%%%%%%%%%%%%%%%%%%%%%%%%%%%%%%%%%%%%%%%%%%%%%%%%%%%%%%%%%
\begin{theorem} The forms $\widetilde{\omega}_j^{(2)}$,
$\widetilde{\omega}_i^{(3)}$ and $\widetilde{\omega}_{i-1}^{(3)}$
will enter in the definition of the $3$-dimensional iterated
integrals. The $2$-form $\omega^{(2)}_i$ will enter in the
definition of the $2$-dimensional iterated integrals. Then
$$I^i_{0\bullet}-I^i_{1\bullet}+I^i_{\bullet 0}-I^i_{\bullet 1}
=2\sum_{j\neq i} I^{i,j}_{\bullet\bullet}.$$
\end{theorem}
%%%%%%%%%%%%%%%%%%%%%%%%%%%%%%%%%%%%%%%%%%%%%%%%%%%%%%%%%%%%%%%%%%%%%%%
We will use a particular case of this theorem, which we will
formulate in a corollary.
%%%%%%%%%%%%%%%%%%%%%%%%%%%%%%%%%%%%%%%%%%%%%%%%%%%%%%%%%%%%%%%%%%%%%%%%%
\begin{corollary} Over a two dimensional complex
manifold, we have
$$I^i_{0\bullet}-I^i_{1\bullet}+I^i_{\bullet 0}-I^i_{\bullet 1}=0.$$
\end{corollary}
%%%%%%%%%%%%%%%%%%%%%%%%%%%%%%%%%%%%%%%%%%%%%%%%%%%%%%%%%%%%%%%%%%%%%%%%%%%%%%%%%
\section{Abstract reciprocity law}
%%%%%%%%%%%%%%%%%%%%%%%%%%%%%%%%%%%%%%%%%%%%%%%%%%%%%%%%%%%%%
\subsection{Definitions}
Consider $X$ a projective smooth surface over
the complex numbers. Let $f_k$ for $k=1,2,3$ be rational
functions on $X$.
%%%%%%%%%%%%%%%%%%%%%%%%%%%%%%%%%%%%%%%%%
\begin{definition} (Divisors $C_i$ and $D_{ij}$)
Denote by $C_i$ for $i=0,1,\dots,N$
the components of the divisors of $f_1$, $f_2$ and $f_3$.
Assume that the curves $C_i$ are non-singular and
that they intersect properly.
Let
$$div(f_k)=\sum_{i=0}^N n_{ki} C_i.$$
Let $x_i$ be a rational function on $X$,
which has zero along $C_i$ of order $1$,
and has no zeroes or poles along $C_j$ for $j\neq i$. Let
$$div(x_i)=C_j+\sum_{j=1}^{M_i}m_{ij}D_{ij}.$$
We choose $x_i$ for $i=0,1,\dots,M$ so that $D_{ij}\neq D_{i'j'}$
for $(i,j)\neq (i',j')$. Let $y_{ij}$ be a rational function on $X$,
which has  zero along $D_{ij}$ of order $1$,
and has no zeroes or poles along $C_{i'}$ for $i'=0,1,\dots,M$
and along $D_{i'j'}$ for $(i',j')\neq (i,j)$.
\end{definition}
%%%%%%%%%%%%%%%%%%%%%%%%%%%%%%%%%%%%%%%%%%%%

We will examine what happens along the curve $C_0$. Choose a small
positive number $\e$. Let $0<\e_0<\e$ and $0<\e_1<\e$. We will
define a tubular neighborhood of radius $\e_0$ around the
following space: the curve $C_0$ minus $\e_1$ neighborhoods of the
intersections of $C_0$ with $C_j$ for $j>0$ and with $D_{ij}$ for
all $i$ and $j$. Near the intersection points we will define tori.
We are going to make precise what means $\e_0$ or $\e_1$
neighborhood. We need this $2$-dimensional region in order to
examine iterated integrals over it. Afterwards, we are going to
cut this $2$-dimensional region so that topologically it looks
like a rectangle, where each side is a path on $X$, which lies on
an algebraic curve.
%%%%%%%%%%%%%%%%%%%%%%%%%%%%%%%%%%%%%%%%%%%%%%%%%%%%%%%%%%%%%%%%%%%%%%%
\begin{definition} (Intersection points $P_{jl}$ and $P_{1jl}$)
Let $$\bigcup_{j=1}^{N}P_{jl}=C_1\cap C_j$$
and let
$$\bigcup_{j=1}^{M_{1}}P_{1jl}=C_1\cap D_{1j}.$$
\end{definition}
%%%%%%%%%%%%%%%%%%%%%%%%%%%%%%%%%%%%%%%%%%%%%%%%%%%%%%%%%
We will construct compact neighborhoods,
which separate the points $P_{jl}$ and
$P_{1jl}$.
%%%%%%%%%%%%%%%%%%%%%%%%%%%%%%%%%%%%%%%%%%%%%%%%%%%%%%%%%%%%%%%%%%%%%%%
\begin{definition} (Neighborhoods $\tilde{T}_{jl}$ and $\tilde{T}_{1jl}$
of the points $P_{jl}$ and $P_{1jl}$)
Let
$$E_j=\{P\in X| |x_0(P)|<\e_0 \mbox{ and } |x_j(P)|<\e_1\}.$$
Define
$$E_{1j}=\{P\in X| \left|
\frac{x_0(P)}{y_{1j}(P)}\right|<\frac{\e_0}{\e_{1j}^{m_{1j}}}
\mbox{ and }  |y_{1j}(P)|<\e_{1j}\}.$$
Denote by $\tilde{T}_{jl}$ the connected component in $E_j$
around the point $P_{jl}$. Also,
denote by $\tilde{T}_{1jl}$ the connected component in $E_{1j}$
around the point $P_{1jl}$.
\end{definition}
%%%%%%%%%%%%%%%%%%%%%%%%%%%%%%%%%%%%%%%%%%%%%%%%%%%%%%%%%%%%%%%%%%%%%%
We will find relations between
the exponents $m_{1j}$ and the small constants $\e_{1j}$
so that the set $E_{1j}$ consists of small compact sets.

If $m_{1j}>0$ then we want $\frac{\e_0}{\e_{1j}^{m_{1j}}}<\e$ and
$\e_{1j}<\e$. This can be achieved by taking $$\e_0<\e^{1+m_{1j}}$$ and
$$\e_{1j}^{m_{1j}}<\min(\e^{m_{1j}},\frac{\e_0}{\e}).$$

If $m_{1j}<0$ then we want $\frac{\e_0}{\e_{1j}^{m_{1j}}}<\e$ and
$\e_{1j}<\e$. This can be achieved by taking $$\e_0<\e$$ and
$$\e_{1j}<\e.$$
Globally, we can choose $\e_0$ so that $\e_0<\e^{1+\max m_{1j}},$
where the maximum is taken over all $j$'s such that $m_{1j}>0$.
If all $m_{1j}<0$ then choose $\e_0<\e$. For $\e_{1j}$, we choose
$\e_{1j}<\e$ if $m_{1j}<0$ and
$\e_{1j}^{m_{1j}}<\min(\e^{m_{1j}},\frac{\e_0}{\e}).$

For fixed $j$, we can choose $\e$
small enough so that the compact sets $E_j$
is a disjoint union of neighborhoods - one for each point $P_{jl}$.
Similarly, for fixed $j$, we can choose $\e$
small enough so that the compact sets $E_{1j}$
is a disjoint union of neighborhoods - one for each point $P_{1jl}$.
Pick the smallest values of $\e$ among the above choices.
We can choose $\e$ even smaller,
so that the compact neighborhoods $\tilde{T}_{jl}$
for all $j$ and $l$ and the compact neighborhoods
$\tilde{T}_{1jl}$ for all $j$ and $l$ are disjoint.
%%%%%%%%%%%%%%%%%%%%%%%%%%%%%%%%%%%%%%%%%%%%%%%%%%%%%%%%%%%%%%%%%%%%%%%%%%%%
\begin{definition} (Tubular neighborhood)
Let
$$Tb=\{P\in X| |x_1|<\e_0\}
-\bigcup_{j,l}\tilde{T}_{jl}-\bigcup_{j,l}\tilde{T}_{1jl}.$$
\end{definition}
%%%%%%%%%%%%%%%%%%%%%%%%%%%%%%%%%%%%%%%%%%%%%%%%%%%%%%%%%%%%%%%%%%%%%%%%%%
\begin{definition} (Foliation of the tubular neighborhood)
Let $u$ be a complex number such that $|u| \leq \e_0$.
Define $$Tb_u=\{P\in X| x_1(P)=u\}.$$
Note that $Tb_0\subset C_0$.
\end{definition}
%%%%%%%%%%%%%%%%%%%%%%%%%%%%%%%%%%%%%%%%%%%%%%%%%%%%%%%%%%%%%%%%%%%%%%%%%%
%%%%%%%%%%%%%%%%%%%%%%%%%%%%%%%%%%%%%%%%%%%%%%%%%%%%%%%%
\begin{definition}
Define $$D_{\e_0}=\{u\in \C| |u|<\e_0\}.$$
\end{definition}
%%%%%%%%%%%%%%%%%%%%%%%%%%%%%%%%%%%%%%%%%%%%%%%%%%%%%%%%%%%%%%%%%%%
\begin{lemma} (fibration)
For all $u<\e_0$, the sequence of topological space
$$Tb_u \rightarrow Tb \rightarrow D_{\e_0}$$
is a fibration, which splits, since $D_{\e_0}$ is contractible.
\end{lemma}
%%%%%%%%%%%%%%%%%%%%%%%%%%%%%%%%%%%%%%%%%%%%%%%%%%%%%%%%%%%%%%%
\proof The map $$\bar{x}_0:Tb \rightarrow D_{\e_0}$$
is a restriction of the map
$$x_0:X\rightarrow \C,$$
which is singular only at finitely many points in the range $D_{\e_0}$.
Thus, by decreasing $\e_0$, we can assume that
$\bar{x}_0$ is singular at most at one point of $D_{\e_0}$.
If it turns out to be non-singular
then we have a fibration.

Assume the map $\bar{x}_0$ has one singular point $u_0$ in $D_{\e_0}$.
If $u_0\neq 0$ then we can decrease $\e_0$ so that $\e_0<u_0$. Then
$\bar{x}_0$ will not have a singular point.

Assume that $\bar{x}_0$ is singular at $0\in D_{\e_0}$.
Consider the boundaries of $Tb_0$ and $Tb_u$.
Note that the connected components
of their boundaries are circles.
Also, each of the two spaces have the same number of
connected components of their boundaries.
We claim that small neighborhoods of their boundaries are homotopic.
Let $\e'$ be a small positive real number.
Then connected components of the domains
$$\{P\in Tb_u|\e_1<|x_j(P)|<(1+\e')\e_1\}$$
and
$$\{P\in Tb_0|\e_1<|x_j(P)|<(1+\e')\e_1\}$$
are cylinders, corresponding to the points $P_{jl}$
for admissible values of $l$.
Similarly,
connected components of the domains
$$\{P\in Tb_u|\e_{1j}<|y_{1j}(P)|<(1+\e')\e_{1j}\}$$
and
$$\{P\in Tb_0|\e_{1j}<|y_{1j}(P)|<(1+\e')\e_{1j}\}$$
are cylinders, corresponding to the points $P_{1jl}$ for
admissible values of $l$. By decreasing $u$ we obtain a
degeneration of the topological surface $Tb_u$ to $Tb_0$. Now one
can use Riemann-Hurwitz Theorem. Since the number of connected
components of the boundaries of
 $Tb_0$ and $Tb_u$ are the same
and also small neighborhoods of the boundaries are homotopic, we
obtain that the topological surfaces  $Tb_0$ and $Tb_u$ have the
same genus. Thus,  $Tb_0$ and $Tb_u$ are homotopic. This proves
that we have a fibration map $\bar{x}_0$.
%%%%%%%%%%%%%%%%%%%%%%%%%%%%%%%%%%%%%%%%%%%%%%%%%%%%%%%
A direct consequence of the above lemma is the following
corollary.
%%%%%%%%%%%%%%%%%%%%%%%%%%%%%%%%%%%%%%%%%%%%%%%%%
\begin{corollary} (homotopy)
For all $u<\e_0$ there is a homotopy map
$$h:D_{\e_0}\times Tb_0 \rightarrow Tb,$$
such that
$$h(u,Tb_0)=Tb_u.$$
\end{corollary}
%%%%%%%%%%%%%%%%%%%%%%%%%%%%%%%%%%%%%%%%%%%%%%%%%%%%%%%%%%%%%%%%%%
%%%%%%%%%%%%%%%%%%%%%%%%%%%%%%%%%%%%%%%%%%
\begin{definition} (Points, loops and paths on $Tb_0$)
Let $Q\in Tb_0$ be an interior point.
Let $R_{jl}$ and $R_{1jl}$ be base points on each of the connected components of
the boundary of $Tb_0$. We choose $R_{jl}$ and $R_{1jl}$ so that they are close to the
intersection points $P_{jl}$ and $P_{1jl}$.

Let $\sigma_{jl}'$ be a loop on $Tb_0$ defined by
$$\sigma_{jl}^0:[0,1]\rightarrow T_{jl}$$
$$t\mapsto (x_0,x_j)^{-1}(0,\e_1 e^{2\pi i t}),$$
starting at $R_{jl}=(x_0,x_j)^{-1}(0,\e_1).$
Let also $\sigma_{1jl}^0$ be a loop on $Tb_0$ defined by
$$\sigma_{1jl}^0:[0,1]\rightarrow T_{1jl}$$
$$t\mapsto (x_0,y_{1j})^{-1}(0,\e_{1j} e^{2\pi i t}),$$
starting at $R_{1jl}=(x_0,x_j)^{-1}(0,\e_{1j}).$

Let $\gamma_{jl}$
be a path from $Q$ to $R_{jl}$ and let $\gamma_{1jl}$
be a path from $Q$ to $R_{1jl}$.

Denote by $$\sigma_{jl}=\gamma_{jl}\sigma_j^0\gamma_{jl}^{-1}$$
and $$\sigma_{1jl}=\gamma_{1jl}\sigma_{1j}^0\gamma_{1jl}^{-1}.$$
Choose the paths $\gamma_{jl}$ and $\gamma_{1jl}$ so that the loop $\delta$,
defined by
$$\delta=\left(\prod_{i=1}^g [\a_i,\b_i]\right)
\left(\prod_{j=1}^N \sigma_{jl}\right)
\left(\prod_{j=1}^{M_1} \sigma_{1jl}\right),$$
is homotopic to the trivial loop at $Q$,
where $g$ is the genus of the curve $Tb_0$ and
$[\a_i,\b_i]$ is the commutator of loops $\a_i$ and $\b_i$
around the handles of $Tb_0$.
Note that the genus of the topological surface $Tb_0$
is the same as the genus of the curve $C_0$.
\end{definition}
%%%%%%%%%%%%%%%%%%%%%%%%%%%%%%%%%%%%%%%%%%%%%%%%%%%%%%%
\begin{definition}(Points, loops and paths on $Tb_u$)
We lift the points $Q$, $R_{jl}$ and $R_{1jl}$ and the loops and
paths $\gamma_{jl}$, $\gamma_{1jl}$, $\sigma_{jl}'$,
$\sigma_{1jl}'$, $\sigma_{jl}$, $\sigma_{1jl}$, $\a_i$ and $\b_i$
from $Tb_0$ to $Tb_u$, using the homotopy $h$ from Corollary 2.8.
Let
$$\begin{tabular}{llll}
$\tilde{Q}=h(u,Q)$\\
\\
$\tilde{R}_{jl}=h(u,R_{jl}),$&$\tilde{R}_{1jl}=h(u,R_{1jl}),$\\
\\
$\tilde{\gamma}_{jl}=h(u,\gamma_{jl}),$&$\tilde{\gamma}_{1jl}=h(u,\gamma_{1jl}),$\\
\\
$\tilde{\sigma}_{jl}'=h(u,\sigma_{jl}'),$&$\tilde{\sigma}_{1jl}'=h(u,\sigma_{1jl}'),$\\
\\
$\tilde{\sigma}_{jl}=h(u,\sigma_{jl}),$&$\tilde{\sigma}_{1jl}=h(u,\sigma_{1jl}),$\\
\\
$\tilde{\a}_i=h(u,\a_i),$&$\tilde{\b}_i=h(u,\b_i)$
\end{tabular}$$
on $Tb_{u}$.
\end{definition}
%%%%%%%%%%%%%%%%%%%%%%%%%%%%%%%%%%%%%%%%%%%%%%%%%%%%%%%%%%%%%%%%%%
\begin{definition} (Loops around the curve $C_0$)
For $u=\e_0 e^{2\pi i t}$, we define the loops
$$\begin{tabular}{llll}
$\tau(t)=h(\e_0 e^{2\pi i t},Q),$& starting at $\tilde{Q}$\\
\\
$\tau_{jl}(t)=h(\e_0 e^{2\pi i t},R_{jl}),$& starting at $\tilde{R}_{jl}$\\
\\
$\tau_{1jl}(t)=h(\e_0 e^{2\pi i t},R_{1jl}),$& starting at $\tilde{R}_{1jl}$
\end{tabular}$$
on $Tb$ around $C_0$.
\end{definition}
%%%%%%%%%%%%%%%%%%%%%%%%%%%%%%%%%%%%%%%%%%%%%%%%%%%%%%%%%%%%%%%%%%
\begin{definition} (Tori $T_{jl}'$, $T_{1jl}'$, $T_{jl}$, $T_{1jl}$ and $T_{i}$)
For $u=\e_0 e^{2\pi i t}$, we define the tori $T_{jl}$ and $T_{1jl}$
$$\begin{tabular}{llll}
$T_{jl}'=\{h(\e_0 e^{2\pi i t},\sigma_{jl}'(s))|0\neq s \neq 1, 0\neq t \neq 1\},$\\
\\
$T_{1jl}'=\{h(\e_0 e^{2\pi i t},\sigma_{1jl}'(s))|0\neq s \neq 1, 0\neq t \neq 1\}.$\\
\\
$T_{jl}=\{h(\e_0 e^{2\pi i t},\sigma_{jl}(s))|0\neq s \neq 1, 0\neq t \neq 1\},$\\
\\
$T_{1jl}=\{h(\e_0 e^{2\pi i t},\sigma_{1jl}(s))|0\neq s \neq 1, 0\neq t \neq 1\}.$\\
\\
$T_{i}=\{h(\e_0 e^{2\pi i t},[\a_i,\b_i](s))|0\neq s \neq 1, 0\neq t \neq 1\},$
\end{tabular}$$
where  $T_{jl}'$, $T_{1jl}'$ are tori near the points $P_{jl}$
and near $P_{1jl}$, respectively.
\end{definition}
%%%%%%%%%%%%%%%%%%%%%%%%%%%%%%%%%%%%%%%%%%%%%%%%%%%%%%%%%%%%%%%%%%
\begin{definition} ('boundary' of the Tori $T_{jl}$, $T_{1jl}$ and $T_{i}$)
By a boundary of any of the above tori,
we mean the loop on the torus,
corresponding to the boundary of the defining region.
More precisely, consider the two loops $\tilde{\sigma}_{jl}'$ and $\tau_{jl}$.
They are generators of the torus $T_{jl}$. Cut the torus $T_{jl}$
along both loops $\tilde{\sigma}_{jl}'$ and $\tau_{jl}$. We obtain a square.
By a boundary of the torus $T_{jl}$, we mean the boundary of this square.
Algebraically, this means that the boundary is a commutator. Explicitly,
$$\begin{tabular}{ll}
$\d T_{jl}'=[\tilde{\sigma}_{jl}',\tau_{jl}],\mbox{ starting at }\tilde{R}_{jl},$\\
\\
$\d T_{1jl}'=[\tilde{\sigma}_{1jl}',\tau_{jl}],\mbox{ starting at }\tilde{R}_{1jl},$\\
\\
$\d T_{jl}=[\tilde{\sigma}_{jl},\tau_{jl}],\mbox{ starting at }\tilde{Q},$\\
\\
$\d T_{1jl}=[\tilde{\sigma}_{1jl},\tau_{jl}],\mbox{ starting at }\tilde{Q},$\\
\\
$\d T_{i}=[[\tilde{\a_i},\tilde{\b_i}],\tau],\mbox{ starting at }\tilde{Q}.$
\end{tabular}$$
\end{definition}
%%%%%%%%%%%%%%%%%%%%%%%%%%%%%%%%%%%%%%%%%%%%%%%%%%%%%%%%%%%%%%

%%%%%%%%%%%%%%%%%%%%%%%%%%%%%%%%%%%%%%%%%%%%%%%%%%%%%%%%%%%%%%%%%%%%%%
\begin{definition} For fixed $u$ let $\mu$ be a loop on $Tb_u$,
starting at the point $\tilde{Q}$,
defined as a map
$$\mu:[0,1]\rightarrow Tb_u.$$
Consider the torus
$T$ as the image of $[0,1]^2$ to $X$
$$T:[0,1]^2 \rightarrow X,$$
$$T(s_0,s_1)=h(\e_0 e^{2\pi i s_0},\mu(s_1)),$$
where $h$ is the homotopy defined in Corollary 2.8.
Note that we have define two generators of the torus $T$, namely,
$\mu$ and $\tau$. Similarly to Definition 2.13,
we define a boundary of $T$, denoted by $\d T$, by a commutator
$$\d T=[\mu,\tau].$$
\end{definition}
%%%%%%%%%%%%%%%%%%%%%%%%%%%%%%%%%%%%%%%%%%%%%%%%%%%%%%%%%%%%%%
\begin{definition}
We are going to define a new homotopy of paths that
parameterizes the torus $T$ from the previous definition.
For $s\in[0,1]$ we define
$$\gamma'_s(t)=\left\{
\begin{tabular}{ll}
$(st,0)$ & $t\in[0,1]$\\
\\
$(s,t-1)$ & $t\in[1,2]$\\
\\
$(s(3-t),1)$ & $t\in[2,3]$\\
\\
$(0,4-t)$ & $t\in[3,4]$
\end{tabular}
\right.$$ Note that for fixed $s\in(0,1)$ the path $\gamma'_s$
starts at $(0,0)$ and continues along the boundary of the
rectangle with sides $s$ and $1$. In the limit $s=1$ the path
$\gamma'_1$ becomes the boundary of the square $[0,1]^2$. And for
$s=0$ the path $\gamma'_0$ becomes just one edge in direction $x_2$
with length $1$. Let $\gamma'_\bullet(s,t)=\gamma'_s(t)$. The
parametrization of the torus $T$ is given by
$$\gamma_\bullet=T\circ\gamma'_\bullet,$$
where
$\circ$ denotes a composition of functions. Define also
$$\gamma_s=T\circ\gamma'_s.$$
\end{definition}
\subsection{Construction of an abstract reciprocity law}

%%%%%%%%%%%%%%%%%%%%%%%%%%%%%%%%%%%%%%%%%%%%%%%%%%%%%%%%%%%%%%%%%
For the torus $T$, we have the following non-commutative Stokes Theorem (Theorem 1.5).
\begin{theorem}
$$\mbox{(a) }\int_{\d T}\frac{df_1}{f_1}\circ\frac{df_2}{f_2}
=\int\int_T \frac{df_1}{f_1}\wedge\frac{df_2}{f_2};$$
$$\mbox{(b) }\int_{\d T}\frac{df_1}{f_1}\circ\frac{df_2}{f_2}\circ\frac{df_3}{f_3}
=\int\int_T\frac{df_1}{f_1}\circ\left(\frac{df_2}{f_2}\wedge\frac{df_3}{f_3}\right)
+\int\int_T \left(\frac{df_1}{f_1}\wedge\frac{df_2}{f_2}\right)\circ\frac{df_3}{f_3}.$$
\end{theorem}
%%%%%%%%%%%%%%%%%%%%%%%%%%%%%%%%%%%%%%%%%%%%%%%%%%%%%%%%%%%%%%%%%%%%%%%%%%%%%%%%%%%%%%%%%
\proof The integrals on the left hands side are defined in
Definition 1.1, using the path
$$\gamma_1=T\circ \gamma'_1:t\mapsto X$$ constructed in Definition 2.15.
The integrals on the righthand side are defined in Definition 1.4,
using the homotopy of paths $\gamma_\bullet=T\circ\gamma_\bullet$,
constructed in Definition 2.15. Applying Theorem 1.5 to
$\gamma_\bullet$, we obtain the above theorem.
%%%%%%%%%%%%%%%%%%%%%%%%%%%%%%%%%%%%%%%%%%%%%%%%%%%%%%%%%%%%%%%%
\begin{definition}(relation in the fundamental group)
Let $$\delta=\left(\prod_{i=1}^g [\a_i,\b_i]\right)
\left(\prod_{j=1}^N \sigma_{jl}\right)
\left(\prod_{j=1}^{M_1} \sigma_{1jl}\right),$$
be homotopic to the trivial loop at $Q$ in $Tb_0$.
Let $\delta_{s_2}$, for $0\leq s_2 \leq 1$
be a homotopy between the loop $\delta$ and the trivial loop at $Q$,
where the homotopy is with fixed starting point $Q$.
So that for $s_2=1$ we have $$\delta_1=\delta$$ and for $s_2=0$
we have $$\delta_0\equiv Q.$$
For each value of $s_2$ and for fixed $u$ consider the loop
$$\tilde{\delta}_{s_2}=h(u,\delta_{s_2})$$
on $Tb_u$, where $h$ is the homotopy defined in Corollary 2.8.
Note that $\tilde{\delta}_{s_2}$ gives a homotopy sitting in $Tb_u$
between $\tilde{\delta}_{1}$ and the trivial loop at $\tilde{Q}$.
\end{definition}
%%%%%%%%%%%%%%%%%%%%%%%%%%%%%%%%%%%%%%%%%%%%%%%%%%%%%%%%%%%%%
\begin{theorem} Using Definition 2.17 of $\tilde{\delta}_1$,
by taking $s_2=1$,
and Definition 2.11 of $\tau$, we have
$$\int_{[\tilde{\delta}_1,\tau]}
\frac{df_1}{f_1}\circ\frac{df_2}{f_2}\circ\frac{df_3}{f_3}
=0.$$
\end{theorem}
%%%%%%%%%%%%%%%%%%%%%%%%%%%%%%%%%%%%%%%%%%%%%%%%%%%%%%%%%%%%
\begin{remark} Theorem 2.18 will be interpreted later
as the sum of the logarithmic symbols is an integer multiple of
$(2\pi i)^3$. Where are the logarithmic symbols? The logarithmic
symbols will be
$$\int_{[\tilde{\sigma}_{jl},\tau]}
\frac{df_1}{f_1}\circ\frac{df_2}{f_2}\circ\frac{df_3}{f_3}.$$
The integrals
$$\int_{[\tilde{\sigma}'_{jl},\tau]}
\frac{df_1}{f_1}\circ\frac{df_2}{f_2}\circ\frac{df_3}{f_3}=0$$
where we use $\tilde{\sigma}'_{jl}$ instead of
$\tilde{\sigma}_{jl}$. And finally, the integer multiple of $(2\pi
i)^3$ will come from
$$\int_{[[\tilde{\a}_i,\tilde{\b}_i],\tau]}
\frac{df_1}{f_1}\circ\frac{df_2}{f_2}\circ\frac{df_3}{f_3}.$$
Actually there is one other source of integer multiples of $(2\pi
i)^3$, which will be called 'extra terms'.
\end{remark}
%%%%%%%%%%%%%%%%%%%%%%%%%%%%%%%%%%%%%%%%%%%%%%%%%%%%%%%%%%%%
We are going to prove this theorem using Theorem 1.6
and Non-commutative Stokes Theorem
for iterated integrals stated in Theorem 1.9.
Before doing that we have to make several definitions.
%%%%%%%%%%%%%%%%%%%%%%%%%%%%%%%%%%%%%%%%%%%%%%%%%%%%%%%%%%%%%%%%
\begin{definition} (For $3$-dimensional non-commutative Stokes Theorem)
For each value of $s_2$ and for fixed $u$ consider the loop
$$\tilde{\delta}_{s_2}=h(u,\delta_{s_2})$$
on $Tb_u$, where $h$ is the homotopy defined in Corollary 2.8. Now
use Definition 1.14. Instead of the loop $\mu$ use
$\tilde{\delta}_{s_2}$, in order to define a torus. Denote the
corresponding torus by $T_{s_2}$. Now use Definition 1.15 to
define parametrization of the torus $T_{s_2}$. For fixed $s_2$,
let
$$\gamma_{\bullet,s_2}=T_{s_2}\circ\gamma'_\bullet.$$
Similarly, for fixed $s$ and $s_2$ let
$$\gamma_{s,s_2}=T_{s_2}\circ\gamma'_s.$$
Now instead of the variable $s$ write the variable $s_1$.
When $s_1$ denote
$$\gamma_{s_1,\bullet}(s_2,t)=\gamma_{s_1,s_2}(t).$$
Let also
$$\gamma_{\bullet,\bullet}(s_1,s_2,t)=\gamma_{s_1,s_2}(t).$$
\end{definition}
%%%%%%%%%%%%%%%%%%%%%%%%%%%%%%%%%%%%%%%%%%%%%%%%%%%%%%%%%%%%%%%%
\proof (of Theorem 1.18) In the notation
of Definition 1.19, the domain of integration is
$\d T_1$. It is enough to show that
$$\int\int_{T_1}
\frac{df_1}{f_1}\circ
\left(\frac{df_2}{f_2}\wedge\frac{df_3}{f_3}\right)
=0$$
and
$$\int\int_{T_1}
\left(\frac{df_1}{f_1}\wedge\frac{df_2}{f_2}\right)
\circ\frac{df_3}{f_3}
=0.$$
Using the notation from Definition 1.19 and from subsection 1.7, we have
$$\int\int_{T_1}
\frac{df_1}{f_1}\circ
\left(\frac{df_2}{f_2}\wedge\frac{df_3}{f_3}\right)=
\int\int_{\gamma_{\bullet,1}}
\frac{df_1}{f_1}\circ
\left(\frac{df_2}{f_2}\wedge\frac{df_3}{f_3}\right).
$$
From Theorem 1.9, we have
$$2\int\int\int_{\gamma_{\bullet,\bullet}}
\frac{df_1}{f_1}\wedge\frac{df_2}{f_2}\wedge\frac{df_3}{f_3}=
\left(
\int\int_{\gamma_{0,\bullet}}
-\int\int_{\gamma_{1,\bullet}}
+\int\int_{\gamma_{\bullet,0}}
-\int\int_{\gamma_{\bullet,1}}
\right)
\frac{df_1}{f_1}\circ
\left(\frac{df_2}{f_2}\wedge\frac{df_3}{f_3}\right).
$$
We have identified the last integral. Note also that
$$\frac{df_1}{f_1}\wedge\frac{df_2}{f_2}\wedge\frac{df_3}{f_3}=0,$$
because $X$ is a $2$-dimensional variety.
The domain $\gamma_{\bullet,0}$ comes from $\tilde{\delta}_{s_2}$ for $s_2=0$,
which is the trivial loop at $\tilde{Q}$.
Then $\gamma_{\bullet,0}$ is a $1$-dimensional region. Therefore,
$$\int\int_{\gamma_{\bullet,0}}
\frac{df_1}{f_1}\circ
\left(\frac{df_2}{f_2}\wedge\frac{df_3}{f_3}\right)
=0.$$
For fixed $s_1$ the domain of integration
$\gamma_{s_1,\bullet}$ lies inside $Tb_u$. Since $Tb_u$ is a subset of
the algebraic curve on $X$ given by the equation $x_0=u$, we obtain that
$$\frac{df_2}{f_2}\wedge\frac{df_3}{f_3}=0.$$
Therefore,
$$\int\int_{\gamma_{0,\bullet}}
\frac{df_1}{f_1}\circ
\left(\frac{df_2}{f_2}\wedge\frac{df_3}{f_3}\right)=0$$
and
$$\int\int_{\gamma_{1,\bullet}}
\frac{df_1}{f_1}\circ
\left(\frac{df_2}{f_2}\wedge\frac{df_3}{f_3}\right)=0.$$
Thus, the integral
$$\int\int_{\gamma_{\bullet,1}}
\frac{df_1}{f_1}\circ
\left(\frac{df_2}{f_2}\wedge\frac{df_3}{f_3}\right)=0.
$$
Similarly, we can prove that
$$\int\int_{T_1}
\left(\frac{df_1}{f_1}\wedge\frac{df_2}{f_2}\right)
\circ\frac{df_3}{f_3} =\int\int_{\gamma_{\bullet,1}}
\left(\frac{df_1}{f_1}\wedge\frac{df_2}{f_2}\right)
\circ\frac{df_3}{f_3} =0.$$ Now we use Theorem 2.16, and we obtain
$$\int\int_{\d T_1}
\frac{df_1}{f_1}\circ\frac{df_2}{f_2}\circ\frac{df_3}{f_3}=0.$$
%%%%%%%%%%%%%%%%%%%%%%%%%%%%%%%%%%%%%%%%%%%%%%%%%%%%%%%%%%%%%%%%%
\begin{definition} ($\mu_i$)
Define loops $\mu_i$ so that $\mu_i$ is either
$[\tilde{\a}_j,\tilde{\b}_j]$ or $\tilde{\sigma}_{jl}$, or
$\tilde{\sigma}_{1jl}$. Let $K$ be the number of all the above
different loops. The integer $K$ is chosen so that for the ordered
product we have
$$\prod_{i=1}^K \mu_i=\tilde{\delta}_1.$$
\end{definition}
%%%%%%%%%%%%%%%%%%%%%%%%%%%%%%%%%%%%%%%%%%%%%%%%%%%%%%%%%%%%%%%%%
\begin{definition} ($\pi_j$)
Define loops $\pi_i$ by
$$\pi_j=\prod_{i=1}^j \mu_i.$$
\end{definition}
%%%%%%%%%%%%%%%%%%%%%%%%%%%%%%%%%%%%%%%%%%%%%%%%%%%%%%%%%%%%%%%%%
Note that $\pi_1=\mu_1$ and $\pi_K=\tilde{\delta}_1$,
which is homotopic to the trivial loop at $\tilde{Q}$.
%%%%%%%%%%%%%%%%%%%%%%%%%%%%%%%%%%%%%%%%%%%%%%%%%%%%%%%%%%%
\begin{lemma}
(a) $[\pi_{j+1},\tau]=\pi_j[\mu_{j+1},\tau]\pi_j^{-1}[\pi_{j},\tau].$

(b) $[\pi_{K},\tau]=\prod_{i=1}^K(\pi_{K-i}[\mu_{K-i+1},\tau]\pi_{K-i}^{-1}).$
\end{lemma}
%%%%%%%%%%%%%%%%%%%%%%%%%%%%%%%%%%%%%%%%%%%%%%%%%%%%%%%%%%
\proof It follows by direct computation of commutators.
%%%%%%%%%%%%%%%%%%%%%%%%%%%%%%%%%%%%%%%%%%%%%%%%%%%%%%%%%%%
\begin{lemma}
$$0=\int_{[\pi_{K},\tau]}
\frac{df_1}{f_1}\circ\frac{df_2}{f_2}\circ\frac{df_3}{f_3}
=\sum_{i=1}^K\int_{\pi_{K-i}[\mu_{K-i+1},\tau]\pi_{K-i}^{-1}}
\frac{df_1}{f_1}\circ\frac{df_2}{f_2}\circ\frac{df_3}{f_3}.$$
\end{lemma}
%%%%%%%%%%%%%%%%%%%%%%%%%%%%%%%%%%%%%%%%%%%%%%%%%%%%%%%%%%%%%%%%%%%%%%
\proof The first equality is Theorem 2.18. For the second
equality, we use composition of paths in iterated integrals stated
in Corollary 1.3. By induction, we can decompose the path
$[\pi_{K},\tau]$ into two paths, then into three paths and finally
into $K$ paths $\pi_{K-i}[\mu_{K-i+1},\tau]\pi_{K-i}^{-1}$ for
$i=1,\dots,K$. Using Corollary 1.3, we can expand the iterated
integral over $[\pi_{K},\tau]$ in terms of sum of products of
iterated integrals over
$\pi_{K-i}[\mu_{K-i+1},\tau]\pi_{K-i}^{-1}$ for $i=1,\dots,K$.
Note that
$$\int_{\pi_{K-i}[\mu_{K-i+1},\tau]\pi_{K-i}^{-1}}\omega_1=0$$
for any differential $1$-form $\omega$.
Therefore, when we have a product of iterated integrals in the expansion,
we must have at least one integral of one differential $1$-form.
That gives zero contribution.
Therefore, there are no products in the expansion over
the smaller paths,
but only a sum of triple iterations
over each of the smaller paths
$\pi_{K-i}[\mu_{K-i+1},\tau]\pi_{K-i}^{-1}$.
This proves the lemma.
%%%%%%%%%%%%%%%%%%%%%%%%%%%%%%%%%%%%%%%%%%%%%%%%%%%%%%%%%%%%%%%%%%%%
Now we are going to examine each of the components
$\int_{\pi_{i}[\mu_{i+1},\tau]\pi_{i}^{-1}}
\frac{df_1}{f_1}\circ\frac{df_2}{f_2}\circ\frac{df_3}{f_3}.$
%%%%%%%%%%%%%%%%%%%%%%%%%%%%%%%%%%%%%%%%%%%%%%%%%%%%%%%%%%%
\begin{lemma}
$$
\begin{tabular}{ll}
$\int_{\pi_{i}[\mu_{i+1},\tau]\pi_{i}^{-1}}
\frac{df_1}{f_1}\circ\frac{df_2}{f_2}\circ\frac{df_3}{f_3}=$
&$\int_{[\mu_{i+1},\tau]}
\frac{df_1}{f_1}\circ\frac{df_2}{f_2}\circ\frac{df_3}{f_3}$\\
\\
&$+\int_{\pi_{i}}\frac{df_1}{f_1}\int_{[\mu_{i+1},\tau]}
\frac{df_2}{f_2}\circ\frac{df_3}{f_3}
+\int_{[\mu_{i+1},\tau]}
\frac{df_1}{f_1}\circ\frac{df_2}{f_2}
\int_{\pi_{i}^{-1}}\frac{df_3}{f_3}
.$
\end{tabular}$$
\end{lemma}
%%%%%%%%%%%%%%%%%%%%%%%%%%%%%%%%%%%%%%%%%%%%%%%%%%%%%%%%%%%%%%
\proof Decompose the path $\pi_{i}[\mu_{i+1},\tau]\pi_{i}^{-1}$ into
$\pi_{i}$, $[\mu_{i+1},\tau]$ and $\pi_{i}^{-1}$. Then use
Corollary 1.3. Finally, use that
$$\int_{[\mu_{i+1},\tau]}\omega=0$$
for any differential $1$-form $\omega.$
%%%%%%%%%%%%%%%%%%%%%%%%%%%%%%%%%%%%%%%%%%%%%%%%%%%%%%%%%%%
\begin{lemma} (a) For $i\leq g$, where $g$ is the genus of $C_0$, we have
$\mu_i=[\tilde{\a}_i,\tilde{\a}_i]$. Then
$$\int_{\mu_{i}}\omega=0$$
and consequently,
$$\int_{\pi_{i}}\omega=0$$
for a $1$ form $\omega$.
(b) When $\mu_i=\tilde{\sigma}'_{jl}$, we have
$$\int_{[\mu_{i},\tau]}
\frac{df_1}{f_1}\circ\frac{df_2}{f_2}=0$$
and
$$\int_{[\mu_{i},\tau]}
\frac{df_1}{f_1}\circ\frac{df_2}{f_2}\circ\frac{df_3}{f_3}=0.$$
\end{lemma}
%%%%%%%%%%%%%%%%%%%%%%%%%%%%%%%%%%%%%%%%%%%%%%%%%%%%%%%%%%%%%%
\proof Part (a) is straight forward. Part (b) will be proven in
subsection 3.1.
%%%%%%%%%%%%%%%%%%%%%%%%%%%%%%%%%%%%%%%%%%%%%%%%%%%%%%%%%%%%%%%%%
\begin{definition}
When $\mu_i=\tilde{\sigma}_{jl}$, we define a new logarithmic symbol
as
$$Log[f_1,f_2,f_3]_{C_0,P_{jl}}^{x_0,Q}=
\lim_{\e\rightarrow 0}\int_{[\mu_{i},\tau]}
\frac{df_1}{f_1}\circ\frac{df_2}{f_2}\circ\frac{df_3}{f_3}.$$
\end{definition}
%%%%%%%%%%%%%%%%%%%%%%%%%%%%%%%%%%%%%%%%%%%%
\begin{definition}
For $i\leq g$, where $g$ is the genus of $C_0$, we have
$\mu_i=[\tilde{\a}_i,\tilde{\b}_i]$. Then we call the integral
$$\int_{[\mu_{i},\tau]}
\frac{df_1}{f_1}\circ\frac{df_2}{f_2}\circ\frac{df_3}{f_3}.$$
the {\bf{commutator}} $[\a,\b]$-{\bf{terms}}
\end{definition}
%%%%%%%%%%%%%%%%%%%%%%%%%%%%%%%%%%%%%%%%%%%%%%%%%%%%%%%%%%%%%%%%%
\begin{definition}
Consider the terms
$$\int_{\pi_{i}}\frac{df_1}{f_1}\int_{[\mu_{i+1},\tau]}
\frac{df_2}{f_2}\circ\frac{df_3}{f_3}$$
and
$$\int_{[\mu_{i+1},\tau]}
\frac{df_1}{f_1}\circ\frac{df_2}{f_2}
\int_{\pi_{i}^{-1}}\frac{df_3}{f_3}$$
in Lemma 2.26. The only case when these terms do not vanish is when
$i>g$ and $\mu_{i+1}=\tilde{\sigma}_{jl}$.
In this case we call the above two integrals
{\bf{extra terms}}.
\end{definition}
%%%%%%%%%%%%%%%%%%%%%%%%%%%%%%%%%%%%%%%%%%%%%%%%%%%%%%%%%%%%%%%%%
\begin{remark}
We are going to show that a commutator $[\a,\b]$-term gives an
integer multiple of $(2\pi i)^3$ in Lemma 3.25 and also
that the extra terms give an integer multiple of $(2\pi i)^3$ in
Theorem 3.30.
\end{remark}
%%%%%%%%%%%%%%%%%%%%%%%%%%%%%%%%%%%%%%%%%%%%%%%%%%%%%%%%%%%%%%%%%
\begin{theorem}(Abstract reciprocity law)
$$\sum_{P_{jl}}Log[f_1,f_2,f_3]_{C_0,P_{jl}}^{x_0,Q}
=-\sum (\mbox{commutator terms})-\sum(\mbox{extra terms}).$$
Note that from the previous remark
the right hand side of the above equation
is an integer multiple of $(2\pi i)^3$.
\end{theorem}
%%%%%%%%%%%%%%%%%%%%%%%%%%%%%%%%%%%%%%%%%%%%%%%%%%%%%%%
\proof For $i$ such that $\mu_i=\tilde{\sigma}_{jl}$,
by Lemma 2.25,
we have that the sum of the Log-symbols
and the extra terms gives the sum of
$$\int_{\pi_{i}[\mu_{i+1},\tau]\pi_{i}^{-1}}
\frac{df_1}{f_1}\circ\frac{df_2}{f_2}\circ\frac{df_3}{f_3},$$
over the above type of indexes $i$.
For $i\leq g$ we have $\mu_i=[\tilde{\a}_i,\tilde{\a}_i]$. Then
by Lemma 2.6 (a), we have
$$\int_{[\mu_{i},\tau]}
\frac{df_1}{f_1}\circ\frac{df_2}{f_2}\circ\frac{df_3}{f_3}
=
\int_{\pi_{i-1}[\mu_{i},\tau]\pi_{i-1}^{-1}}
\frac{df_1}{f_1}\circ\frac{df_2}{f_2}\circ\frac{df_3}{f_3}.$$
For $i$ such that $\mu_i=\tilde{\sigma}'_{jl}$,
by Lemma 2.26(b), we have that
$$\int_{\pi_{i-1}[\mu_{i},\tau]\pi_{i-1}^{-1}}
\frac{df_1}{f_1}\circ\frac{df_2}{f_2}\circ\frac{df_3}{f_3}=0.$$ We
have to show that the sum of the above three integrals is zero. By
Lemma 2.24, we have that this sum is equal to
$$\int_{[\pi_{K},\tau]}
\frac{df_1}{f_1}\circ\frac{df_2}{f_2}\circ\frac{df_3}{f_3},$$
where $\pi_K=\tilde{\delta}_1$. Now from Theorem 2.18,
the above integral is zero.
%%%%%%%%%%%%%%%%%%%%%%%%%%%%%%%%%%%%%%%%%%%%%%%%%%%%%%%%%%%%%%%%%%%%%%
%%%%%%%%%%%%%%%%%%%%%%%%%%%%%%%%%%%%%%%%%%%%%%%%%%%%%%%%%%%%%%%%%%%%%%%%
\section{Construction of a new symbol and a reciprocity law}
\subsection{Local properties. Integrating over the torus $T'_{jl}$}

In this subsection, we use the notation from subsection 2.4.
Consider the torus $T'_{jl}$ from Definition 2.12. We use
Definition 2.14 and 2.15 in order to define iterated integrals
over $T'_{jl}$ and over $\d T'_{jl}$. Let $x_i$ be as in
Definition $2.1$. Then near the point $P_{jl}$ we have local
coordinates $x_0$ and $x_j$. For $k=1,2,3$ let
$$f_k=x_0^{n_{k0}}x_j^{n_{kj}}g_k,$$
where $g_k$ are holomorphic functions near $(x_0,x_j)=(0,0)$.
We are going to simplify the notation.
%%%%%%%%%%%%%%%%%%%%%%%%%%%%%%%%%%%%%%%%%%%
\begin{definition}(Integers $m_k$ and $n_k$, function $x$, $y$ and $g_k$ and torus $T_0$)

We are going to use the integers
$$m_k\mbox{ instead of }n_{k0},$$
$$n_k\mbox{ instead of }n_{kj},$$
and the functions
$$x\mbox{ instead of }x_0,$$
$$y\mbox{ instead of }x_j,$$
(see Definition 2.1 for $n_{k0}$ and for $n_{kj}$).
Also we are going to use
$$T_0\mbox{ instead of }T'_{jl},$$
(see Definition 2.12 for $T'_{jl}$).
We define $g_k$ for $k=1,2,3$ by
$$f_k=x^{m_k}y^{n_k}g_k.$$
Note that $g_k$ is non-zero holomorphic function near
$(x,y)=(0,0)$.
\end{definition}
%%%%%%%%%%%%%%%%%%%%%%%%%%%%%%%%%%%%%%%%%%%%%%%%%
\begin{lemma}
(a) $\int\int_{\d T_0}
\frac{df_1}{f_1}\circ\frac{df_2}{f_2}=\int\int_{T_0}
\frac{df_1}{f_1}\wedge\frac{df_2}{f_2}=-(m_1n_2-m_2n_1)(2\pi i)^2$;

(b) $\int\int_{T_0}
\frac{df_1}{f_1}\circ\left(\frac{df_2}{f_2}\wedge\frac{df_3}{f_3}\right)
=-(m_1+n_1)(m_2n_3-m_3n_2)\frac{(2\pi i)^3}{2}+O(\e)$;

(c) $\int\int_{T_0}
\left(\frac{df_1}{f_1}\wedge\frac{df_2}{f_2}\right)\circ\frac{df_3}{f_3}
=(m_3+n_3)(m_1n_2-m_2n_1)\frac{(2\pi i)^3}{2}+O(\e).$
\end{lemma}
\proof Denote by $O(x,y)$ all the terms of degree at least $1$ in
the variable $x$ or $y$. Denote by $O(1)$ all the terms of degree
at least $0$ or higher. Note that on $T_0$ we have $|x|<\e$ and
$|y|<\e$. For part (a), we have
$$\frac{df_1}{f_1}=m_1\frac{dx}{x}+n_1\frac{dy}{y}+\frac{dg_1}{g_1}=
m_1\frac{dx}{x}+n_1\frac{dy}{y}+dO(x,y),$$ since $g_1$
has no zeroes or poles at $x=0$ or $y=0$. Similarly,
$$\frac{df_2}{f_2}=m_2\frac{dx}{x}+n_2\frac{dy}{y}+dO(x,y).$$
Then
$$\frac{df_1}{f_1}\wedge \frac{df_2}{f_2}
=(m_1n_2-m_2n_1)\frac{dx}{x}\wedge\frac{dy}{y}+
\frac{dx}{x}\wedge dO(x,y)+dO(x,y)\wedge
dO(x,y)+ dO(x,y)\wedge\frac{dy}{y}.$$ Since $x$ and
$y$ vary over loops, we only pick the residues, when we compute
the integral. Therefore,
$$\int_{T_0}\frac{df_1}{f_1}\wedge \frac{df_2}{f_2}
=(m_1n_2-m_2n_1)\int_{T_0}\frac{dx_0}{x_0}\wedge\frac{dx_j}{x_j}
=-(m_1n_2-m_2n_1)(2\pi i)^2.$$

For part (b) denote by
$\gamma_\bullet(s,t)=T_0\circ\gamma'_s(t),$ as in Definition 2.15.
Note
that for $t\notin [1,2]$, the image of $\gamma_\bullet$ is one
dimensional. In order to have a non-zero $2$-dimensional iterated
integral, we must have the $2$-form integrated over
$\gamma_\bullet$ restricted to $(s,t)\in [0,1]\times [1,2]$. All
other intervals for $t$ will have no contribution, since
$\gamma_\bullet$ has $1$-dimensional image. Then the first
$1$-form in the $2$-dimensional iterated integral, must be
integrated over $\gamma_s$ up to a point $(s,t)\in [0,1]\times
[1,2]$. Consider the image of $\gamma_\bullet$ in the
$(x_0,x_j)$-coordinate system.
Then we have to integrate the first $1$-form up to
$(x'_0,x'_j)$ along $\gamma_s$ for some $s$. We obtain
$$\int_{T_0}
\frac{df_1}{f_1}\circ\left(\frac{df_2}{f_2}\wedge\frac{df_3}{f_3}\right)
=\int\int_{[0,1]^2}
\left(\int_{(x_0',x_j')=(0,0)}^{(x_1,x_2)}\frac{df_1}{f_1}(x_1',x_2')\right)
\left(\frac{df_2}{f_2}\wedge\frac{df_3}{f_3}\right)(x_1,x_2).$$
We have
$$\int_{(x_1',x_2')=(0,0)}^{(x_1,x_2)}\frac{df_1}{f_1}(x_1',x_2')
=m_1\frac{dx_0}{x_0}+n_1\frac{dx_j}{x_j}+O(\e).$$ Therefore,
$$\begin{tabular}{ll}
$\int\int_{T_0}
\frac{df_1}{f_1}\circ\left(\frac{df_2}{f_2}\wedge\frac{df_3}{f_3}\right)=$\\
\\
$=\int_{T_0} \left(m_1\frac{dx_0}{x_0}+n_1\frac{dx_j}{x_j}\right)
\circ\left((m_2n_3-m_3n_2)\frac{dx_0}{x_0}\wedge\frac{dx_j}{x_j}\right)+O(\e)=$\\
\\
$=-(m_1+n_1)(m_2n_3-m_3n_2)\frac{(2\pi i)^3}{2}+O(\e).$
\end{tabular}$$

Part (c) can be proven in a similar way as part (b).
%%%%%%%%%%%%%%%%%%%%%%%%%%%%%%%%%%%%%%%%%%%%%%%%%
\begin{corollary} For the torus ${\d T'_{1jl}}$, we have
(a) $\int\int_{\d T'_{1jl}}
\frac{df_1}{f_1}\circ\frac{df_2}{f_2}=\int\int_{T_0}
\frac{df_1}{f_1}\wedge\frac{df_2}{f_2}=0;$

(b) $\int\int_{\d T'_{1jl}}
\frac{df_1}{f_1}\circ\left(\frac{df_2}{f_2}\wedge\frac{df_3}{f_3}\right)
=O(\e)$;

(c) $\int\int_{\d T'_{1jl}}
\left(\frac{df_1}{f_1}\wedge\frac{df_2}{f_2}\right)\circ\frac{df_3}{f_3}
=O(\e).$
\end{corollary}
\proof We are going to use the curves $D_{ij}$ from Definition
2.1. Since the differential forms $\frac{df_k}{f_k}$ have no
residue along $D_{ij}$, we obtain that the corollary as a
consequence of Lemma 3.2.
%%%%%%%%%%%%%%%%%%%%%%%%%%%%%%%%%%%%%%%%%%%%%%%%%
%For part (c) we use the construction from the proof of part(b).
%Note that the last $1$-form in the $2$-dimensional iterated
%integral, must be integrated over $\gamma_s$ from a point
%$(s,t)\in [0,1]\times [1,2]$ to $(0,0)$. Consider the image of
%$\gamma_\bullet$ in the $(x_1,x_2)$-coordinate system, where
%$$z_j=\e_j e^{2\pi i x_j}$$
%for $j=1,2$. Then we have to integrate the first $1$-form from
%$(x_1,x_2)$ to $(0,0)$. We obtain
%$$\begin{tabular}{ll}
%$\int\int_{T_0}
%\left(\frac{df_1}{f_1}\wedge\frac{df_2}{f_2}\right)\circ\frac{df_3}{f_3}=$\\
%\\
%$=\int_{(x_1,x_2)=(1,1)}^{(0,0)}
%\left(\int\int_{[0,x_1]\times[0,x_2]}
%T^*\left(\frac{df_2}{f_2}\wedge\frac{df_3}{f_3}\right)(x_1',x_2')\right)
%T^*\frac{df_1}{f_1}(x_1,x_2).$
%\end{tabular}$$
%We have
%$$\int_{(x_1,x_2)=(x_1',x_2')}^{(0,0)}T^*\frac{df_1}{f_1}(x_1,x_2)
%=-i_1dx_1-j_1dx_2+O(\e).$$ Therefore,
%$$\begin{tabular}{ll}
%$\int\int_{T_0}
%\left(\frac{df_1}{f_1}\wedge\frac{df_2}{f_2}\right)\circ\frac{df_3}{f_3}=$\\
%\\
%$=\int\int_{T_0}
%\left((i_1j_2-i_2j_1)\frac{dz_1}{z_1}\wedge\frac{dz_2}{z_2}\right)
%\circ\left(-i_3\frac{dz_1}{z_1}-j_3\frac{dz_2}{z_2}\right) +O(\e)=$\\
%\\
%$=-(i_3+j_3)(i_1j_2-i_2j_1)\frac{(2\pi i)^3}{2}+O(\e).$
%\end{tabular}$
%%%%%%%%%%%%%%%%%%%%%%%%%%%%%%%%%%%%%%%%%%%%%%%%%%%%%%%%%%%%%%%%%%%
%%%%%%%%%%%%%%%%%%%%%%%%%%%%%%%%%%%%%%%%%%%%%%%%%%%%%%%%%%%%%%%%%%%%
\subsection{Semi-local behavior. Integrating over the torus $T_{jl}$.}
We are going to use the construction of the domain of integration
from the subsection 2.4.
%%%%%%%%%%%%%%%%%%%%%%%%%%%%%%%%%%%%%%%%%%%%%%%%%%%%%%%%%%%%%%%%%
\begin{definition} (Paths $\gamma_0$, $\gamma$, loops $\sigma_0$, $\sigma$, $\tau_0$)
We are going to use
$$\gamma_0\mbox{ instead of }\gamma_{jl},$$
$$\gamma\mbox{ instead of }\tilde{\gamma}_{jl},$$
$$\sigma_0\mbox{ instead of }\tilde{\sigma}'_{jl},$$
$$\sigma\mbox{ instead of }\tilde{\sigma}_{jl},$$
$$\tau_0\mbox{ instead of }\tau_{jl}.$$
See Definition 2.9  for $\gamma_{jl}$. See  Definition 2.10
for $\tilde{\gamma}_{jl},$ $\tilde{\sigma}'_{jl},$
and $\tilde{\sigma}_{jl}.$ See  Definition 2.11 for
$\tau_{jl}.$
\end{definition}
%%%%%%%%%%%%%%%%%%%%%%%%%%%%%%%%%%%%%%%%%%%%%%%
\begin{definition} ($\tau'$)
Let $\tau'$ be a loop on $Tb_u$ starting at $\tilde{Q},$ defined
in the following way. First define a constant $C$ by
$$(x,y)(\tilde{Q})=(\e_1,C).$$
For $t\in[0,1]$ define the loop $\tau'$ so that
$$\tau'(0)=\tilde{Q}$$
and
$$\tau'(t)\subset (x,y)^{-1}(\e_1e^{2\pi i t},C).$$
Note that $\tau'$ does not vary in direction of $y$, while
the loop $\tau$ might vary in that direction.
See Definition 2.11 for $\tau$.
\end{definition}
%%%%%%%%%%%%%%%%%%%%%%%%%%%%%%%%%%%%%%%%%%%%%%%%%%%%%%
\begin{definition} ($T$ and $T'$)
We are going to use
$$T\mbox{ instead of }T_{jl}.$$
Let $T'$ denote a torus with generators
$\tau'$ and $\sigma$.
\end{definition}
%%%%%%%%%%%%%%%%%%%%%%%%%%%%%%%%%%%%%%%%%%%%%%%%%%%%%%
\begin{definition} ($\lambda^t$)
For fixed $t$, define also
$\Lambda^t$ to be a
segment of the straight line, in the complex plane
$x=\e_1 e^{2\pi i t}$ in the coordinate system $(x,y)$,
joining $(x,y)(\tau(t))$
with $(x,y)(\tau_j(t))$.
Let $$\lambda^t\subset(x,y)^{-1}L^t$$
be the corresponding path, joining $\tau(t)$
with $\tau'(t)$.
\end{definition}
%%%%%%%%%%%%%%%%%%%%%%%%%%%%%%%%%%%%%%%%%%%%%%%%%
\begin{proposition}
The following two iterated integrals are equal in the limit
$\e\rightarrow 0$. That is,
$$\lim_{\e\rightarrow 0}\int_{[\sigma,\tau']}
\frac{dx_{i_1}}{x_{i_1}}\circ\frac{dx_{i_2}}{x_{i_2}}\circ\frac{dx_{i_3}}{x_{i_3}}
-\int_{[\sigma,\tau]}
\frac{dx_{i_1}}{x_{i_1}}\circ\frac{dx_{i_2}}{x_{i_2}}\circ\frac{dx_{i_3}}{x_{i_3}}
=0,$$
where the commutator $[\sigma,\tau']$ in the
first integral is the boundary of $T'$
and the commutator $[\sigma,\tau]$ in the
second integral is the boundary of $T$.
\end{proposition}
%%%%%%%%%%%%%%%%%%%%%%%%%%%%%%%%%%%%%%%%%%%%%%%%%%%%%%%%%%%%%%%%%
\begin{remark} It is easier to do explicit computations, using the first
integral, involving $\tau'$.
\end{remark}
%%%%%%%%%%%%%%%%%%%%%%%%%%%%%%%%%%%%%%%%%%%%%%%%%%%%%%%%
Recall that $\tau_0$ and $\sigma_0$ are the two
loops on $T_0$ (see Definitions 2.13 and 3.3).
%%%%%%%%%%%%%%%%%%%%%%%%%%%%%%%%%%%%%%%%%%%%%%%%%%%
\begin{definition} ($S$, $S^{-1}$, $S'$ and $S'^{-1}$)
Let $S$ be the region bounded by
$$\d S=\gamma\tau_0\gamma^{-1}\tau^{-1}.$$
Let
$S^{-1}$ be the region bounded by
$$\d S^{-1}=\gamma^{-1}\tau\gamma\tau^{-1}_0.$$
Let $S'$ be the region bounded by
$$\d S'=\gamma\tau_0\gamma^{-1}\tau'^{-1}.$$
Let also
$S'^{-1}$ be the region bounded by
$$\d S'^{-1}=\gamma^{-1}\tau_j\gamma\tau_0^{-1}.$$
\end{definition}
%%%%%%%%%%%%%%%%%%%%%%%%%%%%%%%%%%%%%%%%%%%%%%%%%
\begin{lemma}
$$\begin{tabular}{llll}
(a) & $\int_{[\sigma,\tau']}
\frac{dx_{i_1}}{x_{i_1}}\circ\frac{dx_{i_2}}{x_{i_2}}\circ\frac{dx_{i_3}}{x_{i_3}}$
&
$=\int_{\d T_0}
\frac{dx_{i_1}}{x_{i_1}}\circ\frac{dx_{i_2}}{x_{i_2}}\circ\frac{dx_{i_3}}{x_{i_3}}+$\\
\\
&&$+\int_{\sigma_0}\frac{dx_{i_1}}{x_{i_1}}
\int\int_{S'^{-1}}\frac{dx_{i_2}}{x_{i_2}}\wedge\frac{dx_{i_3}}{x_{i_3}}+$\\
\\
&&$+\int\int_{S'^{-1}}\frac{dx_{i_1}}{x_{i_1}}\wedge\frac{dx_{i_2}}{x_{i_2}}
\int_{\sigma^{-1}_0}\frac{dx_{i_3}}{x_{i_3}};$\\
\\
(b) & $\int_{[\sigma,\tau]}
\frac{dx_{i_1}}{x_{i_1}}\circ\frac{dx_{i_2}}{x_{i_2}}\circ\frac{dx_{i_3}}{x_{i_3}}$
&
$=\int_{\d T_0}
\frac{dx_{i_1}}{x_{i_1}}\circ\frac{dx_{i_2}}{x_{i_2}}\circ\frac{dx_{i_3}}{x_{i_3}}+$\\
\\
&&$+\int_{\sigma_0}\frac{dx_{i_1}}{x_{i_1}}
\int\int_{S^{-1}}\frac{dx_{i_2}}{x_{i_2}}\wedge\frac{dx_{i_3}}{x_{i_3}}+$\\
\\
&&$+\int\int_{S^{-1}}\frac{dx_{i_1}}{x_{i_1}}\wedge\frac{dx_{i_2}}{x_{i_2}}
\int_{\sigma^{-1}_0}\frac{dx_{i_3}}{x_{i_3}}.$
\end{tabular}$$
\end{lemma}
%%%%%%%%%%%%%%%%%%%%%%%%%%%%%%%%%%%%%%%%%%%%%%%%%%
\proof Consider the domain
$$S\cup\gamma S^{-1}\gamma^{-1}.$$
Note that $S$ and $\gamma S^{-1}\gamma^{-1}$ parameterize the same
domain but with opposite orientations. Then
$$\int_{S\cup\tilde{\gamma}_{jl}S^{-1}\tilde{\gamma}_{jl}^{-1}}
\frac{dx_{i_1}}{x_{i_1}}\circ
\left(\frac{dx_{i_2}}{x_{i_2}}\wedge\frac{dx_{i_3}}{x_{i_3}}\right)
=0.$$
and
$$\int_{S\cup\gamma S^{-1}\gamma^{-1}}
\left(\frac{dx_{i_1}}{x_{i_1}}\wedge\frac{dx_{i_2}}{x_{i_2}}\right)
\circ\frac{dx_{i_3}}{x_{i_3}}
=0.$$

From Lemma 1.6 for composition of domains of $2$-dimensional
iterated integrals for the domain $T$, we have

$$\begin{tabular}{ll}
$\int\int_{T}\frac{dx_{i_1}}{x_{i_1}}\circ
\left(\frac{dx_{i_2}}{x_{i_2}}\wedge\frac{dx_{i_3}}{x_{i_3}}\right)$&
$=\int\int_{S}\frac{dx_{i_1}}{x_{i_1}}\circ
\left(\frac{dx_{i_2}}{x_{i_2}}\wedge\frac{dx_{i_3}}{x_{i_3}}\right)

+\int\int_{S^{-1}}\frac{dx_{i_1}}{x_{i_1}}\circ
\left(\frac{dx_{i_2}}{x_{i_2}}\wedge\frac{dx_{i_3}}{x_{i_3}}\right)+$\\
\\
&$+\int_{\tilde{\gamma}_{jl}}
\frac{dx_{i_1}}{x_{i_1}}
\int\int_{S^{-1}}\frac{dx_{i_2}}{x_{i_2}}\wedge\frac{dx_{i_3}}{x_{i_3}}+$\\
\\
&$+\int\int_{T_0} \frac{dx_{i_1}}{x_{i_1}}\circ
\left(\frac{dx_{i_2}}{x_{i_2}}\wedge\frac{dx_{i_3}}{x_{i_3}}\right)

+\int_{\gamma}\frac{dx_{i_1}}{x_{i_1}}
\int\int_{T_0}\frac{dx_{i_2}}{x_{i_2}}\wedge\frac{dx_{i_3}}{x_{i_3}}+$\\
\\
&$+\int_{\gamma}\frac{dx_{i_1}}{x_{i_1}}
\int\int_{S^{-1}}
\frac{dx_{i_2}}{x_{i_2}}\wedge\frac{dx_{i_3}}{x_{i_3}}$
\end{tabular}$$
Using the above relation for the domain
$S\cup\gamma S^{-1}\gamma^{-1}$,
we obtain
part (b) of the Lemma. Part (a) can be done in a similar way,
when we write $S'$ instead of $S$ and $\tau'$ instead of $\tau$.
%%%%%%%%%%%%%%%%%%%%%%%%%%%%%%%%%%%%%%%%%%%%%%%%%%
Denote by $S'^{-1}-S^{-1}$ the region bounded by
$\tau_j\tau^{-1}$. Then we have the following lemma.
%%%%%%%%%%%%%%%%%%%%%%%%%%%%%%%%%%%%%%%%%%%%%%%%%%%%%%%%%%
\begin{lemma}
$$\begin{tabular}{llll}
$\int_{[\sigma,\tau']}
\frac{dx_{i_1}}{x_{i_1}}\circ\frac{dx_{i_2}}{x_{i_2}}\circ\frac{dx_{i_3}}{x_{i_3}}
-\int_{[\sigma,\tau]}
\frac{dx_{i_1}}{x_{i_1}}\circ\frac{dx_{i_2}}{x_{i_2}}\circ\frac{dx_{i_3}}{x_{i_3}}=$\\
\\
$=\int_{\sigma_0}\frac{dx_{i_1}}{x_{i_1}}
\int\int_{S'^{-1}-S^{-1}}\frac{dx_{i_2}}{x_{i_2}}\wedge\frac{dx_{i_3}}{x_{i_3}}+$\\
\\
$+\int\int_{S'^{-1}-S^{-1}}\frac{dx_{i_1}}{x_{i_1}}\wedge\frac{dx_{i_2}}{x_{i_2}}
\int_{\sigma_0^{-1}}\frac{dx_{i_3}}{x_{i_3}}.$
\end{tabular}$$
\end{lemma}
%%%%%%%%%%%%%%%%%%%%%%%%%%%%%%%%%%%%%%%%%%%%%%%%%%%%%%%%%%
For the region $S'^{-1}-S^{-1}$, we have the following lemma.
%%%%%%%%%%%%%%%%%%%%%%%%%%%%%%%%%%%%%%%%%%%%%%%%%%%%%%%%%%
\begin{lemma}
$$\int\int_{S'^{-1}-S^{-1}}\frac{dx_{i_1}}{x_{i_1}}\wedge\frac{dx_{i_2}}{x_{i_2}}
=\int_{\tau'\tau^{-1}}\frac{dx_{i_1}}{x_{i_1}}\circ\frac{dx_{i_2}}{x_{i_2}}.$$
\end{lemma}
%%%%%%%%%%%%%%%%%%%%%%%%%%%%%%%%%%%%%%%%%%%%%%%%%%%%%%%%%%
\proof  It follows from non-commutative Stokes theorem for
iterated integrals (Theorem 1.5).
%%%%%%%%%%%%%%%%%%%%%%%%%%%%%%%%%%%%%%%%%%%%%%%%%%%%%%%%%%%%%%

\proof (of Proposition 3.8) In order to prove the proposition it
remains to show that
$$\lim_{\e\rightarrow 0}
\int_{\tau'\tau^{-1}}
\frac{dx_{i_1}}{x_{i_1}}\circ\frac{dx_{i_2}}{x_{i_2}}
=0.$$
Afterwards, we can use the above three lemmas.

If $i_1\neq 0$ or $i_2\neq 0$ then the above integral becomes $0$
in the limit $\e\rightarrow 0$, If $i_1=i_2=0$ then
$\frac{dx_{i_1}}{x_{i_1}}\wedge\frac{dx_{i_2}}{x_{i_2}}=0$. And by
Lemma 3.6, we have that
$$\int_{\tau'\tau^{-1}}
\frac{dx_{i_1}}{x_{i_1}}\circ\frac{dx_{i_2}}{x_{i_2}}
=0.$$
%%%%%%%%%%%%%%%%%%%%%%%%%%%%%%%%%%%%%%%%%%%%%%%%%%%%%%%%%%%%%%%%%%%%%%%%%%

Denote by $$\omega_k=\frac{df_k}{f_k}$$ for $k=1,2,3$. Using
Stokes formula for $1$- and $2$-dimensional iterated integrals
(Theorem 1.5), we obtain the following lemma.
%%%%%%%%%%%%%%%%%%%%%%%%%%%%%%%%%%%%%%%%%%%%%%%%%%%%%%%%%%%%%%%%%%
\begin{lemma}
$$\int_{\d T'}
\omega_1\circ\omega_2\circ\omega_3
=\int\int_{T'}
\omega_1\circ(\omega_2\wedge\omega_3)
+\int\int_{T_{jl}}
(\omega_1\wedge\omega_2)\circ\omega_3.$$
\end{lemma}
%%%%%%%%%%%%%%%%%%%%%%%%%%%%%%%%%%%%%%%%%%%%%%%%%%%%%%%%%%%%%%%%%%%%%%%%%%%%%%%
Using Lemma 1.6 for composition of domains of $2$-dimensional
iterated integrals, we obtain the following lemma.
\begin{lemma}
$$\begin{tabular}{ll}
(a)$\int\int_{T'} \omega_1\circ(\omega_2\wedge\omega_3)
=\int\int_{T_0} \omega_1\circ(\omega_2\wedge\omega_3)
+\int_{\gamma_0}\omega_1
\int\int_{T_0}\omega_2\wedge\omega_3
+\int_{\sigma_0}\omega_1\int\int_{S'^{-1}}
\omega_2\wedge\omega_3$\\
\\
(b)$\int\int_{T'} (\omega_1\wedge\omega_2)\circ\omega_3
=\int\int_{T_0} (\omega_1\wedge\omega_2)\circ\omega_3
+\int\int_{T_0}\omega_1\wedge\omega_2\int_{\gamma^{-1}}\omega_3
+\int\int_{S'^{-1}} \omega_1\wedge\omega_2\int_{\sigma_0^{-1}}\omega_3$
\end{tabular}$$
\end{lemma}
\proof
It is the same as the proof of Lemma 3.7.
%%%%%%%%%%%%%%%%%%%%%%%%%%%%%%%%%%%%%%%%%%%%%%%%%%%%%%%%
\begin{lemma}
$$\begin{tabular}{ll}
$(a) \frac{1}{(2\pi
i)^2}\int\int_{T'}\omega_1\circ(\omega_2\wedge\omega_3)=$&
$-(m_1+n_1)(m_2n_3-m_3n_2)\pi i+$\\
\\
&$-(m_2n_3-m_3n_2)\int_{\gamma_0}\frac{dg_1}{g_1}+$\\
\\
&$-m_3n_1\int_{\gamma_0}\frac{dg_2}{g_2}-$\\
\\
&$+m_2n_1\int_{\gamma_0}\frac{dg_3}{g_3}+O(\e);$\\
\end{tabular}$$
$$\begin{tabular}{ll}
$(b) \frac{1}{(2\pi i)^2}\int\int_{T'}
(\omega_1\wedge\omega_2)\circ\omega_3=$&
$(m_3+n_3)(m_1n_2-m_2n_1)\pi i+$\\
\\
&$-(m_1n_2-m_2n_1)\int_{\gamma_0}\frac{dg_3}{g_3}-$\\
\\
&$+m_1n_3\int_{\gamma_0}\frac{dg_2}{g_2}+$\\
\\
&$-m_2n_3\int_{\gamma_0}\frac{dg_1}{g_1}+O(\e).$\\
\end{tabular}$$
\end{lemma}
%%%%%%%%%%%%%%%%%%%%%%%%%%%%%%%%%%%%%%%%%%%%%%%%%%%%%%%%%%%%%%%%%%%%%%
\proof For part (a) we have
$$\begin{tabular}{ll}
$\int\int_{T_0}\omega_1\circ(\omega_2\wedge\omega_3)=$
&$+\int\int_{T_0} \omega_1\circ(\omega_2\wedge\omega_3)+$\\
\\
&$+\int_{\gamma}\omega_1\int\int_{T_0}\omega_2\wedge\omega_3
+\int_{\sigma}\omega_1
\int\int_{S'^{-1}}
\omega_2\wedge\omega_3$
\end{tabular}$$
From Lemma 4.3 part(b), we have
$$\int\int_{T_0}
\omega_1\circ(\omega_2\wedge\omega_3)
=-\frac{1}{2}(2\pi i)^3(m_1+n_1)(m_2n_3-m_3n_2).$$
This takes care of the first
summand in this lemma.
Consider the functions
$f_k$ and $g_k$
in terms of the local coordinates $x$ and $y$
near the point $P_{jl}$.
Using Lemma 4.3 part (a), we obtain
$$\int_{\gamma}\omega_1
\int\int_{T_0}\omega_2\wedge\omega_3=
-(m_2n_3-m_3n_2)(2\pi i)^2\int_{\gamma}\frac{df_1}{f_1}.$$ When
we express $f_1$ in term of $g_1$ and powers of $x$ and $y$,
we obtain
$$
\begin{tabular}{ll}
$\frac{1}{(2\pi
i)^2}\int_{\gamma}\omega_1\int\int_{T_0}\omega_2\wedge\omega_3$&
$=-(m_2n_3-m_3n_2)
\int_{\gamma}\frac{d(x^{m_1}y^{n_1}g_1)}
{x^{m_1}y^{n_1}g_1}=$\\
\\
&$=-(m_2n_3-m_3n_2)
\left(n_1\int_{\gamma}\frac{dy}{y}
+\int_{\gamma}\frac{dg_1}{g_1}\right).$
\end{tabular}
$$
Now we have to compute
$$\int_{\sigma}\omega_1
\int\int_{S'^{-1}}
\omega_2\wedge\omega_3.$$
We have
$$\int_{\sigma}\omega_1=m_1(2\pi i).$$
We can simplify the double integral by using that
$S'^{-1}=\gamma^{-1}\times\tau'$. Also, we can express $f_k$ as
$x^{m_k}y^{n_k}g_k$ and then use that under the integral we
must have $dx/x$, since we integrate the variable $x$ over
the loop $\tau_0$. For the differential $2$-form under the
integral, we have
$$\omega_2\wedge\omega_3=
\left(m_2\frac{dx}{x}+n_2\frac{dy}{y}+\frac{dg_2}{g_2}\right)
\wedge\left(m_3\frac{dx}{x}+n_3\frac{dy}{y}+\frac{dg_3}{g_3}\right).$$
When we pick only the terms that involve $dx/x$, we obtain the
$2$-form
$$(m_2n_3-m_3n_2)\frac{dx}{x}\wedge\frac{dy}{y}
+m_2\frac{dx}{x}\wedge\frac{dg_3}{g_3}
-m_3\frac{dx}{x}\wedge\frac{dg_2}{g_2}.$$
The term $-(m_2n_3-m_3n_2)\frac{dx}{x}\wedge\frac{dy}{y}$
gives
$$\int\int_{S'^{-1}}-(m_2n_3-m_3n_2)\frac{dx}{x}
\wedge\frac{dy}{y}
=-(m_2n_3-m_3n_2)(2\pi i)
\int_{\gamma_0^{-1}}\frac{dy}{y}.$$ This
term cancels with the term
$$\lim_{\e\rightarrow 0}-(m_2n_3-m_3n_2)n_1\int_{\gamma}\frac{dy}{y},$$
coming from
$$\int_{\gamma}\omega_1\int\int_{T_0}\omega_2\wedge\omega_3.$$

The $2$-form $$m_2\frac{dx}{x}\wedge\frac{dg_3}{g_3}$$ gives
$$\int\int_{S'^{-1}}m_2\frac{dx}{x}\wedge\frac{dg_3}{g_3}
=-m_2(2\pi i)
\int_{\gamma_0^{-1}}\frac{dg_3}{g_3}.$$

Finally, the
$2$-form $$m_3\frac{dx}{x}\wedge\frac{dg_2}{g_2}$$ gives
$$\int\int_{S'^{-1}}m_3\frac{dx}{x}\wedge\frac{dg_2}{g_2}
=m_3(2\pi i)
\int_{\gamma_0^{-1}}\frac{dg_2}{g_2}.$$
This proves part
(a) of the lemma.

Part (b) can be proven in the same way.

%%%%%%%%%%%%%%%%%%%%%%%%%%%%%%%%%%%%%%%%%%%%%%%%%%%%%%%%%%%%%%%%%%%%
\subsection{Differential equation}
We will define a partial differential equation, which has no
solution locally, but has a solution only over a fixed path, when
it is reduced to an ordinary differential equation. Let $A_i$ for
$i=1,\dots,n$ be indeterminant constants which do not commute. One
can think of them as constant square matrices. Let $h_i$ for
$i=1,\dots,n$ be rational functions on $X$. Consider the
differential equation
$$dF=F\sum_{i=1}^n A_i \frac{dh_i}{h_i},$$
where the multiplication is matrix multiplication (or
multiplication of indeterminant constants). When this partial
differential equation is restricted to a path, we obtain a
solution, given by a generating series of iterated integrals.
%%%%%%%%%%%%%%%%%%%%%%%%%%%%%%%%%%%%%%%%%%%%%%%%%%%%%%%%%%%%%%%

From Lemma 2.6 it seems that
the differential $\frac{df_k}{f_k}$
should anti-commute.
Thus, our first attempt to relate a differential equation
to a logarithmic version of the Parshin symbol is the following
$$dF=F\sum_{k=1}^3 A_k \frac{df_k}{f_k},$$
where the indeterminant $A_k$ for $k=1,2,3$ anti-commute.

Why this attempt is not good? The factor $x_1$ from $f_1$ and
$x_1$ from $f_2$ lead to summands $\frac{dx_1}{x_1}$ in
$\frac{df_1}{f_1}$ and $\frac{dx_1}{x_1}$ in $\frac{df_2}{f_2}$.
These two copies of $\frac{dx_1}{x_1}$ commute, when iterated.
However, because the constants $A_1$ and $A_2$ anti-commute in the
above differential equation, we obtain that the two copies of
$\frac{dx_1}{x_1}$  - one coming from $f_1$ and the other from
$f_2$ - should anti-commute. Thus, the both commute and
anti-commute, when we use the above differential equation.
%%%%%%%%%%%%%%%%%%%%%%%%%%%%%%%%%%%%%%%%%%%%%%%%%%%%%%%%

Instead of considering the functions $f_k$ for $k=1,2,3$, we must
decompose the three functions into factors. We want these factors
to enter in the choice of local coordinates. (Later we will show
independence of the choices of local coordinates.) Also, we want
similar factors coming from different $f_k$'s to commute and
different factors coming from different functions $f_k$ to
anti-commute. In order to do that we must distinguish between
similar factors, coming from different functions $f_k$'s.
%%%%%%%%%%%%%%%%%%%%%%%%%%%%%%%%%%%%%%%%%%%%%%%%%
\begin{definition} (factorization of the functions
$f_k$ for $k=1,2,3$.) We are going to use the definitions of $x_i$
and the corresponding powers $n_{ki}$ for $i=1,\dots,N$ from
Definition 2.1. Define $x_{N+k}$ for $k=1,2,3$ by
$$x_{N+k}=f_k\prod_{i=1}^N x_i^{-m_{ki}}.$$
Define also the constants $m_{N+l,k}$ for $k,l=1,2,3$ by
$$m_{N+l,k}=\delta_{kl},$$
where $\delta_{kl}$ is the Kroneker delta function. We are going
to use the factorizations
$$f_k=\prod_{i=1}^{N+3}x_i^{m_{ki}}.$$
Then we have
$$\frac{df_k}{f_k}=\sum_{i=1}^{N+3}m_{ki}\frac{dx_i}{x_i}.$$
\end{definition}
%%%%%%%%%%%%%%%%%%%%%%%%%%%%%%%%%%%%%%
\begin{definition}
Consider the differential equation
$$dF=F\sum_{k=1}^3\sum_{i=1}^{N+3}
A_{ki} m_{ki}\frac{dx_i}{x_i}.$$
We are going to use solutions of this differential equation along a path $\gamma$, which we
will denote by
$$F_{\gamma}.$$
\end{definition}
%%%%%%%%%%%%%%%%%%%%%%%%%%%%%%%%%%%%%%%%%
\begin{definition} (Variable $A$ for the purpose of the new symbol)
We set an equivalence among certain monomials in the variables
$A_{ki}$. Let
$$A\sim A_{1,j_1}A_{2,j_2}A_{3,j_3}$$
for all values of $j_1$, $j_2$ and $j_3$.
\end{definition}
%%%%%%%%%%%%%%%%%%%%%%%%%%%%%%%%%%%%%%%%%%%%%%%%%
\begin{definition} (Variable $B$ for the purpose of the logarithm of the Parshin symbol)
Let $\sigma$ be a permutation of $\{1,2,3\}$.
Define an equivalence $$B\sim A_{1,j_1}A_{2,j_2}A_{3,j_3}
\sim A_{\sigma(1),j_1}A_{\sigma(3),j_2}A_{\sigma(3),j_3}.$$
For $j_1\neq j_2\neq j_3 \neq j_1$, let
$$A_{1,j_{\sigma(1)}}A_{2,j_{\sigma(2)}}A_{3,j_{\sigma(3)}}
\sim
sign(\sigma)A_{1,j_1}A_{2,j_2}A_{3,j_3}.$$
For $j_1=j_2\neq j_3$ and for a permutation $\sigma$ of $\{1,2\}$, define
the equivalence
$$A_{1,j_{\sigma(1)}}A_{2,j_{\sigma(1)}}A_{3,j_{\sigma(2)}}
\sim
sign(\sigma)A_{1,j_1}A_{2,j_1}A_{3,j_2}.$$
For $j_1=0$ and $j_2\neq 0$ define
$$B\sim A_{1,j_1}A_{2,j_1}A_{3,j_2}.$$
\end{definition}
%%%%%%%%%%%%%%%%%%%%%%%%%%%%%%%%%%%%%%%%%%%%%%%%%%%%%%%%%%%%%%%%
%%%%%%%%%%%%%%%%%%%%%%%%%%%%%%%%%%%%%%%%%%%%%%%%%%%%%%%%%%%%%%%%%%%
\subsection{The commutators $[\a_i,\b_i]$}
In order to prove a reciprocity law for the logarithm of the new
symbol $Log[f_1,f_2,f_3]^{x_0,\gamma_i}_{C_0,P_i}$, we have to
consider the loops on $C_0$, given by the commutator
$[\a_j,\b_j]$.
%%%%%%%%%%%%%%%%%%%%%%%%%%%%%%%%%%%%%%%%%
\begin{definition}
Define a torus
$$T'_j=h(u,[\a_j,\b_j]),$$
for $|u|=\e_1$.
\end{definition}
%%%%%%%%%%%%%%%%%%%%%%%%%%%%%%%%%%%%%%%%%%%%%%%%%%%%%%%%%%%%%%%%%%%%
\begin{lemma}
Let $\a$ be any loop on $X$ and let $f$ be a rational function on
$X$. Then
$$\int_\a \frac{df}{f}=2\pi in$$
for some integer $n$.
\end{lemma}
\proof Consider $f$ as a function from the variety $X$ to $\C
P^1-\{0,\infty\}$. On $\C P^1-\{0,\infty\}$ we have the
differential form $dz/z.$ Then
$$\int_\a \frac{df}{f}=\int_\a f^*\frac{dz}{z}=\int_{f_*\a}\frac{dz}{z}=2\pi in$$
for some integer $n$.

We are going to use the following lemma for commutators.
\begin{lemma}
(a)$\int_{[\a,\b]}\omega=0;$

(b) $\int_{[\a,\b]}\omega_1\circ\omega_2
=\int_{\a}\omega_1\int_{\b}\omega_2-\int_{\b}\omega_1\int_{\a}\omega_2;$

$$\begin{tabular}{lll} (c)
$\int_{[\a,\b]}\omega_1\circ\omega_2\circ\omega_3$
&=&$\int_{\a}\omega_1\circ\omega_2\int_{\b}\omega_3
-\int_{\b}\omega_1\circ\omega_2\int_{\a}\omega_3+$\\
\\
&&$+\int_{\a}\omega_3\circ\omega_2\int_{\b}\omega_1
-\int_{\b}\omega_3\circ\omega_2\int_{\a}\omega_1-$\\
\\
&&$-\int_{\a}\omega_1\int_{\b}\omega_2\int_{\a}\omega_3
+\int_{\b}\omega_1\int_{\a}\omega_2\int_{\b}\omega_3.$
\end{tabular}$$
\end{lemma}
\proof Part (a) is trivial. Part (b) is Theorem 3.1 in \cite{H}.
Part (c) is Theorem 4. in \cite{H}.
%%%%%%%%%%%%%%%%%%%%%%%%%%%%%%%%%%%%%%%%%%%%%%%%%%%%%%%%%%%
\begin{proposition}The 2-dimensional iterated integrals
$$\int\int_{\d T'_j}
\frac{df_1}{f_1}\circ\frac{df_2}{f_2}\circ\frac{df_3}{f_3}
$$
is an integer times $(2\pi i)^3$.
\end{proposition}
\proof Let $\sigma$ be the loop in the fiber of the torus $T'_j$.
First we prove the following lemma.
%%%%%%%%%%%%%%%%%%%%%%%%%%%%%%%%%%%%%%%%%%%%%%%%%%%%%%%%%%%%%%%%%%%%%%%%%%
\begin{lemma}
$$(2\pi i)^{-3}\int_{[[\alpha,\beta],\sigma]}
\frac{df_1}{f_1}\circ\frac{df_2}{f_2}\circ\frac{df_3}{f_3}$$ is an integer.
\end{lemma}
\proof Using part (c) of the Lemma 5.7 we obtain
$$\begin{tabular}{lll}
$\int_{[[\alpha,\beta],\tau]}
\frac{df_1}{f_1}\circ\frac{df_2}{f_2}\circ\frac{df_3}{f_3}$
&=&$\int_{[\alpha,\beta]}\frac{df_1}{f_1}\circ\frac{df_2}{f_2}\int_{\tau}\frac{df_3}{f_3}
-\int_{\tau}\frac{df_1}{f_1}\circ\frac{df_2}{f_2}\int_{[\alpha,\beta]}\frac{df_3}{f_3}+$\\
\\
&&$+\int_{[\alpha,\beta]}\frac{df_3}{f_3}\circ\frac{df_2}{f_2}\int_{\tau}\frac{df_1}{f_1}
-\int_{\tau}\frac{df_3}{f_3}\circ\frac{df_2}{f_2}\int_{[\alpha,\beta]}\frac{df_1}{f_1}-$\\
\\
&&$-\int_{[\alpha,\beta]}\frac{df_1}{f_1}\int_{\tau}\frac{df_2}{f_2}\int_{[\alpha,\beta]}\frac{df_3}{f_3}
+\int_{\tau}\frac{df_1}{f_1}\int_{[\alpha,\beta]}\frac{df_2}{f_2}\int_{\tau}\frac{df_3}{f_3}.$
\end{tabular}$$

From part (a) of Lemma 5.7, we have that an integral of a $1$-form
over a commutator is zero. Therefore,
$$\int_{[[\alpha,\beta],\tau]}
\frac{df_1}{f_1}\circ\frac{df_2}{f_2}\circ\frac{df_3}{f_3}=
\int_{[\alpha,\beta]}\frac{df_1}{f_1}\circ\frac{df_2}{f_2}\int_{\tau}\frac{df_3}{f_3}
+\int_{[\alpha,\beta]}\frac{df_3}{f_3}\circ\frac{df_2}{f_2}\int_{\tau}\frac{df_1}{f_1}.$$

Now we use part (b) of Lemma 5.7. And we obtain
$$\begin{tabular}{lll}
$\int_{[[\alpha,\beta],\tau]}
\frac{df_1}{f_1}\circ\frac{df_2}{f_2}\circ\frac{df_3}{f_3}$
&=&$\left(\int_{\a}\frac{df_1}{f_1}\int_{\b}\frac{df_2}{f_2}
-\int_{\b}\frac{df_1}{f_1}\int_{\a}\frac{df_2}{f_2}\right)
\int_{\tau}\frac{df_3}{f_3}+$\\
\\
&&$+\left(\int_{\a}\frac{df_3}{f_3}\int_{\b}\frac{df_2}{f_2}
-\int_{\b}\frac{df_3}{f_3}\int_{\a}\frac{df_2}{f_2}\right)
\int_{\tau}\frac{df_1}{f_1}$
\end{tabular}$$
By Lemma 5.6, all of the integrals of $1$-forms are integer
multiples of $2\pi i$. This proves the lemma.
%%%%%%%%%%%%%%%%%%%%%%%%%%%%%%%%%%%%%%%%%%%

\proof (of Proposition 5.15) The boundary of the torus $T'$ is
$[\sigma[\a,\b]]$. By the above lemma the integral is an integer
multiple of $(2\pi i)^3$.
%%%%%%%%%%%%%%%%%%%%%%%%%%%%%%%%%%%%%%%%%%%%%%%%%%%%%%%%%%%%%%%%%%%%
\begin{definition} ($h_k$)
We have that $m_k$ is the order of the function $f_k$ along the
curve $C_0$. Also we have that $x$ is a rational function on $X$
that has zero of order $1$ along $C_0$. Define functions $h_k$ by
$$f_k=x^{m_k}h_k.$$
\end{definition}
%%%%%%%%%%%%%%%%%%%%%%%%%%%%%%%%%%%%%%%%%%%%
\begin{corollary}
For fixed $j$ the contribution from the commutator $[\a_j,\b_j]$
is the iterated integrals
$$
\begin{tabular}{lll}
$\int\int_{\d T'_j}
\frac{df_1}{f_1}\circ\frac{df_2}{f_2}\circ\frac{df_3}{f_3}=$
&$2\pi i m_1\left(\int_{\a_j}\frac{dh_3}{h_3}\int_{\b_j}\frac{dh_2}{h_2}
-\int_{\b_j}\frac{dh_3}{h_3}\int_{\a_j}\frac{dh_2}{h_2}\right)+$\\
\\
&$+2\pi i m_3\left(\int_{\a_j}\frac{dh_1}{h_1}\int_{\b_j}\frac{dh_2}{h_2}
-\int_{\b_j}\frac{dh_1}{h_1}\int_{\a_j}\frac{dh_2}{h_2}\right).$
\end{tabular}
$$
\end{corollary}
%%%%%%%%%%%%%%%%%%%%%%%%%%%%%%%%%%%%%%%%%%%%%%%%%%%%%%%%%%%%%%%%%%%%%
\subsection{Extra terms}
The source of the extra terms are the integrals
$$\int_{\pi_{i}}\frac{df_1}{f_1}\int_{[\mu_{i+1},\tau]}
\frac{df_2}{f_2}\circ\frac{df_3}{f_3}$$
and
$$\int_{[\mu_{i+1},\tau]}
\frac{df_1}{f_1}\circ\frac{df_2}{f_2}
\int_{\pi_{i}^{-1}}\frac{df_3}{f_3}$$ from Definition 2.29. Note
that their difference is the same as the sum of the coefficients
equivalent to $A$ in
$F_{\pi_{i}[\mu_{i+1},\tau]\pi_{i}^{-1}}-F_{[\mu_{i+1},\tau]}$.

In terms of the differential equation, an extra term is
the difference between the coefficients of $A_{1i_1}A_{2i_2}A_{3i_3}$ in
$F_{\pi_i[\mu_{i+1},\tau]\pi_i^{-1}}$ and in $F_{[\mu_{i+1},\tau]}$.
%%%%%%%%%%%%%%%%%%%%%%%%%%%%%%%%%%%%%%%%%%%%%%%%%%
\begin{lemma}
The difference between the coefficients of $A_{1j_1}A_{2j_2}A_{3j_3}$ in
$F_{\pi_i[\mu_{i+1},\tau]\pi_i^{-1}}$ and in $F_{[\mu_{i+1},\tau]}$ is
$$(2\pi i)^2 n_{j_1}n_{2j_2}n_{3j_3}
\left(\int_{\pi_i}\frac{dx_{j_1}}{x_{j_1}}
\int_{[\mu_{i+1},\tau]}
\frac{dx_{j_2}}{x_{j_2}}\circ\frac{dx_{j_3}}{x_{j_3}}
-\int_{\pi_i}\frac{dx_{j_3}}{x_{j_3}}
\int_{[\mu_{i+1},\tau]}
\frac{dx_{j_1}}{x_{j_1}}\circ\frac{dx_{j_2}}{x_{j_2}}
\right).$$
\end{lemma}
%%%%%%%%%%%%%%%%%%%%%%%%%%%%%%%%%%%%%%
\proof Cut the loop $\pi_i[\mu_{i+1},\tau]\pi_i^{-1}$ into $3$ loops
$\pi_i$, $[\mu_{i+1},\tau]$ and $\pi_i^{-1}$.
Then use Corollary 1.3 for composition of paths.
%%%%%%%%%%%%%%%%%%%%%%%%%%%%%%%%%%%%%%%%%%%%%%%%%%%%%%%%%%%%

Recall that $L_j$ is the number of intersection points of $C_0$
and $C_j$.
%%%%%%%%%%%%%%%%%%%%%%%%%%%%%%%%%%%%%%%%%%%%%%%
\begin{definition}
For $k=1,2,3$ consider the orders of vanishing $m_k$ and $m_{kj}$
of $x$ and $x_j$, respectively. (see Definitions 3.23 and 3.25.)
Let
$$D_{1}(j)=\left|\begin{tabular}{ll}
$m_2$ & $n_{2j}$\\
$m_3$ & $n_{3j}$
\end{tabular}\right|,\mbox{ }
D_{2}(i)=\left|\begin{tabular}{ll}
$m_3$ & $n_{3j}$\\
$m_1$ & $n_{1j}$
\end{tabular}\right|,\mbox{ }
D_{3}(i)=\left|\begin{tabular}{ll}
$m_1$ & $n_{1j}$\\
$m_2$ & $n_{2j}$
\end{tabular}\right|.$$
\end{definition}
%%%%%%%%%%%%%%%%%%%%%%%%%%%%%%%%%%%%%%%%%%%%%%%%%%%%
\begin{theorem}
The sum of the coefficients equivalent to $A$ in
$$\sum_i(F_{\pi_i[\mu_{i+1},\tau]\pi_i^{-1}}-F_{[\mu_{i+1},\tau]})$$
is
$$\sum_{j_1<j_2}
(n_{1j_1}D_1(j_2)-n_{3j_1}D_3(j_2))L_{j_1}L_{j_2}.
+\frac{1}{2}\sum_{j_2=1}^N
(n_{1j_2}D_1(j_2)-n_{3j_2}D_3(j_2))L_{j_2}(L_{j_2}-1).
$$
\end{theorem}
%%%%%%%%%%%%%%%%%%%%%%%%%%%%%%%%%%%%%%%
\proof Fix $j_1$, $j_2$ and $j_3$. We are going to use the above
lemma. We consider the sum
$$\sum_i
\int_{\pi_i}\frac{dx_{j_1}}{x_{j_1}}
\int_{[\mu_{i+1},\tau]}
\frac{dx_{j_2}}{x_{j_2}}\circ\frac{dx_{j_3}}{x_{j_3}}.$$
The second integral is zero if $j_2\neq 0$ and $j_3\neq 0$. Let $j_3=0$.
Again the second integral is zero
if $\mu_{i+1}\neq \sigma_{j_2,l}$
for all $l=1,\dots,L_{j_2}$.
Suppose that $\mu_{i+1}= \sigma_{j_2,l}.$  Then
$$\int_{[\mu_{i+1},\tau]}
\frac{dx_{j_2}}{x_{j_2}}\circ\frac{dx_{0}}{x_{0}}
=(2\pi i)^2.$$
Let $l>1$ so that $\mu_{i+1}= \sigma_{j_2,l}.$ Then
$$\pi_i=\pi_{j_2,l-1}=
\left(
\prod_{j'=1}^{j_2-1}
\left(
\prod_{l'=1}^{L_{j'}}\tilde{\sigma}_{j'l'}
\right)\right)
\prod_{l'=1}^{l-1}\tilde{\sigma}_{j_2l'}.
$$
Then we have
$$(2\pi i)^{-1}\int_{\pi_i}\frac{dx_{j_1}}{x_{j_1}}
=\left\{
\begin{tabular}{lll}
$L_{j_1}$ & for $j_1<j_2$\\
\\
$l-1$& for $j_1=j_2$\\
\\
$0$ &  for $j_1>j_2$.
\end{tabular}
\right.
$$
Therefore,
$$(2\pi i)^{-3}\sum_{l=2}^{L_{j_2}}\int_{\pi_{j_2,l-1}}\frac{dx_{j_1}}{x_{j_1}}
\int_{[\mu_{j_2,l},\tau]}
\frac{dx_{j_2}}{x_{j_2}}\circ\frac{dx_{0}}{x_{0}}.
=\left\{
\begin{tabular}{lll}
$L_{j_1}(L_{j_2}-1)$ & for $j_1<j_2$\\
\\
$\frac{1}{2}L_{j_2}(L_{j_2}-1)$& for $j_1=j_2$\\
\\
$0$ &  for $j_1>j_2$.
\end{tabular}
\right.
$$
If $l=1$ in $\mu_{i+1}= \sigma_{j_2,l}=\sigma_{j_2,1}.$ Then
$$\pi_i=\pi_{j_2-1,L_{j_2-1}}=
\prod_{j'=1}^{j_2-1}
\left(
\prod_{l'=1}^{L_{j'}}\tilde{\sigma}_{j'l'}
\right).
$$
The contribution of the extra term for $l=1$ is
$$(2\pi i)^{-3}
\int_{\pi_{j_2,l-1}}\frac{dx_{j_1}}{x_{j_1}}
\int_{[\mu_{j_2,l},\tau]}
\frac{dx_{j_2}}{x_{j_2}}\circ\frac{dx_{0}}{x_{0}}
=\left\{
\begin{tabular}{lll}
$L_{j_1}$ & for $j_1<j_2$\\
\\
$0$ &  for $j_1\geq j_2$.
\end{tabular}
\right.
$$
Therefore,
$$(2\pi i)^{-3}\sum_i\int_{\pi_i}\frac{dx_{j_1}}{x_{j_1}}
\int_{[\mu_{i},\tau]}
\frac{dx_{j_2}}{x_{j_2}}\circ\frac{dx_{0}}{x_{0}}.
=\left\{
\begin{tabular}{lll}
$L_{j_1}L_{j_2}$ & for $j_1<j_2$\\
\\
$\frac{1}{2}L_{j_2}(L_{j_2}-1)$& for $j_1=j_2$\\
\\
$0$ &  for $j_1>j_2$.
\end{tabular}
\right.
$$

Therefore, the coefficient of $A_{1j_1}A_{2j_2}A_{30}$ in
$$\sum_i(F_{\pi_i[\mu_{i+1},\tau]}$$
is
$$
\begin{tabular}{lll}
$n_{1j_1}n_{2j_2}m_3 L_{j_1}L_{j_2}$ & for $j_1<j_2$\\
\\
$\frac{1}{2}n_{1j_2}n_{2j_2}m_3L_{j_2}(L_{j_2}-1)$& for $j_1=j_2$\\
\\
$0$ &  for $j_1>j_2$.
\end{tabular}
$$
Similarly, the coefficient of $A_{1j_1}A_{20}A_{3j_2}$ in
$$\sum_i(F_{\pi_i[\mu_{i+1},\tau]}$$
is
$$
\begin{tabular}{lll}
$-n_{1j_1}n_{3j_2}m_2 L_{j_1}L_{j_2}$ & for $j_1<j_2$\\
\\
$-\frac{1}{2}n_{1j_2}n_{3j_2}m_2L_{j_2}(L_{j_2}-1)$& for $j_1=j_2$\\
\\
$0$ &  for $j_1>j_2$.
\end{tabular}
$$
Therefore, the coefficient of $A_{1j_1}A_{2j_2}A_{30}$  plus the
coefficient of $A_{1j_1}A_{20}A_{3j_2}$ in
$$\sum_i(F_{\pi_i[\mu_{i+1},\tau]}$$
is
$$
\begin{tabular}{lll}
$n_{1j_1}D_1(j_2) L_{j_1}L_{j_2}$ & for $j_1<j_2$\\
\\
$\frac{1}{2}n_{1j_2}D_1(j_2)L_{j_2}(L_{j_2}-1)$& for $j_1=j_2$\\
\\
$0$ &  for $j_1>j_2$.
\end{tabular}
$$
Note that the coefficient of $A_{10}A_{2j_2}A_{3j_1}$ in
$$\sum_i(F_{\pi_i[\mu_{i+1},\tau]}$$ is zero.

Similarly, the coefficient of $A_{10}A_{2j_2}A_{3j_1}$  plus the
coefficient of $A_{1j_2}A_{20}A_{3j_1}$ in
$$\sum_i(F_{[\mu_{i+1},\tau]\pi^{-1}}$$
is
$$
\begin{tabular}{lll}
$-n_{3j_1}D_3(j_2) L_{j_1}L_{j_2}$ & for $j_1<j_2$\\
\\
$-\frac{1}{2}n_{3j_2}D_3(j_2)L_{j_2}(L_{j_2}-1)$& for $j_1=j_2$\\
\\
$0$ &  for $j_3>j_2$.
\end{tabular}$$

The contribution of all the extra terms is
$$\sum_{j_1<j_2}
(n_{1j_1}D_1(j_2)-n_{3j_1}D_3(j_2))L_{j_1}L_{j_2}.
+\frac{1}{2}\sum_{j_2=1}^N
(n_{1j_2}D_1(j_2)-n_{3j_2}D_3(j_2))L_{j_2}(L_{j_2}-1)
$$

%%%%%%%%%%%%%%%%%%%%%%%%%%%%%%%%%%%%%%%%%%%%%
\begin{corollary}
Each of the extra terms is and integer multiple of
$(2\pi i)^3$.
\end{corollary}
%%%%%%%%%%%%%%%%%%%%%%%%%%%%%%%%%%%%%%%%%%%%%%
\subsection{New symbol}
\begin{definition} (Logarithmic version of the new symbol)
For $k=1,2,3$ define $m_k=n_{k1}$, $n_k=n_{kj}$
and $g_{kj}=x_1^{-m_k}x_j^{-n_k}f_k=\prod_{i\neq 1,j}x_i^{n_{ki}}$

Define a logarithm of the new symbol as
$$
\begin{tabular}{ll}
$Log[f_1,f_2,f_3]^{\gamma_{jl}}_{C_0,P_{jl}}
=\frac{1}{2}(2\pi i)^3
((m_1+n_1)D_1-(m_3+n_3)D_3)+$\\
\\
$+(D_1+m_2n_3)\int_{\gamma_{jl}}\frac{dg_{1j}}{g_{1j}}
+D_2\int_{\gamma_{jl}}\frac{dg_{2j}}{g_{2j}}
+(D_3-m_2n_1)\int_{\gamma_{jl}}\frac{dg_{3j}}{g_{3j}},$
\end{tabular}
$$
where
$$D_1=\left|
\begin{tabular}{ll}
$m_2$ & $n_{2}$\\
$m_3$ & $n_{3}$
\end{tabular}\right|,\mbox{ }
D_2=\left|
\begin{tabular}{ll}
$m_3$ & $n_{3}$\\
$m_1$ & $n_{1}$
\end{tabular}\right|,\mbox{ and }
D_3=\left|
\begin{tabular}{ll}
$m_1$ & $n_{1}$\\
$m_2$ & $n_{2}$
\end{tabular}\right|.$$
\end{definition}
%%%%%%%%%%%%%%%%%%%%%%%%%%%%%%%%%%%%%%%%%%%%%%%%%%%%%%%%%%%%%%%
Using Lemma 3.16 and Lemma 3.14, we obtain the following
corollary.
\begin{corollary} In terms of iterated integrals, the logarithm of the new symbol is
$$Log[f_1,f_2,f_3]^{\gamma_{jl}}_{C_0,P_{jl}}
=\int_{[\sigma_{jl},\tau]}\frac{df_1}{f_1}\circ\frac{df_2}{f_2}\circ\frac{df_3}{f_3}.$$
\end{corollary}
%%%%%%%%%%%%%%%%%%%%%%%%%%%%%%%%%%%%%%%%%%%%%%%%%%%%%%%%%%%%%%%%%%
\subsection{Reciprocity law for the logarithm of the new symbol}
\begin{theorem}
$$\sum_{P_{jl}} Log[f_1,f_2,f_3]^{\gamma_{jl}}_{C_0,P_{jl}}=(2\pi i)^3(M+N),$$
where
$$\begin{tabular}{lll}
$M=$&$\sum_{j_1<j_2}
(n_{1j_1}D_1(j_2)-n_{3j_1}D_3(j_2))L_{j_1}L_{j_2}$\\
\\
&$+\frac{1}{2}\sum_{j_2=1}^N
(n_{1j_2}D_1(j_2)-n_{3j_2}D_3(j_2))L_{j_2}(L_{j_2}-1),$\\
\\
$N=$&$(2\pi i)^{-2} m_1\left(\int_{\a_j}\frac{dh_3}{h_3}\int_{\b_j}\frac{dh_2}{h_2}
-\int_{\b_j}\frac{dh_3}{h_3}\int_{\a_j}\frac{dh_2}{h_2}\right)+$\\
\\
&$+(2\pi i)^{-2} m_3\left(\int_{\a_j}\frac{dh_1}{h_1}\int_{\b_j}\frac{dh_2}{h_2}
-\int_{\b_j}\frac{dh_1}{h_1}\int_{\a_j}\frac{dh_2}{h_2}\right),$
\end{tabular}
$$
\end{theorem}
\proof Form Theorem 2.31,
we have that the sum of the Log-symbols is equal to the sum of the
extra tems $N$ and the sum of the commutator terms $M$. The sum
of the commutator terms is given by Corollary 3.27.
And the sum of the extra terms is given by Theorem 3.30. This proves the theorem.
%%%%%%%%%%%%%%%%%%%%%%%%%%%%%%%%%%%%%%%%%%%%%%%%%%%%%%%%%%%%%%%%%%%
%%%%%%%%%%%%%%%%%%%%%%%%%%%%%%%%%%%%%%%%%%%%%%%%%%%%%%%%%%%%%%%%%%%
%%%%%%%%%%%%%%%%%%%%%%%%%%%%%%%%%%%%%%%%%%%%%%%%%%%%%%%%%%%%%%%%%%%
\section{Logarithmic version of the Parshin symbol}
\subsection{Integration over a torus revisited}
We are going to use the differential equation from Definition 3.18
and the equivalence from Definition 3.20.
%%%%%%%%%%%%%%%%%%%%%%%%%%%%%%%%%%%%%%
\begin{proposition}
The sum of the coefficients of the monomials equivalent to $B$ in $F_{[\sigma_{jl},\tau]}$ is
$$-4(m_1m_2n_3+m_1n_2m_3+n_1m_2m_3
-n_1n_2m_3-n_1m_2n_3-m_1n_2n_3)\frac{(2\pi i)^3}{2}.$$
\end{proposition}
%%%%%%%%%%%%%%%%%%%%%%%%%%%%%%%%%%%%
\proof Using Definition 3.20 and Lemma 3.16, we have that the term
$A_{1,0}A_{2,0}A_{3,j}$ contributes $-m_1m_2n_3\frac{(2\pi
i)^3}{2}$. Similarly, the term $A_{1,j}A_{2,j}A_{3,0}$ contributes
$+n_1n_2m_3\frac{(2\pi i)^3}{2}$. When we consider all possible
permutations of the indexes, we obtain the above proposition.

%%%%%%%%%%%%%%%%%%%%%%%%%%%%%%%%%%%%%%%%%%%%%%%%%%%%%%%%%%%
\begin{remark}
When we divide the coefficient in the above lemma by
$-4(2\pi i)^2$ and then exponentiate,
we obtain the sign in the Parshin symbol.
\end{remark}
%%%%%%%%%%%%%%%%%%%%%%%%%%%%%%%%%%%%%%%%%%%%%%%%%%%%%%%%%%%%%%%%%%%%
\subsection{Logarithmic symbol}
Let $Q$ be a base point on $C_1$.
Consider the loop $\sigma_{jl}$ around the point $P_{jl}$.

%%%%%%%%%%%%%%%%%%%%%%%%%%%%%%%%%%%%%%%%%%%%%%%%%%%%%%%%%%%
\begin{definition}
For $k=1,2,3$ define $m_k=n_{k1}$, $n_k=n_{kj}$
and $g_{kj}=\frac{f_k}{x_1^{m_k}x_j^{n_k}}=\prod_{i\neq 1,j}x_i^{n_{ki}}$

Define a logarithm of the Parshin symbol as
$$Log\{f_1,f_2,f_3\}^{\gamma_{jl}}_{C_1,P_{jl}}=(2\pi i)^2\left(\pi i K
+D_1\int_{\gamma_{jl}}\frac{dg_{1j}}{g_{1j}}
+D_2\int_{\gamma_{jl}}\frac{dg_{2j}}{g_{2j}}
+D_3\int_{\gamma_{jl}}\frac{dg_{3j}}{g_{3j}}
\right),$$
where
$$K_j=m_1m_2n_{3}-m_1n_{2}m_3+n_{1}m_2n_{3}
-n_{1}n_{2}m_3+n_{1}m_2m_3-m_{1}n_{2}n_3,$$
$$D_1=\left|
\begin{tabular}{ll}
$m_2$ & $n_{2}$\\
$m_3$ & $n_{3}$
\end{tabular}\right|,\mbox{ }
D_2=\left|
\begin{tabular}{ll}
$m_3$ & $n_{3}$\\
$m_1$ & $n_{1}$
\end{tabular}\right|,\mbox{ and }
D_3=\left|
\begin{tabular}{ll}
$m_1$ & $n_{1}$\\
$m_2$ & $n_{2}$
\end{tabular}\right|.$$
\end{definition}
%%%%%%%%%%%%%%%%%%%%%%%%%%%%%%%%%%%%%%%%%%%%%%%%%%%%%%%%%%%%
\begin{theorem}The coefficient of $-\frac{1}{4}B$
for the integral $F$ along the loop, based at $P_{jl}$, is
$Log\{f_1,f_2,f_3\}^Q_{C_1,P_{jl}}$.
\end{theorem}
%%%%%%%%%%%%%%%%%%%%%%%%%%%%%%%%%%%%%%%%%%%%%%%%%%%%%%%%%%%
\proof The portion $K_j$ from the symbol was considered in
subsection 4.1. The remaining portion is anti-symmetrization of
the  indexes of Lemma 3.15. That gives precisely the sum of the
coefficients equivalent to $-\frac{1}{4}B$.
%%%%%%%%%%%%%%%%%%%%%%%%%%%%%%%%%%%%%%%%%%%%%%%%%%%%%%%%%%%%%
\begin{corollary}
The Parshin symbol of the functions $f_1,f_2,f_3$ at $C_1,P_{jl}$ is
$$\{f_1,f_2,f_3\}_{C_1,P_{jl}}
=\left(f_1^{D_1}f_2^{D_2}f_3^{D_3}\right)(Q)
\cdot exp\left((2\pi i)^{-2} \left(Log\{f_1,f_2,f_3\}_{C_1,P_{jl}}
\right)\right).$$
\end{corollary}
\proof It follows by direct computation.
%%%%%%%%%%%%%%%%%%%%%%%%%%%%%%%%%%%%%%%%%%%%%%%%%%%%%%%%%%%%%%
\subsection{Vanishing of the commutator terms}
%%%%%%%%%%%%%%%%%%%%%%%%%%%%%%%%%%%%%%%%%%%%%%%%%%%%%%%%%%%%%
\begin{definition}
Let $$N_{ijk}=(2\pi i)^{-3}\int_{[[\alpha,\beta],\sigma]}
\frac{dx_{i}}{x_{i}}\circ\frac{dx_{j}}{x_{j}}\circ\frac{dx_{k}}{x_{k}}$$
be the integer from the above lemma.
\end{definition}
%%%%%%%%%%%%%%%%%%%%%%%%%%%%%%%%%%%%%%
\begin{lemma} For $i\neq j\neq k\neq i$ we have
$$N_{ijk}+N_{jki}+N_{kij}=0.$$
\end{lemma}
%%%%%%%%%%%%%%%%%%%%%%%%%%%%%%%%%%%%%%%
\proof It follows by direct computation from the formula for $N_{ijk}$
from Lemma 3.18(c).
%%%%%%%%%%%%%%%%%%%%%%%%%%%%%%%%%%%%%%%%%%%%%%%%%%%%%%%%%%
For $i\neq j$, the coefficient of $A_{1i}A_{2i}A_{3j}$
is $n_{1i}n_{2i}n_{3j}N_{iij}$. Note that
$$A_{1i}A_{2i}A_{3j}\sim -A_{1i}A_{3j}A_{2i}.$$
The coefficient of $A_{1i}A_{2i}A_{3j}$
is $n_{1i}n_{2i}n_{3j}N_{iji}$.
And the coefficient of $A_{3j}A_{1i}A_{2i}$
is $n_{1i}n_{2i}n_{3j}N_{jii}$.
Note that
$$A_{3j}A_{1i}A_{2i}\sim A_{1i}A_{3j}A_{2i}.$$
%%%%%%%%%%%%%%%%%%%%%%%%%%%%%%%%%%%%%%%%%%%%%%%%%%%%
\begin{lemma} For $i\neq j$, we have
$$N_{iij}-N_{iji}+N_{jii}
=4\left(\int_{\a}\frac{dx_{j}}{x_{j}}\int_{\b}\frac{dx_{i}}{x_{i}}
-\int_{\b}\frac{dx_{j}}{x_{j}}\int_{\a}\frac{dx_{i}}{x_{i}}\right)
\int_{\sigma}\frac{dx_{i}}{x_{i}}.$$
\end{lemma}
%%%%%%%%%%%%%%%%%%%%%%%%%%%%%%%%%%%%%%%%%%%%%%%%%%%%%%%%%%%%%%%
\proof
For each of the three summands, we have
$$\begin{tabular}{llll}
$N_{iij}$&$=$&$\int_{[[\alpha,\beta],\sigma]}
\frac{dx_{i}}{x_{i}}\circ\frac{dx_{i}}{x_{i}}\circ\frac{dx_{j}}{x_{j}}=$\\
\\
&$=$&$\left(\int_{\a}\frac{dx_{i}}{x_{i}}\int_{\b}\frac{dx_{i}}{x_{i}}
-\int_{\b}\frac{dx_{i}}{x_{i}}\int_{\a}\frac{dx_{i}}{x_{i}}\right)
\int_{\sigma}\frac{dx_{j}}{x_{j}}+$\\
\\
&&$+\left(\int_{\a}\frac{dx_{j}}{x_{j}}\int_{\b}\frac{dx_{i}}{x_{i}}
-\int_{\b}\frac{dx_{j}}{x_{j}}\int_{\a}\frac{dx_{i}}{x_{i}}\right)
\int_{\sigma}\frac{dx_{i}}{x_{i}}=$\\
\\
&$=$&$\left(\int_{\a}\frac{dx_{j}}{x_{j}}\int_{\b}\frac{dx_{i}}{x_{i}}
-\int_{\b}\frac{dx_{j}}{x_{j}}\int_{\a}\frac{dx_{i}}{x_{i}}\right)
\int_{\sigma}\frac{dx_{i}}{x_{i}}$
\end{tabular}$$

$$\begin{tabular}{llll}
$-N_{iji}$&$=$&$-\int_{[[\alpha,\beta],\sigma]}
\frac{dx_{i}}{x_{i}}\circ\frac{dx_{j}}{x_{j}}\circ\frac{dx_{i}}{x_{i}}=$\\
\\
&$=$&$-\left(\int_{\a}\frac{dx_{i}}{x_{i}}\int_{\b}\frac{dx_{j}}{x_{j}}
-\int_{\b}\frac{dx_{i}}{x_{i}}\int_{\a}\frac{dx_{j}}{x_{j}}\right)
\int_{\sigma}\frac{dx_{i}}{x_{i}}-$\\
\\
&&$-\left(\int_{\a}\frac{dx_{i}}{x_{i}}\int_{\b}\frac{dx_{j}}{x_{j}}
-\int_{\b}\frac{dx_{i}}{x_{i}}\int_{\a}\frac{dx_{j}}{x_{j}}\right)
\int_{\sigma}\frac{dx_{i}}{x_{i}}=$\\
\\
&$=$&$-2\left(\int_{\a}\frac{dx_{i}}{x_{i}}\int_{\b}\frac{dx_{j}}{x_{j}}
-\int_{\b}\frac{dx_{i}}{x_{i}}\int_{\a}\frac{dx_{j}}{x_{j}}\right)
\int_{\sigma}\frac{dx_{i}}{x_{i}}$
\end{tabular}$$

$$\begin{tabular}{llll}
$N_{jii}$&$=$&$\int_{[[\alpha,\beta],\sigma]}
\frac{dx_{i}}{x_{i}}\circ\frac{dx_{i}}{x_{i}}\circ\frac{dx_{j}}{x_{j}}=$\\
\\
&$=$&$\left(\int_{\a}\frac{dx_{ij}}{x_{j}}\int_{\b}\frac{dx_{i}}{x_{i_1}}
-\int_{\b}\frac{dx_{j}}{x_{j}}\int_{\a}\frac{dx_{i}}{x_{i}}\right)
\int_{\sigma}\frac{dx_{i}}{x_{i}}=$\\
\\
&&$+\left(\int_{\a}\frac{dx_{i}}{x_{i}}\int_{\b}\frac{dx_{i}}{x_{i}}
-\int_{\b}\frac{dx_{i}}{x_{i}}\int_{\a}\frac{dx_{i}}{x_{i}}\right)
\int_{\sigma}\frac{dx_{j}}{x_{j}}+$\\
\\
&$=$&$\left(\int_{\a}\frac{dx_{j}}{x_{j}}\int_{\b}\frac{dx_{i}}{x_{i}}
-\int_{\b}\frac{dx_{j}}{x_{j}}\int_{\a}\frac{dx_{i}}{x_{i}}\right)
\int_{\sigma}\frac{dx_{i}}{x_{i}}$
\end{tabular}$$

Thus,
$$N_{iij}-N_{iji}+N_{jii}
=4\int_{\sigma}\frac{dx_{i}}{x_{i}}
\left(\int_{\a}\frac{dx_{j}}{x_{j}}\int_{\b}\frac{dx_{i}}{x_{i}}
-\int_{\b}\frac{dx_{j}}{x_{j}}\int_{\a}\frac{dx_{i}}{x_{i}}\right)
.$$

%%%%%%%%%%%%%%%%%%%%%%%%%%%%%%%%%%%%%%%%%%%%%%%%%%%%%%%%%%%%%

Obviously, when three of the indexes coincide then
$N_{iii}=0$.

%%%%%%%%%%%%%%%%%%%%%%%%%%%%%%%%%%%%%%%%%%%%%%%%%%%%%%%%%%%%%%%%
\begin{proposition}
The sum of the coefficients equivalent to $B$ in
$F_{[[\a_i^u,\b_i^u]\tau]}$ is zero.
\end{proposition}
\proof From Lemma 3.18, we obtain that the only contributions to
the coefficient of $B$ come from $A_{1i}A_{2j}A_{3k}$, when only
two of the indexes $i,j,k$ coincide. From Lemma 3.19, it follows
that the contributions are integer combinations of integrals of
the type
$$\int_{\sigma}\frac{dx_{i}}{x_{i}}
\left(\int_{\a_i^u}\frac{dx_{j}}{x_{j}}\int_{\b_i^u}\frac{dx_{i}}{x_{i}}
-\int_{\b_i^u}\frac{dx_{j}}{x_{j}}\int_{\a_i^u}\frac{dx_{i}}{x_{i}}\right).
$$
The integral $$\int_{\sigma}\frac{dx_{i}}{x_{i}}$$ is not zero
only when $i=1$. However, for $i=1$ the variable $x_1$ becomes the
constant $u$ along the loops $\a_i^u$ and $\b_i^u$, because
$\a_i^u$ and $\b_i^u$ are defined on the curve
$$Tb^u=\{Y\in Tb|x_1(Y)=u\}.$$ Therefore,
$$\int_{\a_i^u}\frac{dx_{1}}{x_{1}}=0$$
and
$$\int_{\b_i^u}\frac{dx_{1}}{x_{1}}=0.$$
Thus,
$$\int_{\sigma}\frac{dx_{i}}{x_{i}}
\left(\int_{\a_i^u}\frac{dx_{j}}{x_{j}}\int_{\b_i^u}\frac{dx_{i}}{x_{i}}
-\int_{\b_i^u}\frac{dx_{j}}{x_{j}}\int_{\a_i^u}\frac{dx_{i}}{x_{i}}\right)=0.
$$
%%%%%%%%%%%%%%%%%%%%%%%%%%%%%%%%%%%%%%%%%%%%%%%%%%%%%%%%%%%%%
\subsection{Vanishing of the extra terms}
The sum of the extra terms is the same as the coefficients equivalent to $B$ in
$F_{\pi_i[\mu_{i+1},\tau]\pi_i^{-1}}-F_{[\mu_{i+1},\tau]}$.
%%%%%%%%%%%%%%%%%%%%%%%%%%%%%%%%%%%%%%%%
\begin{proposition}
The coefficient of $\frac{1}{4}B$ in
$$\sum_i F_{\pi_i[\mu_{i+1},\tau]}+F_{[\mu_{i+1},\tau]\pi_i^{-1}}$$ is
zero.
\end{proposition}
%%%%%%%%%%%%%%%%%%%%%%%%%%%%%%%%%%%%%%%%%%%%%%
We need two lemmas in order to compute the coefficient of $A_{1i_1}A_{2i_2}A_{3i_3}$
%%%%%%%%%%%%%%%%%%%%%%%%%%%%%%%%%%%%%%%%%%%%%%%%%%
\begin{lemma}
The coefficient of $A_{1i_1}A_{2i_2}A_{3i_3}$
in $F_{\pi_i[\mu_{i+1},\tau]\pi_i^{-1}}$ is
$$(2\pi i)^2 n_{i_1}n_{2i_2}n_{3i_3}
\left(\int_{\pi_i}\frac{dx_{i_1}}{x_{i_1}}
\int_{[\mu_{i+1},\tau]}
\frac{dx_{i_2}}{x_{i_2}}\circ\frac{dx_{i_3}}{x_{i_3}}
-\int_{\pi_i}\frac{dx_{i_3}}{x_{i_3}}
\int_{[\mu_{i+1},\tau]}
\frac{dx_{i_1}}{x_{i_1}}\circ\frac{dx_{i_2}}{x_{i_2}}
\right).$$
\end{lemma}
%%%%%%%%%%%%%%%%%%%%%%%%%%%%%%%%%%%%%%
\proof Cut the loop $\pi_i[\mu_{i+1},\tau]\pi_i^{-1}$ into $3$ loops
$\pi_i$, $[\mu_{i+1},\tau]$ and $\pi_i^{-1}$.
Then use Corollary 1.3 for composition of paths.
%%%%%%%%%%%%%%%%%%%%%%%%%%%%%%%%%%%%%%%%%%%%%%
\begin{lemma} We have the following relation
$$\mbox{(a) }\int_{[\mu_{i+1},\tau]}
\frac{dx_{i_2}}{x_{i_2}}\circ\frac{dx_{i_1}}{x_{i_1}}
=-\int_{[\mu_{i+1},\tau]}
\frac{dx_{i_1}}{x_{i_1}}\circ\frac{dx_{i_2}}{x_{i_2}}
.$$
$$\mbox{(b) }\int_{[\mu_{i+1},\tau]}
\frac{dx_{i_1}}{x_{i_1}}\circ\frac{dx_{i_1}}{x_{i_1}}
=0.$$
\end{lemma}
%%%%%%%%%%%%%%%%%%%%%%%%%%%%%%%%%%%%%%%%%%%%%%
\proof By Lemma 3.18  part (b), we have that
$$\int_{[\mu_{i+1},\tau]}
\frac{dx_{i_1}}{x_{i_1}}\frac{dx_{i_2}}{x_{i_2}}
=\int_{\mu_{i+1}}\frac{dx_{i_1}}{x_{i_1}}
\int_{\tau}\frac{dx_{i_2}}{x_{i_2}}
-\int_{\tau}\frac{dx_{i_1}}{x_{i_1}}
\int_{\mu_{i+1}}\frac{dx_{i_2}}{x_{i_2}}.$$ The right hand side is
anti-symmetric on the indexes $i_1$ and $i_2$. This proves part
(a). For part (b), use the above equality when $i_1=i_2$.
%%%%%%%%%%%%%%%%%%%%%%%%%%%%%%%%%%%%%%%%%%%%%%%%%%%%%%%%%%%%
\proof (Proposition 4.7) Let $i_1\neq i_2\neq i_3\neq i_1$.
By Lemmas 4.8 and 4.9, we have that the coefficient of
$A_{3i_3}A_{2i_2}A_{1i_1}$ along the path $\pi_i[\mu_{i+1},\tau]\pi_i^{-1}$ is
$$\begin{tabular}{ll}
$(2\pi i)^2 n_{1i_1}n_{2i_2}n_{3i_3}
\left(
\int_{\pi_i}\frac{dx_{i_3}}{x_{i_3}}
\int_{[\mu_{i+1},\tau]}
\frac{dx_{i_2}}{x_{i_2}}\circ\frac{dx_{i_1}}{x_{i_1}}
-\int_{\pi_i}\frac{dx_{i_1}}{x_{i_1}}
\int_{[\mu_{i+1},\tau]}
\frac{dx_{i_3}}{x_{i_3}}\circ\frac{dx_{i_2}}{x_{i_2}}
\right)=$\\
\\
$(2\pi i)^2 n_{1i_1}n_{2i_2}n_{3i_3}
\left(
-\int_{\pi_i}\frac{dx_{i_3}}{x_{i_3}}
\int_{[\mu_{i+1},\tau]}
\frac{dx_{i_1}}{x_{i_1}}\circ\frac{dx_{i_2}}{x_{i_2}}
+\int_{\pi_i}\frac{dx_{i_1}}{x_{i_1}}
\int_{[\mu_{i+1},\tau]}
\frac{dx_{i_2}}{x_{i_2}}\circ\frac{dx_{i_3}}{x_{i_3}}
\right)=$\\
\\
$(2\pi i)^2 n_{1i_1}n_{2i_2}n_{3i_3}
\left(
+\int_{\pi_i}\frac{dx_{i_1}}{x_{i_1}}
\int_{[\mu_{i+1},\tau]}
\frac{dx_{i_2}}{x_{i_2}}\circ\frac{dx_{i_3}}{x_{i_3}}
-\int_{\pi_i}\frac{dx_{i_3}}{x_{i_3}}
\int_{[\mu_{i+1},\tau]}
\frac{dx_{i_1}}{x_{i_1}}\circ\frac{dx_{i_2}}{x_{i_2}}
\right).$
\end{tabular}$$
Note that this is the same as the coefficient of
$A_{1i_1}A_{2i_2}A_{3i_3}$. However, $A_{3i_3}A_{2i_2}A_{1i_1}\sim
-A_{1i_1}A_{2i_2}A_{3i_3}$. Therefore the two coefficients cancel
in the equivalence.

By Lemmas 4.8 and 4.9(b), when three of the indexes coincide, that
is $i_1=i_2=i_3$, the corresponding of the extra terms are zero.

It remain to examine what happens when two of the indexes
coincide. Compare the coefficients of $A_{1i_1}A_{2i_1}A_{3i_2}$,
$A_{1i_1}A_{3i_2}A_{2i_1}$ and $A_{3i_2}A_{1i_1}A_{2i_1}$. By
Lemma 4.8, the coefficient of $A_{1i_1}A_{2i_1}A_{3i_2}$ in
$F_{\pi_i[\mu_{i+1},\tau]\pi_i^{-1}}$ is
$$\begin{tabular}{ll}
$(2\pi i)^2 n_{1i_1}n_{2i_1}n_{3i_2}
\left(
\int_{\pi_i}\frac{dx_{i_1}}{x_{i_1}}
\int_{[\mu_{i+1},\tau]}
\frac{dx_{i_1}}{x_{i_1}}\circ\frac{dx_{i_2}}{x_{i_2}}
-\int_{\pi_i}\frac{dx_{i_2}}{x_{i_2}}
\int_{[\mu_{i+1},\tau]}
\frac{dx_{i_1}}{x_{i_1}}\circ\frac{dx_{i_1}}{x_{i_1}}
\right)=$\\
\\
$(2\pi i)^2 n_{1i_1}n_{2i_1}n_{3i_2}
\left(
\int_{\pi_i}\frac{dx_{i_1}}{x_{i_1}}
\int_{[\mu_{i+1},\tau]}
\frac{dx_{i_1}}{x_{i_1}}\circ\frac{dx_{i_2}}{x_{i_2}}
\right).$
\end{tabular}$$
The coefficient of $A_{1i_1}A_{3i_2}A_{2i_1}$
in $F_{\pi_i[\mu_{i+1},\tau]\pi_i^{-1}}$ is
$$\begin{tabular}{ll}
$(2\pi i)^2 n_{1i_1}n_{2i_1}n_{3i_2}
\left(
\int_{\pi_i}\frac{dx_{i_1}}{x_{i_1}}
\int_{[\mu_{i+1},\tau]}
\frac{dx_{i_2}}{x_{i_2}}\circ\frac{dx_{i_1}}{x_{i_1}}
-\int_{\pi_i}\frac{dx_{i_1}}{x_{i_1}}
\int_{[\mu_{i+1},\tau]}
\frac{dx_{i_1}}{x_{i_1}}\circ\frac{dx_{i_2}}{x_{i_2}}
\right)=$\\
\\
$-2(2\pi i)^2 n_{1i_1}n_{2i_1}n_{3i_2}
\left(
\int_{\pi_i}\frac{dx_{i_1}}{x_{i_1}}
\int_{[\mu_{i+1},\tau]}
\frac{dx_{i_1}}{x_{i_1}}\circ\frac{dx_{i_2}}{x_{i_2}}
\right).$
\end{tabular}$$
And finally, the coefficient of $A_{3i_2}A_{1i_1}A_{2i_1}$
in $F_{\pi_i[\mu_{i+1},\tau]\pi_i^{-1}}$ is
$$\begin{tabular}{ll}
$(2\pi i)^2 n_{1i_1}n_{2i_1}n_{3i_2}
\left(
\int_{\pi_i}\frac{dx_{i_2}}{x_{i_2}}
\int_{[\mu_{i+1},\tau]}
\frac{dx_{i_1}}{x_{i_1}}\circ\frac{dx_{i_1}}{x_{i_1}}
-\int_{\pi_i}\frac{dx_{i_1}}{x_{i_1}}
\int_{[\mu_{i+1},\tau]}
\frac{dx_{i_2}}{x_{i_2}}\circ\frac{dx_{i_1}}{x_{i_1}}
\right)=$\\
\\
$(2\pi i)^2 n_{1i_1}n_{2i_1}n_{3i_2}
\left(
-\int_{\pi_i}\frac{dx_{i_1}}{x_{i_1}}
\int_{[\mu_{i+1},\tau]}
\frac{dx_{i_2}}{x_{i_2}}\circ\frac{dx_{i_1}}{x_{i_1}}
\right)=$\\
\\
$(2\pi i)^2 n_{1i_1}n_{2i_1}n_{3i_2}
\left(
\int_{\pi_i}\frac{dx_{i_1}}{x_{i_1}}
\int_{[\mu_{i+1},\tau]}
\frac{dx_{i_1}}{x_{i_1}}\circ\frac{dx_{i_2}}{x_{i_2}}
\right).$
\end{tabular}$$

The only case when the above integrals are not zero is when
$i_2=0$. Set $i_1=j$. Note that
$$\int_{[\mu_{i+1},\tau]}
\frac{dx_{j}}{x_{j}}\circ\frac{dx_{0}}{x_{0}}=(2\pi i)^2$$
for $\mu_{i+1}=\tilde{\sigma}_{jl}$.
Then
$\int_{\pi_i}\frac{dx_{j}}{x_{j}}=2\pi il.$

However, the sum of the extra terms coming from
$A_{1i_1}A_{2i_1}A_{3i_2}$, $A_{1i_1}A_{2i_2}A_{3i_1}$ and
$A_{1i_2}A_{2i_1}A_{3i_1}$ is zero. Thus, there is no contribution
from the extra terms in the logarithm of the Parshin symbol.
%%%%%%%%%%%%%%%%%%%%%%%%%%%%%%%%%%%%%%%%%%%%%%%%%%%%%%%%%%%
\subsection{Reciprocity law}

\begin{theorem}
For the logarithm of the Parshin symbol we have the following
reciprocity
$$\sum_{j,l}Log\{f_1,f_2,f_3\}^{\gamma_{jl}}_{C_1,P_{jl}}=0.$$
\end{theorem}
\proof  Form Theorem 2.31,
we have that the sum of the Log-symbols applied to an iteration
$$\frac{dx_i}{x_i}\circ\frac{dx_j}{x_j}\circ\frac{dx_k}{x_k}$$
is equal to the sum of the extra terms and the sum of the
commutator terms. We consider the sum of the coefficients
equivalent to $-1/4 B$. The sum of the commutator terms is zero by
Proposition 4.9. And the sum of the extra terms is zero by
Proposition 4.10. This proves the theorem.
%%%%%%%%%%%%%%%%%%%%%%%%%%%%%%%%%%%%%%%%%%%%%%%%%%%%%%%%%%%%%%%%%%%%%%%%%%%%%%%%%
%\subsection{Independence of the choice of local coordinates}

%%%%%%%%%%%%%%%%%%%%%%%%%%%%%%%%%%%%%%%%%%%%%%%%%%%%%%%%%%%%%%%%%%%%%%%%%%%%%%%%%
\section{Refinement of the Parshin symbol}
\subsection{Logarithmic version of a refinement of the Parshin symbol}
%%%%%%%%%%%%%%%%%%%%%%%%%%%%%%%
\begin{definition}
Define a logarithm of the refinement of the Parshin symbol as the
difference
$$Log(f_1,f_2,f_3)^{\gamma_{jl}}_{C_0,P_{jl}}=
Log[f_1,f_2,f_3]^{\gamma_{jl}}_{C_0,P_{jl}}
-Log\{f_1,f_2,f_3\}^{\gamma_{jl}}_{C_0,P_{jl}}.$$
\end{definition}
%%%%%%%%%%%%%%%%%%%%%%%%%%%%%%%%%%%%%%%%
\begin{theorem}
$$Log(f_1,f_2,f_3)^{\gamma_{jl}}_{C_0,P_{jl}}
=(2\pi i)^2\left(\pi i(m_2n_1n_3-n_2m_1m_3)
+m_2n_3\int_{\gamma_{jl}}\frac{dg_{1j}}{g_{1j}}
-m_2n_1\int_{\gamma_{jl}}\frac{dg_{3j}}{g_{3j}},
\right)$$
\end{theorem}
\proof It follows directly from Definitions 3.28 and 4.3.
%%%%%%%%%%%%%%%%%%%%%%%%%%%%%%%%%%%%%%%%%%
As a direct consequence we obtain a refinement of the logarithm of
the Parshin symbol.
\begin{corollary}
$$Log(f_1,f_2,f_3)^{\gamma_{jl}}_{C_0,P_{jl}}
+Log(f_2,f_3,f_1)^{\gamma_{jl}}_{C_0,P_{jl}}
+Log(f_3,f_1,f_2)^{\gamma_{jl}}_{C_0,P_{jl}}=
Log\{f_1,f_2,f_3\}^{\gamma_{jl}}_{C_0,P_{jl}}.$$
\end{corollary}
%%%%%%%%%%%%%%%%%%%%%%%%%%%%%%%%%%%%%%%%
\subsection{Logarithmic reciprocity law}
\begin{theorem} We have the following reciprocity law
$$\sum_{jl}Log(f_1,f_2,f_3)^{\gamma_{jl}}_{C_0,P_{jl}}=(2\pi i)^3(M+N),$$
where
$$
\begin{tabular}{ll}
$M=$
&$(2\pi i)^3\sum_{j_1<j_2}
(n_{1j_1}D_1(j_2)-n_{3j_1}D_3(j_2))L_{j_1}L_{j_2}$\\
\\
&$+(2\pi i)^3\sum_{j_2=1}^N
(n_{1j_2}D_1(j_2)-n_{3j_2}D_3(j_2))\frac{1}{2}L_{j_2}(L_{j_2}-1)$\\
\\
$N=$&$(2\pi i)^{-2} m_1\left(\int_{\a_j}\frac{dh_3}{h_3}\int_{\b_j}\frac{dh_2}{h_2}
-\int_{\b_j}\frac{dh_3}{h_3}\int_{\a_j}\frac{dh_2}{h_2}\right)+$\\
\\
&$ +(2\pi i)^{-2} m_3\left(\int_{\a_j}\frac{dh_1}{h_1}\int_{\b_j}\frac{dh_2}{h_2}
-\int_{\b_j}\frac{dh_1}{h_1}\int_{\a_j}\frac{dh_2}{h_2}\right).$
\end{tabular}
$$
\end{theorem}
\proof It follows from the reciprocity laws for
$Log[f_1,f_2,f_3]^{\gamma_{jl}}_{C_0,P_{jl}}$ and
$Log\{f_1,f_2,f_3\}^{\gamma_{jl}}_{C_0,P_{jl}}$, stated in
Theorems 3.34 and 4.13, respectively.
%%%%%%%%%%%%%%%%%%%%%%%%%%%%%%%%%%%%%%%%%%%%%%%
\subsection{Refinement of the Parshin symbol and a reciprocity law}
\begin{definition}
Define a refinement of the Parshin symbol as
$$(f_1,f_2,f_3)^{x_0}_{C_0,P_{jl}}=
\frac{g_{1j}(Q)^{m_2n_3}}{g_{3j}(Q)^{m_2n_1}}
exp\left((2\pi i)^{-2}Log(f_1,f_2,f_3)^{\gamma_{jl}}_{C_0,P_{jl}}\right).$$
\end{definition}
%%%%%%%%%%%%%%%%%%%%%%%%%%%%%%%%%%%%%%%%%%%%%%%%%%%%%%%%
\begin{theorem}
Define a refinement of the Parshin symbol as
$$(f_1,f_2,f_3)^{x_0}_{C_0,P_{jl}}=
(-1)^{m_2n_1n_3-n_2m_1m_3}
\frac{g_{1j}(P_{jl})^{m_2n_3}}{g_{3j}(P_{jl})^{m_2n_1}}.$$
\end{theorem}
%%%%%%%%%%%%%%%%%%%%%%%%%%%%%%%%%%%%%%%%%
\begin{remark}
For the local coordinates at $P_{jl}$ we use the variables $x_0$
and $x_j$. The symbol is invariant of the choices of $x_j$.
However, it is dependent on $x_0$. Let us explain to what extend
the refinement of the Parshin symbol depends on the choice of
$x_0$. The variable $x_0$ vanishes along the curve $C_0$. And we
use the same $x_0$ for each of the points $P_{jl}$, where we
compute the symbol. If we have a different choice of $x_0$ then we
obtain a different symbol. But the new choice of $x_0$ is the same
for each point $P_{jl}$.
\end{remark}
%%%%%%%%%%%%%%%%%%%%%%%%%%%%%%%%%%%%%%%%%%%%%%%%%%%%%%%
\begin{theorem} Using the above definition, we have
$$\prod_{jl}(f_1,f_2,f_3)^{x_0}_{C_0,P_{jl}}=1.$$
\end{theorem}
%%%%%%%%%%%%%%%%%%%%%%%%%%%%%%%%
\proof It follows by exponentiating the reciprocity law from Theorem 4.5.
\subsection{Example}
%%%%%%%%%%%%%%%%%%%%%%%%%%%%%%%%%%%%%%%%%%%%%%%%%%%%%%%%%%%%%%
Let $X=\C P^1\times \C P^1$ with projective coordinates
$(x_0:x_1)\times(y_0,y_1)$. For $n=1,2,3$ let
$$f_n=(x_1-ax_0)^{i_n}(x_1-bx_0)^{j_n}(x_1-cx_0)^{k_n}
\left(\frac{y_1}{y_0}\right)^{l_n},$$
where $i_n+j_n+k_n=0$ for $n=1,2,3$. Let $C$ be the variety
$y_1=0$. Consider the functions $f_n$ in affine coordinate system
by setting $x_0=1$ $y_0=1$. Let
$$f_n=(x_1-a)^{i_n}(x_1-b)^{j_n}(x_1-c)^{k_n}y_1^{l_n}.$$
The points which will have non-zero symbol are $P_a=(a,0)$,
$P_b=(b,0)$ and $P_c=(c,0)$. Let the base point $Q$ be with
coordinates $Q=(0,0)$.

First we compute the symbol $(f_1,f_2,f_3)^{y_1}_{C,P_a}$. Let
$$z=z(x_1)=\frac{x_1-a}{x_1-b}.$$ Then
$$x_1=\frac{bz-a}{z-1}.$$
Also
$$x_1-b=\frac{b-a}{z-1}$$
and
$$x_1-c=\frac{(b-c)z+c-a}{z-1}.$$ Then
$$f_n(z,y_1)=z^{i_n}y_1^{l_n}\left(\frac{b-a}{z-1}\right)^{i_n+j_n}
\left(\frac{(b-c)z+c-a}{z-1}\right)^{k_n}.$$ Then at $P_a$, we
have
$$g_n(z,y_1)=\left(\frac{b-a}{z-1}\right)^{i_n+j_n}
\left(\frac{(b-c)z+c-a}{z-1}\right)^{k_n}.$$ Using that
$i_n+j_n+k_n=0$, we obtain
$$g_n(z,y_1)=\left(\frac{b-a}{z-1}\right)^{-k_n}
\left(\frac{(b-c)z+c-a}{z-1}\right)^{k_n}.$$ At the point $P_a$,
which is with coordinates $(z,y_1)=(0,0)$, we have
$$g_n(P_a)=\left(\frac{a-c}{a-b}\right)^{k_n}.$$
Therefore,
$$g_3(P_a)^{i_1l_2}g_1(P_a)^{-i_3l_2}=
\left(\frac{a-b}{a-c}\right)^{l_2
\left|\begin{tabular}{ll} $i_1$ & $k_1$\\
$i_3$ & $k_3$
\end{tabular}\right|
}
$$
For the sign of the symbol, we have $(-1)^{i_1i_3l_2-l_1l_3i_2}$.
Finally, the symbol at $P_a$ is
$$(f_1,f_2,f_3)_{C,P_a}^{y_1}=(-1)^{i_1i_3l_2-l_1l_3i_2}\left(\frac{a-b}{a-c}\right)^{l_2
\left|\begin{tabular}{ll} $i_1$ & $k_1$\\
$i_3$ & $k_3$
\end{tabular}\right|
}.
$$
For the symbol at the points $P_b$ and $P_c$, we have to permute
cyclicly $a,b,c$ and use the same cyclic permutation for the
powers $i_n,j_n,k_n$. We obtain
$$(f_1,f_2,f_3)_{C,P_b}^{y_1}=(-1)^{j_1j_3l_2-l_1l_3j_2}\left(\frac{b-c}{b-a}\right)^{l_2
\left|\begin{tabular}{ll} $j_1$ & $i_1$\\
$j_3$ & $i_3$
\end{tabular}\right|
}
$$
and
$$(f_1,f_2,f_3)_{C,P_c}^{y_1}=(-1)^{k_1k_3l_2-l_1l_3k_2}\left(\frac{c-a}{c-b}\right)^{l_2
\left|\begin{tabular}{ll} $k_1$ & $j_1$\\
$k_3$ & $j_3$
\end{tabular}\right|
}.
$$
Let $$A=\frac{a-b}{a-c}$$ and
$$B=\frac{b-c}{b-a},$$
which enter in the first two symbols. Then for the fraction in the
last symbol, we have
$$\frac{c-a}{c-b}=-(AB)^{-1}.$$
The product of the three symbols is
$$
\begin{tabular}{ll}
$(f_1,f_2,f_3)_{C,P_a}^{y_1}(f_1,f_2,f_3)_{C,P_b}^{y_1}(f_1,f_2,f_3)_{C,P_c}^{y_1}=$\\
\\
$=(-1)^{(i_1i_3+j_1j_3+k_1k_3)l_2-(i_2+j_2+k_2)l_1l_3}A^{l_2
\left|\begin{tabular}{ll} $i_1$ & $k_1$\\
$i_3$ & $k_3$
\end{tabular}\right|
}B^{l_2
\left|\begin{tabular}{ll} $j_1$ & $i_1$\\
$j_3$ & $i_3$
\end{tabular}\right|
}(-AB)^{-l_2
\left|\begin{tabular}{ll} $k_1$ & $j_1$\\
$k_3$ & $j_3$
\end{tabular}\right|
}.$
\end{tabular}
$$
Combining the power of $A$ and using that $i_n+j_n+k_n=0$, we
obtain
$$\begin{tabular}{ll}$l_2
\left|\begin{tabular}{ll} $i_1$ & $k_1$\\
$i_3$ & $k_3$
\end{tabular}\right|
-l_2
\left|\begin{tabular}{ll} $k_1$ & $j_1$\\
$k_3$ & $j_3$
\end{tabular}\right|
=l_2
\left|\begin{tabular}{ll} $-j_1-k_1$ & $k_1$\\
$-j_3-k_3$ & $k_3$
\end{tabular}\right|
-l_2
\left|\begin{tabular}{ll} $k_1$ & $j_1$\\
$k_3$ & $j_3$
\end{tabular}\right|=$\\
\\
$=l_2
\left|\begin{tabular}{ll} $-j_1$ & $k_1$\\
$-j_3$ & $k_3$
\end{tabular}\right|
-l_2
\left|\begin{tabular}{ll} $k_1$ & $j_1$\\
$k_3$ & $j_3$
\end{tabular}\right|=0.$
\end{tabular}$$
Similarly, the powers of $B$ cancel. For the sign we have
$$(-1)^{(i_1i_3+j_1j_3+k_1k_3)l_2-(i_2+j_2+k_2)l_1l_3}(-1)^{-l_2
\left|\begin{tabular}{ll} $k_1$ & $j_1$\\
$k_3$ & $j_3$
\end{tabular}\right|
}$$ For the power of $(-1)$ modulo $2$, we have
$$\begin{tabular}{ll}
$(i_1i_3+j_1j_3+k_1k_3)l_2-(i_2+j_2+k_2)l_1l_3+l_2(k_1j_3+j_1k_3)=$\\
\\
$l_2(i_1i_3+(j_1+k_1)(j_3+k_3))-(i_2+j_2+k_2)l_1l_3$\\
\\
$l_2(i_1i_3+(-i_1)(-i_3))+0\cdot l_1l_3=$\\
\\
$=0$ mod $2$
\end{tabular}$$
Therefore,
$$(f_1,f_2,f_3)_{C,P_a}^{y_1}(f_1,f_2,f_3)_{C,P_b}^{y_1}(f_1,f_2,f_3)_{C,P_c}^{y_1}=1.$$
%%%%%%%%%%%%%%%%%%%%%%%%%%%%%%%%%%%%%%%%%%%%%%%%%%%%%%%%%%%%%%%%%%%%%%%%%%%%%%%%%
%%%%%%%%%%%%%%%%%%%%%%%%%%%%%%%%%%%%%%%%%%%%%%%%%%%%%%%%%%%%%%%%%%%%%%%%%%%%%%%%%
%%%\begin{thebibliography}
%%%%%%%%%%%%%%%%%%%%%%%%%%%%%%%%%%%%%%%%%%%%%%%%%%%%%%%%%%%
\renewcommand{\em}{\textrm}

\begin{small}

\renewcommand{\refname}{ {\flushleft\normalsize\bf{References}} }
    
\end{small}

\end{document}